\documentclass{amsart}

\begin{document}

\title[Levi-Civita]{Arithmetic Levi-Civita connection}
\bigskip

\def \cO{\mathcal O}
\def \ra{\rightarrow}
\def \bZ{{\mathbb Z}}
\def \cP{\mathcal V}
\def \cH{{\mathcal H}}
\def \cB{{\mathcal B}}
\def \d{\delta}

\newtheorem{THM}{{\!}}[section]
\newtheorem{THMX}{{\!}}
\renewcommand{\theTHMX}{}
\newtheorem{theorem}{Theorem}[section]
\newtheorem{corollary}[theorem]{Corollary}
\newtheorem{lemma}[theorem]{Lemma}
\newtheorem{proposition}[theorem]{Proposition}
\theoremstyle{definition}
\newtheorem{definition}[theorem]{Definition}
\theoremstyle{remark}
\newtheorem{remark}[theorem]{Remark}
\newtheorem{example}[theorem]{\bf Example}
\numberwithin{equation}{section}

\address{Department of Mathematics and Statistics\\University of New Mexico \\ Albuquerque, NM 87131, USA}
\email{buium@math.unm.edu} 
\subjclass[2010]{11E95,20G25, 53B20}
\maketitle

\bigskip

\medskip
\centerline{\bf Alexandru Buium}
\bigskip

\begin{abstract}
This paper is part of a series of papers where an arithmetic analogue of classical differential geometry is being developed. In this arithmetic differential geometry functions are replaced by integer numbers, derivations are replaced by Fermat quotient operators, and connections (respectively curvature) are replaced by certain adelic (respectively global) objects attached to symmetric matrices with integral coefficients. Previous papers were devoted to an arithmetic analogue of the Chern connection. The present paper is devoted to an arithmetic analogue of the Levi-Civita connection. 
\end{abstract}

\section{Introduction}

\subsection{Aim of the paper} 
From a technical viewpoint the present
 paper is devoted to proving the existence and uniqueness
of certain remarkable Frobenius lifts on the various $p$-adic completions of the $\bZ$-group scheme $GL_n$;
each such collection of Frobenius lifts will be   attached to a given symmetric matrix with integer coefficients. We will then consider the problem of defining and computing commutators of these Frobenius lifts as $p$ varies. 

From a conceptual viewpoint,
the above collection of  Frobenius lifts attached to a symmetric integral matrix can be viewed as an arithmetic analogue, for the spectrum of the integers, of the Levi-Civita connection attached to a metric on a manifold; the collection of commutators of these Frobenius lifts can then be viewed as an arithmetic analogue of curvature. As such
this paper can be viewed as  part of a series of papers \cite{curvature1, adel2, adel3, curvature2}
devoted to  developing an arithmetic analogue of {\it classical differential geometry}; this analogue can be referred to as {\it arithmetic differential geometry}. However, for the convenience of the reader,  the present paper is written so as to be entirely self-contained. 

Expressed in a naive  form, the main idea in the papers \cite{curvature1, adel2, adel3, curvature2}
is to  replace functions on smooth manifolds by integer numbers, to replace coordinates by prime numbers,  and to  replace differentiation acting on functions by ``arithmetic derivative" operators acting on numbers.
In this setting the ``arithmetic derivative" of an integer $n\in \bZ$ with respect to a prime $p$ is taken to be equal to the {\it Fermat quotient} $\frac{n-n^p}{p}$. In order to turn   this  idea  into a functioning theory one needs to ``geometrize" it in the same sense in which Lie and Cartan ``geometrized" differential equations. 
We refer to the monographs \cite{book} and \cite{foundations} for a comprehensive introduction to this program. Cf. also \cite{char, duke, crelle, local} for some purely arithmetic applications
of this theory.  This whole line of research is, of course, part of the general, well established, effort to unveil and exploit the  analogies between numbers and functions.

The papers \cite{curvature1, adel2, adel3, curvature2} were mainly concerned with an arithmetic analogue of Chern connection. On the other hand, in the monograph \cite{foundations}, first steps were taken, in  a special case, towards developing an arithmetic analogue of the Levi-Civita connection. The Levi-Civita story in \cite{foundations} has, however, at least two  limitations.
First, if one fixes a prime $p$ then the story in \cite{foundations} 
 only deals, in some sense, with an analogue of metrics of cohomogeneity one: indeed, for $p$ fixed,  there is only one ``arithmetic derivative"
there acting on the coefficients of the metric. 
Second, for varying $p$, the story in \cite{foundations} only deals with ``metrics with constant coefficients". Although, in arithmetic,  such metrics still lead  to non-zero curvature, restricting attention to such metrics is a drastic limitation.

In the present paper we would like to revisit from scratch the arithmetic  Levi-Civita story in \cite{foundations} by putting it in a more general context: in this context both limitations referred to above will disappear.  The first limitation (where $p$ is fixed) will be overcome by considering several ``arithmetic derivatives" corresponding to the  primes dividing $p$ in a  number field. The second limitation (where $p$ varies) will be overcome by constructing {\it algebraizing correspondences} for our arithmetic Levi-Civita connection
(in the same sense in which we constructed algebraizing correspondences for the Chern connection  in \cite{curvature2, foundations}). In this second context the curvature of Levi-Civita connection in dimension $n$ will take values in the ${\mathbb Q}$-{\it algebra of correspondences}
on the function field of $GL_n$ which is an infinite dimensional Lie ${\mathbb Q}$-algebra. 
One can ask if this infinite dimensional Lie algebra can be replaced by a finite dimensional one.
This turns out to be possible in special cases; we will tackle this problem elsewhere. 

\subsection{Organization of the paper} Section 2 contains the  definitions and statements of our main results on Levi-Civita connections. 
Section 3 contains the definitions and statements of our main results on curvature.
Both sections 2 and 3 also include a discussion (and proofs) of consequences of our main results.
Section 4 contains the proofs of our main results. Section 5 is an Appendix devoted to revisiting 
(from a somewhat non-conventional angle) the classical differential geometric setting. This Appendix is not logically necessary for the understanding of the paper. However, since the Appendix provides the main motivation/blueprint for the arithmetic story, {\it the reader is strongly 
encouraged  to read the Appendix before starting in on the body of the paper}.

\subsection{Main results} Our main results will be:

1) An existence and uniqueness theorem for our (adelic) Levi-Civita connections (which can be viewed as an  analogue of the ``Fundamental Theorem of Riemannian Geometry"); cf. Theorems \ref{levi Civita} and \ref{mock}. 

 2) Algebraization  theorems for these connections (allowing to define curvature for varying $p$); cf. Theorems \ref{algebraization}, \ref{algebraization1} and Proposition \ref{galben}.
 
 3)  A series of congruences mod $p$ for our connections and curvature in case $p$ is fixed (leading, in particular, to analogues of the classical symmetries of the Riemann tensor); cf.
 Propositions \ref{voices}, \ref{congruences}, \ref{elisabeth}, and Corollaries
 \ref{seinfeld1}, \ref{seinfeld2}, \ref{jaja}.
 
  4) A series of computations of curvature (in particular non-vanishing theorems for curvature), especially for ``conformal metrics," in dimension $n=2$; cf. Propositions  \ref{fishatnoon}, 
  \ref{maro}, and
   Corollary \ref{musca}.

\subsection{General conventions and notation} 
 For background on schemes and  formal schemes we refer to Chapter 2 in \cite{hartshorne}. For background on local fields and number fields  we refer to Chapter 1 of \cite{LangNumberTheory}. 
 
Unless otherwise stated all rings are commutative and unital. 
When commutativity is not assumed we will talk about {\it not necessarily commutative} rings; in this case
homomorphisms and antihomomorphisms will be unital 
and, to simplify notation and terminology, antihomomorphisms will often be referred to, again,  as {\it homomorphisms}. Also, in this case, we will often use the same letter to denote a ring and its opposite; the context will always indicate the precise meaning of our notation.
By a {\it Lie ring} we understand a Lie $\bZ$-algebra. 
Any (not necessarily commutative) ring can be viewed as a Lie ring with respect to the commutator.
 For any (not necessarily commutative) monoid $M$ we denote by $\bZ M$ the (not necessarily commutative) monoid ring on $M$.
 For any  ring $A$ we denote by ${\mathfrak g}{\mathfrak l}_n(A)$ the (not necessarily commutative) ring  of $n\times n$ matrices with coefficients in $A$ and we denote  by $GL_n(A)$ the group of invertible elements of ${\mathfrak g}{\mathfrak l}_n(A)$.
For any rings $A$ and $B$ and any set theoretic map 
$f:A\ra B$
we still denote by 
$f:{\mathfrak g}{\mathfrak l}_n(A)\ra {\mathfrak g}{\mathfrak l}_n(B)$ the induced map; so for 
 any $n\times n$ matrix $a=(a_{ij})$ with entries $a_{ij}\in A$ we let $f(a)=(f(a_{ij}))$ the $n \times n$ matrix with entries $f(a_{ij})$. For a matrix $a=(a_{ij})\in {\mathfrak g}{\mathfrak l}_n(A)$ we denote by $a^t$ the transpose of $a$; if in addition $p\in \bZ$ is a prime then we denote by $a^{(p)}=(a_{ij}^p)$ the matrix with entries $a_{ij}^p$.
 
All schemes and formal schemes are assumed separated; formal schemes are  assumed Noetherian.

 For $A$ a ring (or $X$ a Noetherian scheme) and a fixed prime $p\in \bZ$ (always assumed given in our context) we always denote by $\widehat{A}$ (respectively $\widehat{X}$)
 the $p$-adic completion of $A$ (respectively $X$). As a rule ring homomorphisms $A\ra B$ and the corresponding morphisms between their spectra $Spec\ B\ra Spec\ A$, or between the formal spectra $Spf\ \widehat{B}\ra Spf\ \widehat{A}$, 
 will  be denoted by the same letters.

\subsection{Acknowledgements} The author is  indebted to Lars Hesselholt and Yuri I. Manin  for inspiring suggestions. The present work was partially supported by the Max-Planck-Institut f\"{u}r Mathematik in Bonn, by the Institut des Hautes \'{E}tudes Scientifiques in Bures sur Yvette,  and by  the Simons Foundation (award 311773).

\section{Connections}
We start by   recalling some basic terminology from \cite{book, foundations}.
 
\subsection{$p$-adic connections}
 Let $A$ be a ring  and let $p\in \bZ$ be an odd prime.
 By a {\it Frobenius lift} on $A$ we understand a ring endomorphism $\phi=\phi^A:A\ra A$ reducing mod $p$ to  the $p$-power 
 Frobenius  $A/pA\ra A/pA$. By a {\it $p$-derivation} 
we understand a map of sets 
$\d=\d^A:A\ra A$ such that $\d 1=0$ and for all $a,b\in A$,

\medskip

1) $\d(a+b)=\d a + \d b+\sum_{k=1}^{p-1} p^{-1}\left(\begin{array}{c} p\\ k\end{array}\right) a^kb^{p-k}$,

\medskip

2) $\d(ab)=a^p\d b + b^p \d a + p (\d a) (\d b)$.

\medskip

\noindent  If $\d:A\ra A$ is a $p$-derivation then the map  $\phi=\phi^A:A\ra A$ defined by
 \begin{equation}
\label{lf}
\phi(a)=a^p+p\d a\end{equation}
 is  a Frobenius lift and we say that $\phi$ and $\d$ are {\it attached} to each other.
Conversely 
if $A$ is $p$-torsion free (i.e., $p$ is a non-zero divisor in $A$) then any Frobenius lift $\phi:A\ra A$ defines a unique $p$-derivation
$\d:A\ra A$ by the formula \ref{lf}.

 Assume $X$ is a Noetherian scheme. A {\it Frobenius lift} on $X$ (respectively on  $\widehat{X}$) will mean a scheme (respectively a formal scheme) endomorphism whose reduction mod $p$ is the $p$-power Frobenius.
 A {\it  $p$-derivation}  on $\widehat{X}$  will mean a map of sheaves of sets $\d^{\widehat{X}}:\cO_{\widehat{X}}\ra \cO_{\widehat{X}}$ which is a $p$-derivation on each open set.
We usually denote by $\phi^{\widehat{X}}:\widehat{X}\ra \widehat{X}$ the attached Frobenius lift.  We sometimes write $\d^X$, $\phi^X$, or even $\d$, $\phi$, instead of $\d^{\widehat{X}}$, $\phi^{\widehat{X}}$.
If $X$ is affine to give a $p$-derivation on $\widehat{X}$ is the same as to give a $p$-derivation on $\cO(\widehat{X})$.
If $\cO_X$ is $p$-torsion free then to give a $p$-derivation on $\widehat{X}$ is the same as to give a Frobenius lift on $\widehat{X}$.

If $X$ and $\d^X$ are as above and $Y\subset X$ is a closed subscheme we say that 
$Y$ is {\it $\d^X$-horizontal} (respectively {\it $\phi^X$-horizontal}) if the ideal defining $\widehat{Y}$ in $\widehat{X}$ is sent into itself by $\d^X$ (respectively by $\phi^X$). If $Y$ is $\d^X$-horizontal we have an induced $p$-derivation $\d^Y$ on $\widehat{Y}$.
If $\cO_Y$ is $p$-torsion free then $Y$ is $\d^X$-horizontal if and only if it is $\phi^X$-horizontal.

\begin{definition}
Let $X$ be a scheme of finite type over a Noetherian ring $\cO$ and let $\d^{\cO}$ be a $p$-derivation on $\cO$. A {\it $p$-adic connection} on $X$ is a $p$-derivation 
$\d=\d^X=\d^{\widehat{X}}$ on
$\widehat{X}$ which extends the $p$-derivation $\d=\d^{\cO}$ on $\cO$.
\end{definition}

As a rule, in this paper, we will only be interested in $p$-adic connections $\d^G$ on smooth group schemes $G$ over $\cO$. 
They should be viewed as arithmetic analogues of  connections in principal bundles  in
the sense of  classical differential geometry.  The case when a subgroup scheme $H\subset G$ is $\d^G$-horizontal should be viewed as an analogue of reduction of the structure group to a subgroup in the classical case.
As noted in \cite{foundations},
asking, as in classical differential geometry,  that our $p$-adic connections be  ``translation invariant" is a condition that is {\it almost  never} satisfied in the arithmetic theory.
But asking from our $p$-adic connections that they  be ``metric" or  ``torsion free" (cf. the definitions below) is reasonable and leads to an interesting theory. We shall follow this path  in what follows.

 \subsection{Metric connections}

Let  $p\in \bZ$ be an odd prime and let $\cO$ be any complete discrete valuation ring with maximal ideal generated by $p$ and perfect residue field. Such a ring possesses a unique Frobenis lift $\phi=\phi^{\cO}$ and hence a unique $p$-derivation $\d=\d^{\cO}$. If such an $\cO$ is given we will say we are in the  {\it local situation}. 

Assume now we are in this local situation.

 Let $x=(x_{ij})$ be an $n\times n$  matrix of indeterminates. We may consider the group scheme
 $$G=GL_n=Spec\ \cO[x,\det(x)^{-1}];$$
 so for the ring of global functions we have 
$\cO(G)=\cO[x,\det(x)^{-1}]$.
According to our terminology  a {\it $p$-adic connection} on $G$ is a $p$-derivation $\d^G$ on $\widehat{G}$ (equivalently  on $\cO(\widehat{G})$) extending the $p$-derivation $\d^{\cO}$.
We also set
$${\mathfrak g}:={\mathfrak g}{\mathfrak l}_n:=Spec\ \cO[x].$$

An example of $p$-adic connection is the {\it trivial $p$-adic connection}, $\d_0^G$,  defined by $\d_0^G x=0$; i.e., $\d_0^G x_{ij}=0$.
Its attached Frobenius lift $\phi_0^{G}$ satisfies $\phi_0^G (x)=x^{(p)}$; i.e., $\phi_0^G(x_{ij})=x_{ij}^p$.

 Consider next a {\it symmetric} matrix $q\in G(\cO)=GL_n(\cO)$, i.e., $q^t=q$. 
 We view $q$ as an arithmetic analogue of a metric.

\begin{definition}
 The {\it Christoffel symbol} (of the first kind) of a
 $p$-adic  connection $\d^G$ on $G$ relative to $q$
  is the matrix
  \begin{equation}
\label{pen}\Gamma:=\d^G x^t \cdot \phi(q)\cdot x^{(p)}\in {\mathfrak g}(\cO(\widehat{G})).\end{equation}
\end{definition}

Our Christoffel symbols can be viewed as analogues of the Christoffel symbols in classical differential geometry; cf. our Appendix and Remark \ref{ozu} below.

 Now to the matrix $q$ one can attach the map of schemes over $\cO$,
$$\cH_q:G\ra G,$$
 defined on the level of rings by the ring homomorphism (still denoted by) 
$$\cH_q:\cO(G)\ra \cO(G),\ \ \cH_q(x):=x^tqx,$$
i.e.,
$$\cH_q(x_{ij})=q_{kl}x_{ki}x_{lj},$$
with the repeated indices $k,l$ summed over. Note that here (and later) we  adopt the Einstein summation notation for indices that are not necessarily appearing both in upper and lower 
positions; no confusion should arise from this.
 We continue to denote by 
$$\cH_q:\widehat{G}\ra \widehat{G}$$ the induced map of formal schemes. The identity matrix in $GL_n(\cO)$ will always be denoted by $1_n$, or simply by $1$ if $n$ is understood from context.

\begin{definition}\label{qui1}
 A $p$-adic  connection $\d^G$ on $G$ with attached  Frobenius lift $\phi^G$ is said to be {\it metric} with respect to $q$ 
  if the following diagram is commutative:
 \begin{equation}
 \label{got}
 \begin{array}{rcl}
 \widehat{G} & \stackrel{\phi^G}{\longrightarrow} & \widehat{G} \\
 \cH_q \downarrow &\ &\downarrow \cH_q\\
 \widehat{G} & \stackrel{\phi_{0}^G}{\longrightarrow} & \widehat{G}\end{array}
\end{equation}
Alternatively, following \cite{adel2, foundations},
 we say that $\phi^G$ is
{\it $\cH_q$-horizontal} with respect to the trivial $p$-adic connection $\phi^G_0$.
\end{definition}

The above condition defining metric connections should be viewed as an arithmetic analogue of the classical concept of {\it metric connection}; cf. \ref{gott} in our Appendix.

\begin{remark}
\label{preozu}
 Explicitly let $\d^G$ be any $p$-adic connection on $G$, 
 and consider the $n\times n$ matrices 
 $$\Delta,\ \ \Lambda,\ \ A,\ \ B,\ \ S$$
 with entries in $\cO(\widehat{G})$ defined by the equalities
$$\begin{array}{rcl}
\d^G x & =: &\Delta,\\
\ & \ & \ \\
\phi^G(x)& =: & x^{(p)}\Lambda=x^{(p)}+p\Delta,\\
\ & \ & \ \\
A & := & x^{(p)t}\cdot \phi(q)\cdot x^{(p)},\\
\ & \ & \ \\
 B & := & (x^tqx)^{(p)}\\
 \ & \ & \ \\
 S & := & (x^{(p)})^{-1}\cdot \Delta.\end{array}$$ 
Note that, in particular, 
$$\Lambda\equiv 1\ \  \text{mod}\ \ \ p.$$
Then the Christoffel symbol of the first kind  is given by
$$\Gamma:=\Delta^t \cdot \phi(q)\cdot x^{(p)}=
\Delta^t \cdot (x^{(p)t})^{-1}\cdot A=S^t\cdot A=\frac{1}{p}(\Lambda^t-1)\cdot A.$$
Also the commutativity of \ref{got} is easily seen to be equivalent to the matrix equality
\begin{equation}
\label{LALB}
\Lambda^t A \Lambda=B.\end{equation}
\end{remark}

\begin{remark}\label{ozu}  An invariant interpretation of our Christoffel symbol
 can be given  in terms of concepts 
 introduced in \cite{foundations}, section 3.7.  
 This invariant interpretation   is not especially helpful when it comes to our  proofs and will not be used in what follows. However,
 for the reader familiar with \cite{foundations}, section 3.7, we briefly review this interpretation. First we recall that $\widehat{\mathfrak g}$ can be canonically identified as a $p$-adic formal scheme (but not as a group) with what in loc. cit. is referred to as the {\it arithmetic Lie algebra}  of $G$, denoted there by $L_{\d}(G)$, and defined as the kernel of the natural projection
  $\pi:J^1(G)\ra \widehat{G}$ from the {\it first $p$-jet space} $J^1(G)$ of $G$ to $\widehat{G}$. Then the matrix $S$ in Remark \ref{preozu} corresponds to the   ``quotient" of the two sections of $\pi$ defined by the $p$-adic connection $\d^G$ and the {\it trivial} $p$-adic connection $\d^G_0$ on $G$, respectively. This ``quotient" is computed in $J^1(G)$ and lies in $L_{\d}(G)$ so it defines a map $\widehat{G}\ra \widehat{\mathfrak g}$ and hence an element of 
 ${\mathfrak g}(\cO(\widehat{G}))$. 
 Alternatively $S$ can be obtained by taking the ``quotient" of the images of the ``identity" $id\in G(\cO(\widehat{G}))$ via the two set theoretic maps
 $G(\cO(\widehat{G}))\ra G(W_2(\cO(\widehat{G})))$ induced by the two ring homomorphisms
 $\cO(\widehat{G}))\ra W_2(\cO(\widehat{G}))$ corresponding to $\d^G$ and $\d^G_0$, respectively. Here $W_2$ stands for the functor of Witt vectors of length $2$.
 Next note that $\Gamma$ is obtained from  $S^t$ via right multiplication by $A$; this right multiplication operation plays the role of  ``lowering the indices" in the classical setting; cf. \ref{cris}. By the way $A$ also has an intrinsic interpretation since $A=\phi^G_0({\mathcal H}_q(x))$. The above makes our Christoffel symbol an analogue of the classical one; indeed, in the classical case, a similar description of the Christoffel symbol can be given in which the first $p$-jet space  above is replaced by the first jet space in the sense of Lie and  Cartan and the matrix  $A=x^{(p)t}\cdot \phi(q)\cdot x^{(p)}$ is replaced by $q$; cf. the Appendix. \end{remark}

\begin{remark}
Metric connections satisfy congruences that are reminiscent of identities in classical differential geometry. We explain this in what follows. 
Let $Z$ be the center of $G$ and let $T$ be the diagonal maximal torus of $G$.  
Consider an ideal $J\subset \cO(\widehat{G})$ and assume  one of the following $2$ situations:

\medskip

1) $J$ is the ideal defining $Z$;

2)  $J$ is the ideal  defining $T$ and $q\in T(\cO)$.

\medskip

\noindent In situation 1) $J$ is of course generated by 
$$\{x_{ii}-x_{jj},\ x_{ij}\ ;\ i,j=1,...,n,\ i\neq j\}$$
while in situation 2) $J$ is generated by
$$\{x_{ij};\ i,j=1,...,n,\ i\neq j\}.$$
Taking determinants in \ref{LALB}
we immediately get that, in either of the situations 1) or 2) above,
\begin{equation}
\det(\Lambda)\equiv \{\det(1_n+p(q^{(p)})^{-1}\cdot \d q)\}^{-1/2}\ \ \ \text{mod}\ \ \ J,
\end{equation}
where the $-1/2$ root is the one that is $\equiv 1$ mod $p$.
Explicitly if 
$$\eta:=\frac{1}{p}\{\det(1_n+p(q^{(p)})^{-1}\cdot \d q)-1\},$$
then
$$\det(\Lambda)\equiv \sum_{k=0}^{\infty} \left(\begin{array}{c} - 1/2\\ k\end{array}\right) p^k\eta^k\ \ \ \text{mod}\ \ J.$$
In particular, since $\Lambda=1+p(x^{(p)})^{-1}\Delta$, and since one has
$$\det(1+pM)\equiv 1+p\cdot \text{tr}(M)\ \ \ \text{mod}\ \ \ p^2$$  for any matrix $M$ with coefficients in any ring,  we get, in our case,
\begin{equation}
\text{tr}((x^{(p)})^{-1}\Delta)\equiv -\frac{1}{2}\cdot \text{tr}((q^{(p)})^{-1}\cdot \d q)\ \ \ \text{mod}\ \ \ (p,J).
\end{equation}
This congruence is analogous to an identity in classical differential geometry; cf. \ref{someyes}.
 
Assume now we are in situation 2) above and, in addition $q$ has entries in the ring $\bZ_p$ of $p$-adic integers.
Then we get the congruence
\begin{equation}
\det(\Lambda)\equiv \left(\frac{\det(q)}{p}\right)\cdot \det(q)^{\frac{p-1}{2}}\ \ \ \text{mod}\ \ \ J,
\end{equation}
where $\left(\frac{\ }{p}\right)$ is the Legendre symbol.
\end{remark}

\subsection{Levi-Civita connections: the global situation}
We introduce the concepts of {\it torsion free}  connections and {\it Levi-Civita} connections. We will first place ourselves in a global situation. Later we will go back to the local situation.

 Let $F$ be a number field which is Galois over ${\mathbb Q}$ 
 and let $M\in \bZ$ be an even integer divisible by the discriminant of $F$. The giving of   the data $F,M$ will be referred to as the {\it global situation}. In such a situation we
 let $\cO_F$ be the ring of integers of $F$ and set
$$\cO:=\cO_{F,M}:=\cO_F[1/M].$$
One can consider the Galois group ${\mathfrak S}(F)$ of $F/{\mathbb Q}$ and the natural map
$$Spec\ \cO\ra {\mathfrak S}(F)$$
sending any non-zero prime ${\mathfrak P}\in Spec\ \cO$ into the Frobenius element $\phi_{\mathfrak P}\in {\mathfrak S}(F)$ at ${\mathfrak P}$; the image of the zero prime will play no role and may be taken to be  the identity.

Assume now we are in this global situation.

\begin{definition}
\label{soldering}
 By a {\it vertical gauge} we will understand a map 
 \begin{equation}
 \label{vertical gauge}
 \{1,...,n\}\ra {\mathfrak S}(F),\ \ i\mapsto \sigma_i,\ \ \sigma_1=\text{id}.\end{equation}
 \end{definition}

 Morally a vertical gauge can be used to ``move vertically" in $Spec\ \cO$, above each prime in $\bZ$, using the Galois action. Here is how this works.

 Let  $p$ be a prime in $\bZ$
not dividing $M$ and let 
${\mathfrak P}$ be a prime ideal in $\cO_F$ dividing $p$.
Assume   that a vertical gauge is given and set
$${\mathfrak P}_i=\sigma_i {\mathfrak P},$$
so ${\mathfrak P}_1={\mathfrak P}$. Note that the ${\mathfrak P}_i$'s are not necessarily distinct.
Let $\cO_{{\mathfrak P}_i}$ be the localization of $\cO_F$ at ${\mathfrak P}_i$ and let, as usual, $\widehat{\cO_{{\mathfrak P}_i}}$ denote the $p$-adic completion of this localization. 
So  $\widehat{\cO_{{\mathfrak P}_i}}$ is in the {\it local situation} considered previously.
In particular $\widehat{\cO_{{\mathfrak P}_i}}$ has a unique Frobenius lift $\phi^{\widehat{\cO_{{\mathfrak P}_i}}}$ which, for simplicity, we denote by $\phi^i$. 
Consequently $\widehat{\cO_{{\mathfrak P}_i}}$ has a unique $p$-derivation $\d^{\widehat{\cO_{{\mathfrak P}_i}}}$ which, for simplicity, we denote by $\d^i$. 
Clearly $\phi^i$ sends $\cO_F$ into itself and the restriction of $\phi^i$ to $\cO_F$, further extended to an automorphism of $F$, is the usual Frobenius element  $\phi_{{\mathfrak P}_i} \in {\mathfrak S}(F)$, which we continue to denote by $\phi^i$. Of course, $\phi^i$ induces a Frobenius lift on $\cO_{{\mathfrak P}_i}$ but does not
 generally induce a Frobenius lift on $\cO_F$. If $F$ is abelian over ${\mathbb Q}$ then, of course, $\phi^i$ {\it does} induce a Frobenius lift on $\cO_F$.

As a matter of notation we will sometimes simply write $\d$ and $\phi$ instead of $\d^1$ and $\phi^1$; but we will {\it never} abbreviate $\d^i,\phi^i$ by $\d,\phi$ if $i\neq 1$.

Consider
the general linear groups over $\cO_{F,M}$ and $\widehat{\cO_{{\mathfrak P}_i}}$ respectively:
$$\begin{array}{rclll}
G & = & GL_n & = & Spec\ \cO_{F,M}[x,\det(x)^{-1}]\\
\ & \ & \ & \ & \ \\
G^i & = & GL_n\otimes \widehat{\cO_{{\mathfrak P}_i}} &=&
Spec\ \widehat{\cO_{{\mathfrak P}_i}}[x,\det(x)^{-1}].\end{array}$$
Note that the $G^i$'s are not necessarily distinct.
We have   induced isomorphisms (still denoted by) $\sigma_i:\widehat{\cO_{{\mathfrak P}_1}}\ra \widehat{\cO_{{\mathfrak P}_i}}$ extending uniquely to isomorphisms (still denoted by)
\begin{equation}
\label{drilling}
\sigma_i:\cO(\widehat{G^1})\ra \cO(\widehat{G^i}),\ \ \sigma_i(x)=x.\end{equation}

\begin{definition}
A {\it vertical connection} on $G$ at ${\mathfrak P}$ is an $n$-tuple $(\d^1,...,\d^n)$  where $\d^i$ is a $p$-adic connection
on $G^i$. 
The {\it Christoffel symbols} (of the first kind) relative to $q$ of a vertical connection are 
the Christoffel symbols of the first kind $\gamma_1,...,\gamma_n$ of $\d^1,...,\d^n$ relative to $q$, i.e., the matrices
\begin{equation}
\label{roro}
\gamma_i:=\d^i x^t \cdot \phi^i(q) \cdot x^{(p)}\in {\mathfrak g}{\mathfrak l}_n(\cO(\widehat{G^i})),\ \ \ i=1,...,n.\end{equation}
If  $\gamma_i=:(\gamma_{ijk})$ we say that the $n$-tuple 
$(\d^1,...,\d^n)$
is {\it torsion free}  (relative to $q$) if for all $i,j,k=1,...,n$ we have the following equalities in $\cO(\widehat{G^1})$:
\begin{equation}
\label{tortor}
\sigma_i^{-1}\gamma_{ijk}=\sigma_j^{-1}\gamma_{jik}.\end{equation}
\end{definition}

\begin{remark}\label{cih}\

1) If $F$ is abelian our notion of vertical connection above coincides with the one in \cite{foundations}.

2) If we define
\begin{equation}
\gamma_i'=\d^i x^t \cdot \phi^i(q),\ \ \ \gamma'_i=(\gamma'_{ijk}),
\end{equation}
then \ref{tortor} holds if and only if
\begin{equation}
\label{tortorprime}
\sigma_i^{-1}\gamma_{ijk}'=\sigma_j^{-1}\gamma_{jik}'.\end{equation}

3)
 The symmetry \ref{tortor} is an analogue of the  symmetry \ref{tf} in the definition of classical torsion freeness.

4) We will usually denote by $(\phi^1,...,\phi^n)$ the Frobenius lifts attached to $(\d^1,...,\d^n)$. So $\phi^i:\cO(\widehat{G^i})\ra \cO(\widehat{G^i})$ is a Frobenius lift, not to be mixed up with one of the maps $\sigma_i$ in \ref{drilling} which are never Frobenius lifts.
A confusion in notation may arise if $\phi^i=\sigma_j$ as elements in ${\mathfrak S}(F)$; in order to avoid  this confusion, when using the letter $\phi$ we will always mean a Frobenius lift and not one of the maps in \ref{drilling}.
\end{remark}

Assume we are in the global situation, 
we are given a vertical gauge, a prime $p\not| M$ and a prime ${\mathfrak P}|p$. Here is our first main result; it is an analogue of the ``Fundamental Theorem of Riemannian Geometry"; cf. Theorem \ref{fundfund} in  the Appendix.

 \begin{theorem}
 \label{levi Civita}
 Assume $q\in GL_n(\cO_{F,M})$, $q^t=q$. Then there exists  a unique vertical connection 
$(\d^1,...,\d^n)$ on $G$  at ${\mathfrak P}$ such that the following hold:

1)   $\d^i$ is  metric with respect to $q$ for all $i$;
  
  2) $(\d^1,...,\d^n)$ is torsion free relative to $q$.
  
 \end{theorem}

 \begin{definition}\label{levilevi}
The vertical connection $(\d^1,...,\d^n)$
in Theorem \ref{levi Civita}
 is called the  {\it vertical Levi-Civita connection}  attached to $q$ at ${\mathfrak P}$.  
 \end{definition}
 
 Next we want to vary $p$ and ${\mathfrak P}$. Assume we are in a global situation $F, M$. We make the following:
 
 \begin{definition}
 \label{papu}
 A {\it transversal gauge} consists of a set ${\mathcal V}$ of primes in $\bZ$ not dividing $M$ together with a map
 \begin{equation}
 \label{transversal gauge}
 {\mathcal V}\ra Spec\ \cO,\ \ \ p\mapsto {\mathfrak P}(p)\end{equation}
 such that ${\mathfrak P}(p)|p$ for all $p$.
 Given $q\in GL_n(\cO_{F,M})$, a vertical gauge, and a transversal gauge,
 the family of all vertical Levi-Civita connections 
 $(\d^1_p,...,\d^n_p)$ on $G$ at ${\mathfrak P}(p)$, where
  $p$ varies 
  in ${\mathcal V}$, will be referred to as the {\it mixed} (or {\it adelic}) {\it Levi-Civita connection}  attached to $q$; one can view it as a family 
  \begin{equation}
  (\d^i_p)\end{equation}
  depending on $2$ indices $i$ and $p$. The family 
  \begin{equation}
  \label{transversall}
  (\d^1_p)\end{equation}
  depending on one index $p$ only will be referred to as the {\it transversal  Levi-Civita connection}  attached to $q$. 
 Let us say that a vertical gauge \ref{vertical gauge}
  is {\it perfect} if the map \ref{vertical gauge}
   is  bijective. 
   Let us say that a transversal gauge \ref{transversal gauge} is {\it perfect}
  if the map
  \begin{equation}
  \label{ppop}
   {\mathcal V}\ra Spec\ \cO \ra {\mathfrak S}(F),\ p\mapsto \phi_{{\mathfrak P}(p)}
   \end{equation}
 is bijective; if this is the case the inverse ${\mathfrak S}(F)\ra {\mathcal V}$ of the map \ref{ppop} will be denoted by  $\sigma\mapsto p({\sigma})$. If a perfect vertical gauge and a perfect transversal gauge are given then the transversal Levi-Civita connection
 \ref{transversall} can  be viewed as a family 
 \begin{equation}
 \label{asa}
 (\d^1_{p(\sigma_i)})\end{equation}
  indexed by $i\in \{1,...,n\}$. 
 \end{definition}
 
  The transversal Levi-Civita connection can be viewed
  as an analogue of the  ``standard" Levi-Civita connection in classical differential geometry; this is not clear a priori and will be explained in the Appendix by way of introducing differential geometric analogues of our vertical and mixed Levi-Civita connections.
  In some sense the vertical and mixed versions of these connections turn out to be more fundamental than then transversal one.
 
  \medskip
  
  Assume the notation of Theorem \ref{levi Civita} and consider the matrices
  \begin{equation}
  \label{cil}
  C_i:=- x^{(p)t}\cdot \sigma_i^{-1}\d^i q \cdot x^{(p)} +\frac{1}{p}\{(x^t \cdot \sigma_i^{-1}q \cdot x)^{(p)}-x^{(p)t}\cdot \sigma_i^{-1}q^{(p)}\cdot x^{(p)}\}.\end{equation}
  If $C_i=(C_{ijk})$ then, clearly,  $$C_{ijk}=C_{ikj}.$$

We will show:

\begin{proposition}
\label{voices}
The following congruences hold in $\cO(\widehat{G^1})$:
\begin{equation}
\label{greene}
\begin{array}{rcll}
\sigma_i^{-1}\gamma_{ijk} & \equiv & \frac{1}{2}(C_{ijk}+C_{jik}-C_{kij}) & \text{mod}\ \ p,\\
\  & \  & \  &\  \\
\sigma_i^{-1}\gamma_{ijk}  & \equiv & - \frac{1}{2}(\sigma_i^{-1} \d^{i} q_{jk}+\sigma_j^{-1} \d^{j}  q_{ik}-
\sigma_k^{-1}\d^{k}  q_{ij}) & \text{mod}\ \ (p,x-1).
\end{array}
\end{equation}
\end{proposition}

Recall that $1$ is the identity matrix so  $(p,x-1)$ is the ideal generated by $p$, $x_{ii}-1$, and $x_{ij}$
for $i\neq j$. The formula \ref{greene} should be viewed as an analogue of the classical
expression for the Levi-Civita connection in classical Riemannian geometry; cf. \ref{windy} in the Appendix.

\subsection{Levi Civita connection: the local situation}
We will (directly) deduce Theorem \ref{levi Civita} and Proposition \ref{voices} from corresponding local results; cf. Theorem \ref{LC} and Proposition \ref{congruences} below. We need some notation. 

Assume in what follows that we are in the local situation; so $\cO$ is a complete discrete valuation ring with maximal ideal generated by $p$ and perfect residue field, viewed  as equipped with its unique Frobenius lift $\phi=\phi^{\cO}$ and its unique $p$-derivation $\d=\d^{\cO}$.
Set, in this situation,
$$G=GL_n=Spec\ \cO[x,\det(x)^{-1}].$$

\begin{definition}
Assume $q_1,...,q_n\in GL_n(\cO)$, $q_i^t=q_i$, let  
$(\d^{G}_1,...,\d^{G}_n)$  be an $n$-tuple of $p$-adic connections  on $G$ and let 
 $\Gamma_1,...,\Gamma_n$ be  the Christoffel symbols of the first kind of
$\d^G_1,...,\d^G_n$ with respect to $q_1,...,q_n$ respectively; explicitly,
$$
\Gamma_i=:\Gamma_i[q_1,...,q_n]:=\d_i^G x^t \cdot \phi(q_i)\cdot x^{(p)},\ \ \ i=1,...,n.
$$
 Set $\Gamma_i=(\Gamma_{ijk})$. We say that $(\d^{G}_1,...,\d^{G}_n)$ is {\it torsion free} relative to $(q_1,...,q_n)$ if for all $i,j,k$, we have
$$\Gamma_{ijk}=\Gamma_{jik}.$$
\end{definition}

For  invariant formulations of the above condition  see Remark \ref{ozi}.

We will prove the following:

\begin{theorem}
\label{LC}
Assume $q_1,...,q_n\in GL_n(\cO)$, $q_i^t=q_i$. Then there exists  a unique $n$-tuple
$(\d^{G}_1,...,\d^{G}_n)$ of $p$-adic connections  on $G$ such that  the following hold:

1)   $\d^{G}_i$ is  metric with respect to $q_i$ for all $i$;
  
  2) $(\d^{G}_1,...,\d^{G}_n)$  is torsion free relative to $(q_1,...,q_n)$.
\end{theorem}

\begin{definition}
\label{locallc}
The tuple $(\d^{G}_1,...,\d^{G}_n)$ is called the {\it  Levi-Civita connection}
on $G$ over $\cO$ attached to the tuple $(q_1,...,q_n)$.\end{definition}

\begin{remark}
\label{ozi}

Let 
$\d_i^G x=:\Delta_i$  let $\phi^G_i$ be the corresponding Frobenius lifts, and write
$$\phi^G_i(x)=\Phi_i=x^{(p)}+p\Delta_i=x^{(p)}\Lambda_i,\ \ \ A_i=x^{(p)t}\phi(q_i)x^{(p)},$$
$$\Delta_i=(\Delta_{ijk}),\ \ \ \Lambda_i=(\Lambda_{ijk}).$$
Then the following hold:

1) Condition 2 in Theorem \ref{LC} is equivalent to the condition
$$(\phi(q_i)\cdot \Delta_i)_{kj}=(\phi(q_j)\cdot \Delta_j)_{ki}$$
and also to the condition
$$(A_i(\Lambda_i-1))_{kj}=(A_j(\Lambda_j-1))_{ki}.$$

2) If $q_1=...=q_n$ then 
Condition 2 in Theorem \ref{LC} is equivalent to the condition
$$\Delta_{ikj}= \Delta_{jki}$$
and also to the condition
$$(\Lambda_i-1_n)_{kj}=(\Lambda_j-1_n)_{ki}.$$
These conditions are also equivalent to the commutativity of the following diagrams of formal schemes over $\cO$:
\begin{equation}
\label{risj}
\begin{array}{rcl}
\widehat{G} & \stackrel{s_i}{\longrightarrow} & \widehat{\mathfrak g}\\
s_j \downarrow & \ & \downarrow r_j\\
\widehat{\mathfrak g} & \stackrel{r_i}{\longrightarrow} & \widehat{{\mathbb A}^n}
\end{array}
\end{equation}
where 
$${\mathfrak g}:=Spec\ \cO[x],\ \ {\mathbb A}^n=Spec\ \cO[z_1,...,z_n],$$
$$r_i(z_k)=x_{ki},\ \ s_i(x)=(x^{(p)})^{-1}\cdot \Delta_i.$$
The commutativity of the diagrams \ref{risj} 
is analogous to the commutativity of the diagrams 
\ref{risju} that define torsion freeness in classical differential geometry. Also the commutativity of \ref{risj}
has an invariant meaning in terms of concepts 
 introduced in \cite{foundations}, section 3.7. 
 Indeed, as already mentioned in Remark \ref{ozu}, $s_i$  is the ``quotient" of the two sections of the projection $J^1(G)\ra \widehat{G}$ corresponding to the $p$-adic connection $\d_i^G$ and the  trivial $p$-adic connection $\d^G_0$ on $G$.

3) If  all $q_i$ are scalar matrices, $q_i=d_i\cdot 1_n$ then Condition 2 in Theorem \ref{LC} is equivalent   to the condition
$$\phi(d_i)\cdot (\Lambda_i-1_n)_{kj}=\phi(d_j)\cdot (\Lambda_j-1_n)_{ki}.$$

 4) Theorems \ref{LC} and \ref{levi Civita} are related as follows. Assume, 
 for the remainder of  this Remark only, that we place ourselves in the global situation, with a given vertical gauge and a given transversal gauge.
 Moreover let $q\in GL_n({\mathcal O}_{F,M})$, $q^t=q$.  Then our proofs will show that the vertical Levi-Civita connection $(\d^1,...,\d^n)=(\d^1_p,...,\d^n_p)$ attached to $q$ at ${\mathfrak P}(p)$ (cf. Theorem \ref{levi Civita}) is given by 
 $$\d^i=\sigma_i\circ \d_i^{G^1}  \circ \sigma_i^{-1},$$
 where $(\d_{p1},...\d_{pn}):=(\d_1^{G^1},...,\d_n^{G^1})$ is the  Levi-Civita connection on $G^1$ over $\widehat{\cO_{{\mathfrak P}(p)}}$ attached to 
 \begin{equation}
 \label{secret}
 (\sigma_1^{-1}q,...,\sigma_n^{-1}q)\end{equation} 
 as in Theorem \ref{LC}. 
 The family of all $(\d_{p1},...,\d_{pn})$ as $p$ varies in ${\mathcal V}$
 will be referred to as the {\it mixed} (or {\it adelic}) {\it  Levi-Civita connection} on $G$ attached to
 \ref{secret}. It can be viewed as a family 
 $$(\d_{pi})$$
 indexed by $2$ indices $p$ and $i$. The family 
 $$(\d_{p1})$$ indexed by $p$ only will be called the  {\it transversal Levi-Civita connection}
 on $G$ attached to \ref{secret}. 
  Again, these concepts are analogous to concepts in classical differential geometry; cf. our Appendix.\end{remark}

In the following discussion we are, again, in the  local situation, i.e., over a complete discrete valuation ring $\cO$ with maximal ideal generated by $p$ and perfect residue field.
One can ask about the dependence of the Christoffel symbols on $(q_1,...,q_n)$; the answer is that 
Christoffel symbols are ``universal $\d$-functions" of $(q_1,...,q_n)$ of order $1$ in the sense of \cite{char,book,foundations}.
Rather than recalling the general definition of $\d$-functions from loc.cit. we will explain this in an ad hoc manner. To do so let us consider $n\times n$ symmetric matrices $$s_1,...,s_n,s'_1,...,s'_n$$ with indeterminate entries on and above the diagonal
and consider the ring
\begin{equation}
\label{ncis}
\bZ_p[x,s_1,...,s_n,s'_1,...,s'_n,\det(x)^{-1},\det(s_1)^{-1},...,\det(s_n)^{-1}]^{\widehat{\ }}.
\end{equation}
Its elements 
$$f(x,s_1,...,s_n,s'_1,...,s'_n)$$
can be specialized by replacing $s_i,s'_i$ with symmetric matrices in $GL_n(\cO)$ to give elements
in $\cO[x,\det(x)^{-1}]^{\widehat{\ }}$.
Then we will prove:

\begin{proposition}
\label{unichris}
There exist $n\times n$ matrices 
$$\Gamma_i^{\text{univ}}(x,s_1,...,s_n,s'_1,...,s'_n)$$
 with entries in the ring \ref{ncis}, depending only on $p$ and $n$ (but not on $\cO$)
such that for any symmetric matrices $q_1,...,q_n\in GL_n(\cO)$,  the Christoffel symbols of the  Levi-Civita connection attached to $(q_1,...,q_n)$ are given by
$$\Gamma_i[q_1,...,q_n]=\Gamma_i^{\text{univ}}(x,q_1,...,q_n,\d q_1,...,\d q_n).$$
\end{proposition}

The Christoffel symbols satisfy some remarkable congruences. Indeed
assume the situation in Theorem \ref{LC} and consider the matrices
\begin{equation}
\label{ci}
C_i:=- x^{(p)t}\cdot \d q_i \cdot x^{(p)} +\frac{1}{p}\{(x^tq_ix)^{(p)}-x^{(p)t}q_i^{(p)}x^{(p)}\}.\end{equation}
If $C_i=(C_{ijk})$ then, clearly, 
$$C_{ijk}=C_{ikj}.$$
We will show:

\begin{proposition}
\label{congruences}
 The following congruences hold in $\cO(\widehat{G})$:
$$
\begin{array}{rcll}
\Gamma_{ijk} & \equiv & \frac{1}{2}(C_{ijk}+C_{jik}-C_{kij}) & \text{mod}\ \ p\\
\  & \  & \  &\ \\
\Gamma_{ijk}  & \equiv & - \frac{1}{2}(\d q_{ijk}+\d q_{jik}-\d q_{kij}) & \text{mod}\ \ (p,x-1).
\end{array}
$$
\end{proposition}

\subsection{Case $n=1$}
If in Theorem \ref{LC} we assume $n=1$ then Condition 2 in that theorem is, of course, automatically satisfied. Also $x$ is {\it one indeterminate} and we can write $\d^G_1=\d^G$, $q_{11}=d$. Then  Condition 1 is trivially seen to be equivalent to the condition that the Frobenius lift $$\phi^G:\cO[x,x^{-1}]^{\widehat{\ }}\ra \cO[x,x^{-1}]^{\widehat{\ }}$$ attached to the $p$-adic connection $\d^G$ satisfy
\begin{equation}
\phi^G(x)=\left(\frac{d^p}{\phi(d)}\right)^{1/2}\cdot x^p,
\end{equation}
where the square root is  chosen to be $\equiv 1$ mod $p$, i.e.,
\begin{equation}
\left(\frac{d^p}{\phi(d)}\right)^{1/2}=\left(1+p\frac{\d d}{d^p}\right)^{-1/2}:=
\sum_{k=0}^{\infty} \left( \begin{array}{c} - 1/2\\ k\end{array}\right)
p^k\left(\frac{\d d}{d^p}\right)^k.
\end{equation}
So in case $d\in \bZ_p^{\times}$ we have the formula
\begin{equation}
\phi^G(x)=\left(\frac{d}{p}\right)\cdot d^{\frac{p-1}{2}}\cdot x^p,
\end{equation}
where $\left(\frac{\ }{p}\right)$ is the Legendre symbol. 

Note that in case $n=1$ the  Levi-Civita connection introduced in Definition \ref{locallc} coincides with the {\it real Chern connection} introduced in \cite{foundations}, Introduction (or Definition 4.25). On the other hand, if $n\geq 2$, the Levi-Civita and the real Chern connection are different objects.

\subsection{Case $n=2$}
In this case there is an  analog of ``compatibility with complex structure" and ``conformal coordinates" which we discuss next.

We start by considering group schemes over $\cO$ and group scheme homomorphims over $\cO$,
\begin{equation}
\label{teatea}
\begin{array}{cclcc}
\ & \ & G''' & \ & \ \\
\ & \ & \uparrow \det & \ & \ \\
G'' & \ra & G' & \ra & G\\
\ & \ & \downarrow \det^{\perp} & \ & \ \\
\ & \ & G''' & \ & \ 
\end{array}\end{equation}
defined as follows. (The horizontal homomorphisms will be closed immersions; the homomorphism $\det^{\perp}$ will only be defined for $\sqrt{-1}\in \cO$.)

We start be letting 
$$G=GL_2=Spec\ \cO[x,\det(x)^{-1}],$$
with $x$ a $2\times 2$ matrix of indeterminates.
Then we let $G':=GL_1^c$ be the centralizer subgroup scheme in $G$ of the matrix
\begin{equation}
c=\left(\begin{array}{cc}
0 & 1\\ -1 & 0\end{array}\right).
\end{equation}
The matrix $c$ can be viewed as an analogue of ``complex structure" and $GL_1^c$ can be viewed as the ``complexified $GL_1$". 
One has
$$G'=Spec\ \cO[\alpha,\beta,(\alpha^2+\beta^2)^{-1}]$$
for $\alpha, \beta$ two indeterminates, with $G'$ embedded into $G$ via
 $$x\mapsto \left(\begin{array}{cc}\alpha  & \beta\\ - \beta & \alpha\end{array}
\right).$$
Note that a $2\times 2$ matrix $q$ that is symmetric and in $GL_1^c$ must be scalar; such a matrix can be viewed as an analogue of a  ``conformal metric".
Next we set 
$$G'''=GL_1=Spec\ \cO[z,z^{-1}],$$ 
we consider the group scheme homomorphism $$\det:G'\ra G'''$$ induced by
 $$z\mapsto \alpha^2+\beta^2,$$
 and we  consider the kernel 
\begin{equation}
G'':=U^c_1:=\text{Ker}(\det:G'\ra G'''),
\end{equation}
which can be viewed as the  ``complexified unitary group in dimension $1$." Of course,
$$G''=Spec\ \frac{\cO [\alpha,\beta]}{(\alpha^2+\beta^2-1)}.$$

Finally, if we assume, in addition, that $\sqrt{-1}\in \cO$,  we may consider the group scheme homomorphism
$$\text{det}^{\perp}: G'\ra G'''$$
induced  by 
$$z\mapsto s:= \frac{\alpha+\sqrt{-1}\beta}{\alpha-\sqrt{-1}\beta}.$$
This ends our definition of the objects in \ref{teatea}.
We will prove:

\begin{theorem}
\label{mock}
Assume
 $$d_1,d_2\in \cO^{\times},\ \ \ q_1=d_1\cdot 1_2\in G(\cO),\ \ \ q_2=d_2\cdot 1_2\in G(\cO)$$
 and let $(\d_1^G,\d_2^G)$
be the  Levi-Civita connection on $G$ over $\cO$ attached to $(q_1,q_2)$. Then:

\medskip

1) $G'$ is $\d^G_i$-horizontal for $i=1,2$.

2) $G''$ is $\d^G_i$-horizontal for $i=1,2$ if and only if $\d d_1=\d d_2=0$.

3) Assume $d_1=d_2=:d \in \bZ$, $d\neq \pm 1$, $d\not\equiv 0$ mod $p$.
Then there is no closed connected proper subgroup scheme of $G'$ that is $\d^G_i$-horizontal for $i=1,2$.
\end{theorem}

\begin{definition}
\label{LCprime}
The pair 
$(\d_1^{G'},\d_2^{G'})$ of $p$-derivations on $G'$  induced by $(\d_1^G,\d_2^G)$ 
(which exist by assertion 1 in  Theorem \ref{mock})
is called the {\it  Levi-Civita connection} on $G'=GL_1^c$ over $\cO$ attached to $(q_1,q_2)$.
\end{definition}

\begin{remark}
\label{fruti}
The Frobenius lifts $\phi_1^{G'},\phi_2^{G'}$ attached to the  Levi-Civita connection on $G'$ attached to $(q_1,q_2)$ do not commute in general; this will be seen when we discuss curvature, cf. Remark \ref{tutti fruti}. Also note that assertion 3 in Theorem \ref{mock} intutively says that the Levi-Civita connection on $G'$, induced from that on $G$, does not induce, in its turn, a connection on any connected proper subgroup of $G'$; this can be viewed as  an ``irreducibility" (or a ``transitivity") statement.
\end{remark}

The objects in Theorem \ref{mock} can be described explicitly. Indeed assume  the situation and notation in that Proposition  and let 
 $(\phi_1^{G'},\phi^{G'}_2)$ be the corresponding Frobenius lifts on $\widehat{G'}$. Let
$$\epsilon:=\phi \left(\frac{d_2}{d_1}\right)\in \cO^{\times},$$
and set
$$\theta_i:= \frac{d_i^p(\alpha^2+\beta^2)^p}{\phi(d_i)(\alpha^{2p}+\beta^{2p})}\in \cO(\widehat{G'})^{\times},\ \ \ i=1,2.
$$
Then the system
\begin{equation}
\begin{array}{rcl}
\label{hon}
x_1^2-2\epsilon x_2+\epsilon^2 x_2^2 & = & \theta_1-1\\
\  & \  & \  \\
x_1^2+2\epsilon x_1+\epsilon^2 x_2^2 & = & \epsilon^2(\theta_2-1)
\end{array}
\end{equation}
with unkowns $x_1,x_2$
is trivially seen to have a unique solution
\begin{equation}
\label{bs}
(v_1,v_2)\end{equation}
in the set
$$p\cO(\widehat{G'})\times p\cO(\widehat{G'}).$$
The solution can be computed explicitly, of course,
the way one finds the intersection of two circles in analytic geometry: 
 one takes   the difference of the equations in \ref{hon},
which is a linear equation,
$$2\epsilon (x_1+x_2)=\epsilon^2(\theta_2-1)-(\theta_1-1),$$
one solves the latter for $x_2$,  one substitutes in one of the equations \ref{hon}, and one solves  the resulting quadratic equation by the quadratic formula;   the radical involved needs to be expressed as a $p$-adic series. Define now
\begin{equation}
\label{as}
u_2:= 1+\epsilon^{-1}v_1,\  \ u_1:=1-\epsilon v_2\in \cO(\widehat{G'}).
\end{equation}

\begin{proposition}
\label{fishatnoon}
We have the following equality of matrices with coefficients in $\cO(\widehat{G'})$:
\begin{equation}
\label{fid}
\phi_i^{G'}\left(\begin{array}{cc}\alpha & \beta\\ - \beta & \alpha\end{array}\right)=
\left(\begin{array}{cc}
\alpha^p & \beta^p \\ - \beta^p & \alpha^p\end{array}\right) \cdot \left( \begin{array}{rr} u_i & v_i\\ - v_i & u_i\end{array}\right),\ \ \ i=1,2.\end{equation}
\end{proposition}

\begin{remark}
\label{fishallday}
Assume that in the above discussion we have $d_1=d_2$. Then the formulas simplify as follows. One has 
$$\theta_1=\theta_2,\ \ \ u_1=u_2,\ \ \ v_1=-v_2$$
and if 
$$u:=u_1,\ \ \ v:=v_1,\ \ \ \theta:=\theta_1,\ \ \ d:=d_1,\ \ \ \eta:=(\theta-1)/p$$ then
$$u=1+v,\ \ \ u^2+v^2=\theta,\ \ \ 2v^2+2v+(1-\theta)=0,$$
and
$$v=-\frac{1}{2}+\frac{1}{2}(2\theta-1)^{1/2}:=-\frac{1}{2}+\frac{1}{2}\sum_{k=0}^{\infty}
\left(\begin{array}{c} 1/2\\ k\end{array}\right) 2^kp^k\eta^k.$$
Note the following congruence:
$$
v\equiv -\frac{p}{2}\cdot \frac{\d d}{d^p}\ \ \ \text{mod}\ \ \ (p^2,\alpha-1,\beta).
$$
This implies the congruences
\begin{equation}
\label{fidx}
\d_1^{G'}\left(\begin{array}{cc}\alpha & \beta\\ - \beta & \alpha\end{array}\right)\equiv -\frac{1}{2}
\left( \begin{array}{cc} \frac{\d d}{d^p} & \frac{\d d}{d^p}\\ 
\ &  \ \\ 
- \frac{\d d}{d^p}& \frac{\d d}{d^p}\end{array}\right),\ \ \ \text{mod}\ \ \ (p,\alpha-1,\beta),\end{equation}
\begin{equation}
\label{fidxx}
\d_2^{G'}\left(\begin{array}{cc}\alpha & \beta\\ - \beta & \alpha\end{array}\right)\equiv -\frac{1}{2}
\left( \begin{array}{cc} \frac{\d d}{d^p} & - \frac{\d d}{d^p}\\
\ &  \ \\ 
\frac{\d d}{d^p} & \frac{\d d}{d^p}\end{array}\right),\ \ \ \text{mod}\ \ \ (p,\alpha-1,\beta).\end{equation}

\end{remark}

\begin{proposition} \label{tirg}
Assume $d_1=d_2=d$. 

1) Consider the Frobenius lift 
$\phi^{G'''}:\widehat{G'''}\ra \widehat{G'''}$ 
on $G'''=Spec\ \cO[z,z^{-1}]$
defined by $\phi^{G'''}(z)=\frac{d^p}{\phi(d)}\cdot z^p$.
Then  the following diagrams are commutative:
\begin{equation}
\label{separatist}
\begin{array}{rcl}
\widehat{G'} & \stackrel{\phi^{G'}_i}{\longrightarrow} & \widehat{G'}\\
\det \downarrow & \  & \downarrow \det\\
\widehat{G'''} & \stackrel{\phi^{G'''}}{\longrightarrow} & \widehat{G'''}
\end{array}
\end{equation}

2) Assume that $\sqrt{-1}\in \cO$. 
Then there are (unique) Frobenius lifts $\phi^{G'''}_i$, $i=1,2$, on $\widehat{G'''}$ making the following diagrams commutative:
\begin{equation}
\label{separatistule}
\begin{array}{rcl}
\widehat{G'} & \stackrel{\phi^{G'}_i}{\longrightarrow} & \widehat{G'}\\
\text{det}^{\perp} \downarrow & \  & \downarrow \text{det}^{\perp}\\
\widehat{G'''} & \stackrel{\phi_i^{G'''}}{\longrightarrow} & \widehat{G'''}
\end{array}
\end{equation}

3)   $\phi_1^{G'}$ and $\phi_2^{G'}$ commute
if and only if $\phi_1^{G'''}$ and $\phi_2^{G'''}$ commute.
\end{proposition}

\begin{remark}
The Frobenius lift $\phi^{G'''}$  in assertion 1 of the Proposition coincides with the Frobenius lift attached to what in \cite{foundations}, Introduction (or Definition 4.33) was called the {\it complex Chern connection} on $GL_1$ attached to $q$. 
\end{remark}

\begin{remark}
As already mentioned  we will later see that $\phi_1^{G'}$ and $\phi_2^{G'}$ do not commute in general, cf. Remark \ref{tutti fruti}; so  $\phi_1^{G'''}$ and $\phi_2^{G'''}$, too, do not commute in general.
\end{remark}

\begin{remark}
Proposition \ref{tirg} is an analogue of a situation encountered in classical Riemannian geometry; in particular the commutator  of the Frobenius lifts
$$[\phi^{G'''}_1,\phi_2^{G'''}]:\widehat{\cO(G''')}\ra \widehat{\cO(G''')},$$
viewed as a function of $d$, 
should be viewed as an arithmetic analogue of the Laplacian composed with the logarithm,
cf. formula \ref{gogugogu} in the Appendix.
\end{remark}

\section{Curvature} 

In what follows we would like to define the
 curvature of vertical and mixed Levi-Civita connections; once one knows how to deal with the mixed case one can deal, of course,  with the transversal case as well.
 Recall that the vertical context refers to the case when we fix a prime $p$ and we ``vary" the primes of $F$ above $p$; in this context the definition of curvature is straightforward and 
 we will derive some basic congruences for its components that are reminiscent of formulae from classical differential geometry. The mixed context refers to the case when the prime $p$ is allowed to vary  while we still allow a ``vertical" variation of primes of $F$ above $p$; in this context the definition of curvature is more subtle: it is based on ``algebraization by correspondences" in a sense similar to  \cite{curvature2, foundations}. The two pictures corresponding to the two contexts above turn out to be different in general.

 \subsection{Vertical context}

 Assume we are in the global situation with data $F,M$, assume we are given a vertical gauge,  a prime $p \not| M$, and a prime ${\mathfrak P}|p$. We will use, in what follows,  the notation introduced after Definition \ref{soldering}.

\begin{definition} 
 Let $q\in GL_n(\cO_{F,M})$,  let $(\d^1,...,\d^n)$ be the  vertical Levi-Civita connection attached to $q$ at ${\mathfrak P}$, and let $(\phi^1,...,\phi^n)$ be the attached Frobenius lifts. The {\it curvature} of the (vertical) Levi-Civita connection is the family $(\varphi_{ij})$ where $i,j=1,...,n$ and $\varphi_{ij}:\cO(\widehat{G^1})\ra \cO(\widehat{G^1})$ are the (additive) maps
 \begin{equation}
 \label{zz}
 \varphi_{ij}:=\frac{1}{p}\{\sigma_i^{-1}\phi^i\sigma_i\sigma_j^{-1}\phi^j\sigma_j-
 \sigma_j^{-1}\phi^j\sigma_j\sigma_i^{-1}\phi^i\sigma_i\}.\end{equation}
 \end{definition}

  Our definition \ref{zz} of curvature is analogous to the classical definition of curvature; cf. \ref{cur} in the Appendix.

Let $C_i=(C_{ijk})$ be as in \ref{cil}, let $(q^{ij})$ be the inverse of the matrix $q=(q_{ij})$, and let $(x^{ij})$ be the inverse of the matrix $x=(x_{ij})$. Set $\Phi_{ij}:=\varphi_{ij}(x)$ and let 
$\Phi_{ijmk}$ be the entries of the matrices $\Phi_{ij}$, so
$$\Phi_{ij}=(\Phi_{ijmk});$$
we refer to $\Phi_{ijmk}$ as the {\it components} of the curvature.
 We will prove:
 
 \begin{proposition}
 \label{elisabeth}
  The following congruences hold in $\cO(\widehat{G^1})$:
\begin{equation}
\label{burton}
\begin{array}{rcll}
\Phi_{ijmk} & \equiv & \frac{1}{2}(\sigma_j^{-1}q^{ms})^{p^2} (x^{rs})^{p^2} 
(C_{jkr}+C_{kjr}-C_{rjk})^p\\
\ &\ & \ & \ \\
\ & \  & 
- \frac{1}{2}(\sigma_i^{-1}q^{ms})^{p^2} (x^{rs})^{p^2} 
(C_{ikr}+C_{kir}-C_{rik})^p & \ \\
\ &\ & \ & \ \\
\ & \  &   \text{mod}\ \ p, &\ \\
\ & \ & \ & \ \\
\Phi_{ijmk}   & \equiv \  &  \frac{1}{2}(\sigma_i^{-1}q^{mr})^{p^2} 
(\sigma_i^{-1}\d^i q_{kr}+\sigma_k^{-1}\d^k q_{ir}- \sigma_r^{-1}\d^r q_{ik})^p   & \ \\
\ &\ & \ & \ \\
\ & \  & 
-\frac{1}{2} (\sigma_j^{-1}q^{mr})^{p^2} 
(\sigma_j^{-1}\d^j q_{kr}+\sigma_k^{-1}\d^k q_{jr}- \sigma_r^{-1}\d^r q_{jk})^p 
&  \ \\
\ &\ & \ & \ \\
\ & \  & 
\text{mod}\ \ (p,x-1),
\end{array}\end{equation}
where the repeated indices $r,s$ are summed over. \end{proposition}

 It is worth noting that only the coefficients of the ``metric" and their first ``arithmetic derivatives" occur in these congruences. If, instead of congruences, one is interested in equalities then $\Phi_{ijmk}$ will ``involve" the ``arithmetic derivatives up to order $2$ of the metric" as in the case of classical Riemannian geometry. More precisely 
  let us consider $n\times n$ symmetric matrices $$s_1,...,s_n,s'_1,...,s'_n,s''_1,...,s''_n$$ with indeterminate entries on and above the diagonal
and consider the ring
\begin{equation}
\label{nciss}
\bZ_p[x,s_1,...,s_n,s'_1,...,s'_n,s''_1,...,s''_n,\det(x)^{-1},\det(s_1)^{-1},...,\det(s_n)^{-1}]^{\widehat{\ }}.
\end{equation}
Its elements 
$$f(x,s_1,...,s_n,s'_1,...,s'_n,s''_1,...,s''_n)$$
can be specialized by replacing $s_i,s'_i,s''_i$ with symmetric matrices in $GL_n(\cO)$ to give elements
in $\cO[x,\det(x)^{-1}]^{\widehat{\ }}$.
For simplicity we change the notation by writing  $\d:=\d^{\cO}$ and $\d^2=\d\circ\d$ on $\cO$.
Then Proposition \ref{unichris} easily implies that the curvature of the Levi-Civut\`{a} connection is given in terms of $q$ by some ``universal $\d$-functions of order $2$" in the following sense:

\begin{corollary}
\label{unichriss}
There exist $n\times n$ matrices 
$$
\Phi^{\text{univ}}_{ijmk}(x,s_1,...,s_n,s'_1,...,s'_n,s''_1,...,s''_n)$$
 with entries in the ring \ref{nciss}, depending only on $p$ and $n$ (but not on $F,M$ or the vertical gauge or the transversal gauge)
such that for any symmetric matrix $q\in GL_n(\cO)$,  the 
components $\Phi_{ijmk}=\Phi_{ijmk}[q]$ of the curvature of
the Levi-Civita connection attached to $q$ are given by
$$\Phi_{ijmk}[q]=\Phi_{ijmk}^{\text{univ}}(x,q_1,...,q_n,\d q_1,...,\d q_n,\d^2 q_1,...,\d^2 q_n),$$
where $q_i:=\sigma_i^{-1}q$.
\end{corollary}

\bigskip

A  significant simplification  occurs if
 one assumes 
\begin{equation}
\label{sigmaqequalq}
\sigma_i(q)=q,\ \ \ i=1,...,n,
\end{equation}
 in which case we say $q$ is {\it vertical gauge invariant}.
 So in this case, with
$\d=\d^1$, we have
\begin{equation}
\sigma_i^{-1}\d^i q=\sigma_i^{-1}\d^i\sigma_i\sigma_i^{-1}q=\d^1 q=\d q.
\end{equation}
Define
\begin{equation}
\label{defofC}C:=- x^{(p)t}\cdot \d q \cdot x^{(p)} +\frac{1}{p}\{(x^t \cdot q \cdot x)^{(p)}-x^{(p)t}\cdot q^{(p)}\cdot x^{(p)}\},\ \ \ C=(C_{jk})\end{equation}
and note that
$$C_{jk}=C_{kj}.$$
Then, recalling $C_i$ from \ref{cil}, we have $C_i=C$ and $C_{ijk}=C_{jk}$ for all $i,j,k$; also 
$$C_{jk}\equiv -\d q_{jk}\ \ \ \text{mod}\ \ \ (p,x-1)$$ for all $j,k$.
Define $R_{ijmk}$ as the $(m,k)$ entry of the matrix
 $$x^{(p^2)t}q^{(p^2)}\Phi_{ij};$$ 
in other words,
$$ x^{(p^2)t}q^{(p^2)}\Phi_{ij}=(R_{ijmk}).$$
We view $R_{ijmk}$ as an analogue of the classical covariant Riemann tensor; cf. \ref{notation monkey} in the Appendix.
Then the congruences \ref{burton} aquire the following simpler form.

\begin{corollary} 
\label{seinfeld1}
Assume  $q$ is vertical gauge invariant. Then the following congruences hold in $\cO(\widehat{G^1})$:
\begin{equation}
\label{tonton}
\begin{array}{rcll}
R_{ijmk} & \equiv &
\frac{1}{2}(C_{ik}+C_{jm}-C_{jk}-C_{im})^p
 &\text{mod}\ \ p\\
\ &\ & \ & \ \\
R_{ijmk}  & \equiv & 
 \frac{1}{2}(\d q_{jk}+\d q_{im}-\d q_{ik}- \d q_{jm})^p&  \text{mod}\ \ (p,x-1).
\end{array}
\end{equation}\end{corollary}

The congruences mod $(p,x-1)$ in \ref{tonton} are analogous to the formulae 
 for the covariant Riemannian tensor in  ``normal coordinates"; cf. \ref{traintrain} in the Appendix. Note however that, 
in arithmetic, one ``loses" one derivative;
so, in some sense,  the case when $q$ is vertical gauge invariant behaves 
(modulo the ``loss of one derivative") 
as if ``coordinates are already normal at $p$." 

\medskip

As a consequence of \ref{tonton} we get:

\begin{corollary} \label{seinfeld2}
Assume $q$  is vertical gauge invariant. Then the following congruences hold in $\cO(\widehat{G^1})$:
\begin{equation}
\label{macout}
\begin{array}{rcll}
R_{ijkm} & \equiv  &- R_{ijmk} & \text{mod}\ \ \ p,\\
\ & \ & \ & \ \\
R_{ijkm} &\equiv & - R_{jikm} &\text{mod}\ \ \ p,\\
\ & \ & \ & \ \\
R_{mijk}+R_{mjki}+R_{mkij} & \equiv  & 0&  \text{mod}\ \ \ p,\\
\ & \ & \ & \ \\
R_{ijkm} &\equiv & R_{kmij}& \text{mod}\ \ \ p.
\end{array}
\end{equation}\end{corollary}

The congruences \ref{macout} are, of course, analogous to the classical symmetries of the covariant Riemann tensor; cf.  \ref{macouttt}
in the Appendix. The first three congruences in \ref{macout} follow directly from \ref{tonton} while the fourth is a well known formal consequence of the first three.
If, in addition, we define
\begin{equation}
\Psi_{ij}:=(x^{(p^2)})^{-1}\Phi_{ij},\ \ \ \ 
R_{ik}=\Psi_{jijk},
\end{equation}
with the repeated index $j$ summed over then
one gets:

\begin{corollary}\label{jaja}
  Assume $q$ is vertical gauge invariant. Then the following congruences hold in $\cO(\widehat{G^1})$:
\begin{equation}
R_{ik}\equiv R_{ki}\ \ \ \text{mod}\ \ \ p.
\end{equation}
\end{corollary}

One can regard $R_{ik}$ as an analogue of the Ricci tensor (cf. \ref{teago} in the Appendix). 

As in the case of the classical Ricci tensor, Corollary \ref{jaja} follows by noting that  if one sets
 $$u:=(x^{(p^2)})^{-1}(q^{(p^2)})^{-1}( x^{(p^2)t})^{-1}$$
 and if one uses the last symmetry in \ref{macout}, one gets
$$R_{ik}=((x^{(p^2)})^{-1}\Phi_{ji})_{jk}=u_{jm}R_{jimk}\equiv u_{mj}R_{mkji}=R_{ki}\ \ \ \text{mod}\ \ \ p,$$
where the repeated indices $j,m$ are summed over.
\bigskip

Here is what formula \ref{tonton} gives for $n=2$ and ``conformal coordinates":

\begin{corollary}
\label{musca}
Assume $n=2$ and $q=d\cdot 1_2$.
Then for $\Phi_{12}=\varphi_{12}(x)$ we have
$$\Phi_{12}\equiv \left(\begin{array}{cc} 0 & \left(\frac{\d d}{d^p}
\right)^p\\
\ & \ \\
 - \left(\frac{\d d}{d^p}\right)^p & 0\end{array}\right)\ \ \ \text{mod}\ \ \ (p,x-1).$$
\end{corollary}

\begin{remark}
\label{tutti fruti}
Consider the situation:
$$n=2,\ \ q=d\cdot 1_2,\ \ \ \d d\not\equiv 0\ \ \text{mod}\ \  p.$$
 Then Corollary \ref{musca} implies
$$\varphi_{12}\neq  0$$
so 
if $(\d_1^G,\d_2^G)$ is the  Levi-Civita connection on $G$ attached to $(q,q)$ then 
the Frobenius lifts $\phi_1^G,\phi_2^G$  on $\widehat{G}=\widehat{GL_2}$ attached to 
$\d_1^G,\d_2^G$
do not commute. Actually the Corollary implies the stronger condition,
$$\varphi_{12}\not\equiv 0\ \ \ \text{mod}\ \ (x_{11}-x_{22},x_{12}+x_{21}),$$
which shows that already
 the Frobenius lifts $\phi_1^{G'},\phi_2^{G'}$  on $\widehat{G'}=\widehat{GL_1^c}$ attached to 
 $\d_1^{G'},\d_2^{G'}$, where the latter is
 the  Levi-Civita connection on $G'$ 
attached to $(q,q)$, do not commute; cf. Remark \ref{fruti}. 
\end{remark}

\subsection{Mixed context}

This context will involve {\it correspondences}.
We start by recalling some  terminology from 
\cite{MacLane}, p. 283; cf. also
\cite{curvature2, foundations}.

\begin{definition}
Let ${\mathcal C}$ be a category. We assume that for each two morphisms
with the same target there is a  a fiber product; for each two such morphisms we fix once and for all a fiber product  and hence a corresponding cartesian square.

1)
 By a {\it correspondence} (or a {\it span}) from an object  $X_1$ of ${\mathcal C}$ to another object $X_2$ of ${\mathcal C}$ we understand a diagram 
$\Gamma$ of schemes
$$
\begin{array}{ccccc}
\ &\ & Y &\  & \  \\
\  &  \stackrel{\pi}{\swarrow} & \  & \stackrel{\varphi}{\searrow} &\ \\
X_1 & \  & \  & \  & X_2
\end{array}
$$
We also write 
$$\Gamma=(Y,\pi,\varphi),\ \ \ Y=Y_{\Gamma}.$$
For a correspondence $\Gamma$ as above we define its {\it transpose}
by
$$\Gamma^t=(Y,\varphi,\pi).$$
If $X_1=X_2=X$ we say the correspondence above is a {\it correspondence on $X$}.
A correspondence $\Gamma=(Y,\pi,\varphi)$
 is called {\it strictly symmetric} if $\Gamma=\Gamma^t$, i.e., $X_1=X_2$ and $\pi=\varphi$.
 
 2) 
If $\Gamma=(Y_{\Gamma},\pi,\varphi)$ is a correspondence from $X_1$ to $X_2$ and  $\Gamma'=(Y_{\Gamma'},\pi',\varphi')$
is a   correspondence from $X_2$ to $X_3$ their {\it composition} is the correspondence
$$\Gamma'\circ \Gamma:=(Y_{\Gamma'\circ \Gamma},\pi\circ \pi'',\varphi'\circ \varphi'')$$
 where the above data are defined by the following diagram in which the square is the fixed cartesian square attached to the two morphisms into $X_2$:
$$\begin{array}{rclcc}
Y_{\Gamma'\circ \Gamma}
 & \stackrel{\varphi''}{\ra} & Y_{\Gamma'} & \stackrel{\varphi'}{\ra} & X_3\\
\pi'' \downarrow & \  & \downarrow \pi' & \  & \  \\
Y_{\Gamma} & \stackrel{\varphi}{\ra} & X_2 & \  & \  \\
\pi\downarrow & \  & \   & \ & \  \\
X_1 & \  & \  & \  & \  \end{array}$$

 3) A {\it  morphism}   $\Gamma'\ra \Gamma$  between correspondences $\Gamma'=(Y',\pi',\varphi')$ and $\Gamma=(Y,\pi,\varphi)$ from $X_1$ to $X_2$ is a morphism of shemes
$v:Y'\ra Y$
such that 
$$\varphi'=\varphi\circ v,\ \ \ \pi'=\pi\circ v.$$
For any $X_1$ and $X_2$ one can consider the category
${\mathcal C}(X_1,X_2)$
 whose objects are correspondences from $X_1$ to $X_2$ and whose morphisms are morphisms between such correspondences. 
As $X_1$ and $X_2$ vary the categories ${\mathcal C}(X_1,X_2)$
actually fit into a bicategory, cf. \cite{MacLane}, p. 283: the $0$-cells of the bicategory in question are the objects of ${\mathcal C}$, the $1$-cells are the correspondences, and the $2$-cells are the morphisms between correspondences. We will not need the framework of bicategories in what follows (because we will soon ``mod out" by isomorphisms); all we need is that for any correspondences $\Gamma, \Gamma',\Gamma''$ from $X_1$ to $X_2$ to $X_3$ to $X_4$ we have a natural isomorphism
$$\Gamma''\circ (\Gamma'\circ \Gamma)\simeq
(\Gamma''\circ \Gamma') \circ \Gamma.$$

 4)
 To any morphism  $u:X'\ra X$ in ${\mathcal C}$
 one can attach the correspondence 
 $$\Gamma_u:=(X',id,u)$$
 from $X'$ to $X$.
 For any further morphism $v:X''\ra X'$ we have a natural isomorphism
 $$\Gamma_{u\circ v}\simeq \Gamma_u\circ \Gamma_v.$$

  5) For any correspondence $\Gamma=(Y,\pi,\varphi)$ on $X$ and any morphism  $u:X'\ra X$
one can define the {\it pull-back} $u^*\Gamma$ of $\Gamma$ via $u:X'\ra X$ as the correspondence on $X'$ given by
$$\Gamma':=u^*\Gamma:=\Gamma^t_u\circ \Gamma \circ \Gamma_u,$$
where the first $\circ$ is performed first, say.
We have
$$\Gamma'\simeq \Gamma\times_X X':=(Y',\pi',\varphi')$$
 where
$$
Y':=(X'\times_{u,X,\pi} Y)\times_{pr_2,Y,pr_1}(Y\times_{\varphi, X, u}X'),
$$
 the subscripts indicate the maps used to construct the fiber products, $pr_1,pr_2$ are the obvious first and second projections to $Y$, and $\pi',\varphi'$ are defined by the obvious projections to $X'$. Furthermore for any further morphism $v:X''\ra X'$ we have a natural isomorphism
$$(\Gamma\times_X X')_{X'}X''\simeq \Gamma\times_X X'',\ \ \ \text{i.e.,}\ \ \ v^*u^*\Gamma\simeq (u\circ v)^*\Gamma.$$
The pull-back operation is compatible with composition in the sense that for any correspondences $\Gamma$ and $\Gamma'$ on $X$ and for any morphism $u:X'\ra X$ there is a natural morphism (which is not generally an isomorphism!),
$$(u^*\Gamma)\circ (u^*\Gamma')\ra u^*(\Gamma\circ \Gamma').$$

6) Assume that all objects and morphisms in a correspondence
$\Gamma=(Y,\pi,\varphi)$ from $X_1$ to $X_2$ are over another object $Z$ and let  $Z'\ra Z$ be a morphism in ${\mathcal C}$. Then one can define the {\it pull-back} of $\Gamma$ via $Z'\ra Z$ as the correspondence from $X_1\times_Z Z'$ to $X_2\times_Z Z'$ given by
\begin{equation}
\label{deschisaa}
(Y\times_Z Z',\pi\times_Z Z',\varphi\times_Z Z').\end{equation}
The operation of pull-back commutes (up to isomorphism, in the obvious sense) with composition. This pull-back operation is, of course, different from  the pull back of $\Gamma$ via $X'\ra X=X_1=X_2$ in 5) above.

7) Assume we are given a correspondence $\Gamma'=(Y',\pi',\varphi')$ on $X'$ and a morphism
$u:X'\ra X$. Then we define the {\it push forward} $u_*\Gamma'$ of $\Gamma'$ via $u$ as the correspondence
$$u_*\Gamma':=(Y',u\circ \pi',u\circ \varphi).$$
For any correspondence $\Gamma$ on $X$ and any morphism $u:X'\ra X$ there is a canonical morphism of correspondences on $X$,
$$u_*u^*\Gamma\ra \Gamma.$$
If $\Gamma$ is a correspondence on $X$, $\Gamma'$ is a correspondence on $X'$, and we are given a morphism $u:X'\ra X$  then by a {\it  morphism} from $\Gamma'$ to $\Gamma$ ({\it extending} $u$)
we will understand a morphism of correspondences
on $X$,
$$u_*\Gamma'\ra \Gamma.$$
The push forward operation is compatible with composition in the sense that for any correspondences $\Gamma$ and $\Gamma'$ on $X'$ and for any morphism $u:X'\ra X$ there is a natural morphism (which is not generally an isomorphism!),
$$u^*(\Gamma \circ \Gamma')\ra u_*(\Gamma)\circ u_*(\Gamma').$$
 
\end{definition}

\medskip

{\it From now on, unless otherwise explicitly stated,  we will assume ${\mathcal C}$ is the category of schemes and so correspondences will always be correspondences in the category of schemes.}

\medskip

Assume now that we are in the global situation and we are given a vertical gauge and a transversal gauge.
For $p\in {\mathcal V}$ and ${\mathfrak P}={\mathfrak P}(p)$ we denote by $X^{\widehat{\mathfrak P}}$,  $A^{\widehat{\mathfrak P}}$, the ${\mathfrak P}$-adic completions of schemes $X$ or rings $A$ over $\cO_{F,M}$ and we continue to denote by $\widehat{X},\widehat{A}$ their $p$-adic completions. Recall that ${\mathfrak P}_1={\mathfrak P}$.
Set
\begin{equation}
\label{popp}
\begin{array}{rccll}
G & = & GL_n & = & Spec\ \cO_{F,M}[x,\det(x)^{-1}],\\
\  & \  & \ & \  \\
G_{\mathfrak P} & := &  G\otimes \widehat{\cO_{\mathfrak P}} & = & Spec\ \widehat{\cO_{\mathfrak P}}[x,\det(x)^{-1}],\\
\  & \  & \ & \  \\
\widehat{G_{\mathfrak P}}
  & = & G^{\widehat{\mathfrak P}} & = &
Spf\ \widehat{\cO_{\mathfrak P}}[x,\det(x)^{-1}]^{\widehat{\ }}.
\end{array}\end{equation}
In our previous notation, of course,  $G_{\mathfrak P}=G^1$. 
Then we will prove that our Levi-Civita connection admits an {\it algebraization by correspondences} in the following sense. 

\begin{theorem}
\label{algebraization}
Let $q_1,...,q_n\in GL_n(\cO_{F,M})$, $q_i^t=q_i$. Let $(\d_1^{G_{\mathfrak P}},...,\d_n^{G_{\mathfrak P}})$ be the  Levi-Civita connection on $G_{\mathfrak P}$ over $\widehat{\cO_{\mathfrak P}}$ attached to $(q_1,...,q_n)$ and let  $(\phi_1^{G_{\mathfrak P}},...,\phi_n^{G_{\mathfrak P}})$ be the attached Frobenius lifts on $\widehat{G_{\mathfrak P}}$.
Then there exists an $n$-tuple of correspondences on $G$,
\begin{equation}
\label{YpG}
\begin{array}{ccccc}
\ &\ & Y_{p/G} &\  & \  \\
\  &  \stackrel{\pi_p}{\swarrow} & \  & \stackrel{\varphi_{pi}}{\searrow} &\ \\
G & \  & \  & \  & G
\end{array}
\end{equation}
where $i=1,...,n$,  such that the following hold:

1) The map $\pi_p:Y_{p/G}\ra G$ is \'{e}tale and $Y_{p/G}$ is affine and irreducible.

2) There is a connected component ${\mathcal Y}_{p/G}$ of  $Y_{p/G}^{\widehat{\mathfrak P}}$ such that the induced map
 $\pi_p^{\widehat{\mathfrak P}}:{\mathcal Y}_{p/G}\ra G^{\widehat{\mathfrak P}}$ is an isomorphism.

3) For each $i$ the 
following induced diagram is commutative:
\begin{equation}
\label{calYpG}
\begin{array}{ccccc}
\ &\ & {\mathcal Y}_{p/G} &\  & \  \\
\  &  \stackrel{\pi_p^{\widehat{\mathfrak P}}}{\swarrow} & \  & \stackrel{\varphi_{pi}^{\widehat{\mathfrak P}}}{\searrow} &\ \\
G^{\widehat{\mathfrak P}} & \  & \stackrel{\phi_i^{G_{\mathfrak P}}}{\longrightarrow}  & \  & G^{\widehat{\mathfrak P}}
\end{array}
\end{equation}
\end{theorem}

Intuitively the correspondences \ref{YpG} give an {\it algebraization}
of our Frobenius lifts $\phi_i^{G_{\mathfrak P}}$.
For a given $q$ the $n$-tuple of correspondences  \ref{YpG} with properties 1, 2, 3 in the theorem is, of course, far from being unique. However, for any given vertical gauge and transversal gauge, the proof of the theorem will provide, for any  $p$ and  $q$,  a {\it canonical} construction for such an $n$-tuple of correspondences \ref{YpG} on $G$.
Once we have at our disposal such a canonical $n$-tuple of correspondences on $G$ there is a general recipe to define curvature as a family of elements in the {\it ring of correspondences} on the field $E=F(x)$ of rational functions of $G$; cf. \cite{curvature2, foundations}. We quickly review in what follows this recipe; we will also add some new terminology, constructions,  and notation.

\begin{definition}\

1) Let $E$ be a field of characteristic zero and let 
$${\mathcal C}(E)={\mathcal C}(Spec\ E,Spec\ E)$$
 be the category of correspondences on $Spec\ E$.
Following the terminology and notation in \cite{foundations} we define a subcategory ${\mathcal C}_0(E)$ of ${\mathcal C}(E)$ as follows. The objects of the category ${\mathcal C}_0(E)$ are correspondences
\begin{equation}
\label{YE}
\begin{array}{ccccc}
\ &\ & Y &\  & \  \\
\  &  \stackrel{\pi}{\swarrow} & \  & \stackrel{\varphi}{\searrow} &\ \\
Spec\ E & \  & \  & \  & Spec\ E
\end{array}
\end{equation}
where $\pi$ is \'{e}tale (equivalently $Y$ is the spectrum of a finite  product of fields that are finite extensions of $E$ via $\pi$) and  $\varphi$ is a finite morphism of schemes (hence also \'{e}tale).
The finiteness of $\varphi$ is automatic if $E$ is finitely generated over ${\mathbb Q}$ which will always be the case in our applications. We say that the correspondence \ref{YE} is {\it irreducible} if and only if $Y$ is irreducible, i.e., the spectrum of a field. For any correspondence \ref{YE} in ${\mathcal C}_0(E)$ we can write $Y$ as a disjoint union of irreducible schemes $Y_i$ and the  correspondences defined by $Y_i$ will be referred to as  the {\it irreducible components} of \ref{YE}.
 A morphism in the category ${\mathcal C}_0(E)$ is, by definition, 
a morphism $v$ of correspondences  with $v$ surjective.  If the two  correspondences are irreducible then the {\it degree} of the morphism $v$ is defined as the degree of $v$ as a morphism of schemes.
A correspondence $\Gamma$ in ${\mathcal C}_0(E)$  will be called {\it categorically reduced} if 
any morphism of correspondences $\Gamma\ra \Gamma'$ in ${\mathcal C}_0(E)$ is an isomorphism.  
An irreducible correspondence $\Gamma$ in ${\mathcal C}_0(E)$ as in \ref{YE} is  categorically reduced if 
and only if
the induced morphism
$$\pi\times \varphi:Y\ra Spec\ E \times Spec\ E$$
is a closed embedding. 
A (not necessarily irreducible) correspondence in ${\mathcal C}_0(E)$  is  categorically reduced if 
and only if
all its irreducible components are categorically reduced and no two of its irreducible components are isomorphic in ${\mathcal C}_0(E)$. 
Any irreducible correspondence $\Gamma$ in ${\mathcal C}_0(E)$ has a morphism to an irreducible categorically reduced correspondence $\Gamma'$ which is uniquely determined by $\Gamma$ up to isomorphism in ${\mathcal C}_0(E)$; if $\Gamma$ is as in \ref{YE} with $Y=Spec\ L$ then one can take $\Gamma'=(Spec\ L',\pi',\varphi')$ where
$L'$ is the compositum of $\pi(E)$ and $\varphi(E)$ in $L$, equivalently
$L'$ is the image of $E\otimes_{\mathbb Q}E$ in $L$ via the natural homomorphism defined by $\pi$ and $\varphi$.

\medskip

2)
We denote by $C_+=C_+(E)$ the set of isomorphism classes of objects in the category ${\mathcal C}_0(E)$  to which we add one more element, denoted by $0$. 
Also we denote by $1$ the class of the {\it identity correspondence} (with $Y=Spec\ E$ and maps $\pi,\varphi$ equal to the identity). 
Then $C_+=C_+(E)$ comes equipped with the following operations:

\medskip

$\bullet$ transposition (coming from interchanging $\pi$ and $\varphi$);

$\bullet$ addition (coming from disjoint union of the $Y$'s);

$\bullet$ multiplication (coming from  composition of correspondences). 
\medskip

With respect to these operations $C_+(E)$ becomes a semiring with involution in the sense that

\medskip

$\bullet$ addition and multiplication are associative and addition is commutative;

$\bullet$  multiplication is left and right distributive with respect to addition;

$\bullet$  transposition is an anti-involution;

$\bullet$  $0$ is a neutral element for addition, $0^t = 0$, $0 \cdot x=0$ for all $x$;

$\bullet$  $1$ is a neutral element for composition and $1^t = 1$.

\medskip

3) The semiring $C_+$ has the additive cancellation property so it can be canonically embedded into the (associative, not necessarily commutative) ring
$C=C(E)$, 
$$C := (C_+ \times C_+)/ \sim$$ where
$$(c_1,c_2)\sim (c_3,c_4)\ \ \ \text{if and only if}\ \ \ 
c_1+c_4=c_2+c_3.$$
The ring $C(E)$ is called the {\it ring of correspondences} on $E$.
We often view $C(E)$ as a {\it Lie ring} with respect to the
 {\it commutator},
$$[c_1,c_2]:=c_1c_2-c_2c_1\in C(E),\ \ \ c_1,c_2\in C(E).$$ 
The involution $c\mapsto c^t$ on $C_+$ induces an involution $c\mapsto c^t$ on the ring $C$.
Also $C$ has a structure of ordered ring with set of positive elements the set $C_+\backslash \{0\}$.

\medskip

4)
A non-zero element of $C(E)$ is called {\it irreducible} if it is in $C_+$ and it cannot be written as a sum of two non-zero elements of $C_+$.
Of course  the class of a correspondence \ref{YE} is irreducible if and only if  the correspondence \ref{YE} is irreducible.
So any non-zero element in $C_+$ can be written uniquely as a $\bZ$-linear combination with positive coefficients of irreducible elements.

\medskip

5)   Consider the following inclusions
 $${\mathfrak S}(E)\subset {\mathfrak R}(E)\subset 
 {\mathfrak M}(E)\subset {\mathfrak L}(E)$$
 where 
 $${\mathfrak S}(E)=\text{Aut}_{\text{ring}}(E)$$ is the group of ring automorphisms of $E$,
 $${\mathfrak R}(E)=\text{End}_{\text{ring}}(E)$$ is the monoid of ring endomorphisms of $E$, 
 $${\mathfrak L}(E)=\text{End}_{\text{gr}}(E)$$ is the (not necessarily commutative) ring of additive group endomorphisms of $E$,  and 
 $${\mathfrak M}(E):=\text{End}_{\text{morita}}(E)$$
  is the ring of {\it virtual Morita endomorphisms} to be defined in what follows. (The terminology will be justified presently.) Let us 
 say that an element $\chi\in {\mathfrak L}(E)$ is a {\it Morita endomorphism}
 if there exists an integer $n\geq 1$  and a ring homomorphism $\rho:E\ra {\mathfrak g}{\mathfrak l}_n(E)$ such that
$$\chi(a)=\text{tr}(\rho(a)),\ \ \ a\in E.$$
Let us say that an element $\chi\in {\mathfrak L}(E)$ is a {\it virtual Morita endomorphism} if it is a difference in ${\mathfrak L}(E)$ of two Morita endomorphisms.
Denote by ${\mathfrak M}_+(E)$ the set of all Morita endomorphisms
and by ${\mathfrak M}(E)$ the set of virtual Morita endomorphisms. The subset ${\mathfrak M}_+(E)$ of ${\mathfrak L}(E)$ is closed under addition and multiplication and contains the unit element, i.e.,  ${\mathfrak M}_+(E)$ is a subsemiring of ${\mathfrak L}(E)$. So  ${\mathfrak M}(E)$ is a subring of ${\mathfrak L}(E)$. 

Next note that there are 
natural ring (anti)homomorphisms
\begin{equation}
\label{twoto}
\bZ{\mathfrak R}(E)\ra C(E)\ra {\mathfrak L}(E).\end{equation}
We will usually drop the prefix {\it anti} in what follows.
The first homomorphism sends a ring endomorphism $\sigma:E\ra E$ into the correspondence \ref{YE} with $\pi$ the identity and $\varphi$ induced by $\sigma$. The second homomorphism sends the class 
$$c\in C(E)$$ of a correspondence \ref{YE} into the group homomorphism
$$c^*\in {\mathfrak L}(E),\ \ \ c^*:E\stackrel{\varphi}{\longrightarrow} \cO(Y)\stackrel{\text{tr}_{\pi}}{\longrightarrow}  E$$
 where $\varphi:E\ra \cO(Y)$ is induced by the map $\varphi$ and $\text{tr}_{\pi}:\cO(Y)\ra E$ is the trace of  the map $\pi:E\ra \cO(Y)$ induced by $\pi$. The composition \ref{twoto}  is the natural map induced by the inclusion
${\mathfrak R}(E)\subset   {\mathfrak L}(E)$
and note that \ref{twoto}  is injective by the ``linear independence of characters" \cite{LangAlgebra}, p. 283.
 For 
$\sigma$ a field automorphism of $E$ the images $c$ and $c^{-1}$ of $\sigma$ and $\sigma^{-1}$ in
$C(E)$ satisfy $c^{-1}=c^t$. 
Also clearly
the 
homomorphism $C(E)\ra {\mathfrak L}(E)$ 
maps $C_+(E)$ into ${\mathfrak M}_+(E)$; hence the image of 
$C(E)\ra {\mathfrak L}(E)$
is contained in ${\mathfrak M}(E)$.

\medskip

6) An element $c\in C_+(E)$ is called {\it strictly symmetric} if it can be represented by a strictly symmetric correspondence. Two elements $c_1,c_2\in C_+(E)$ are called {\it compatible} if one can write
$$c_1c_2^t=c_3+c_4,$$
with $c_3,c_4\in C_+(E)$ and $c_3$ strictly symmetric. The relation of compatibility is, of course, symmetric, and trivially seen to be reflexive; it is not transitive in general.

\medskip

7) There are unique ring homomorphisms (the {\it left and right degree maps})
$$\text{deg}_l,\text{deg}_r:C(E)\ra \bZ$$
 given by attaching to a correspondence $(Y,\pi,\varphi)$ the positive integers  $\text{deg}(\pi)$ and $\text{deg}(\varphi)$
respectively. For $c\in C(E)$ the pair $(\text{deg}_l(c),\text{deg}_r(c))$ will be referred to a the {\it bidegree} of $c$.
\medskip

8)
Consider   the $\bZ$-linear span $J$ in $C(E)$ of all the elements of the form
$$c'-d\cdot c$$
where $c$ and $c'$ are the classes of two irreducible correspondences $\Gamma$ and $\Gamma'$ 
between which there is a morphism $\Gamma'\ra \Gamma$ of degree $d$.
One checks that $J$ is a bilateral ideal in $C(E)$.
Define the {\it  ring of categorically reduced correspondences} by 
$${\mathfrak C}(E):=C(E)/J$$
and denote by $c\mapsto \tilde{c}$ the projection $C(E)\ra {\mathfrak C}(E)$.
 Clearly, ${\mathfrak C}(E)$ has a $\bZ$-module basis consisting of   the images of the irreducible categorically reduced correspondences in  ${\mathcal C}_0(E)$; we shall refer to this basis as the {\it canonical basis} of ${\mathfrak C}(E)$. Clearly the projection $C(E)\ra {\mathfrak C}(E)$ has a distinguished section in the category of abelian groups sending any member of the canonical basis into the class of the corresponding correspondence; this section is not, however, a ring homomorphism.

The ring homomorphism $C(E)\ra {\mathfrak L}(E)$, $c \mapsto c^*$, is easily seen to factor through a homomorphism
\begin{equation}
\label{tipat3}
{\mathfrak C}(E)\ra {\mathfrak L}(E)\end{equation}
and in particular we still have an injective ring homomorphism
\begin{equation}
\label{tipat1}\bZ{\mathfrak R}(E)\ra {\mathfrak C}(E).\end{equation}
Let ${\mathfrak C}_+(E)$ be the image in ${\mathfrak C}(E)$ of $C_+(E)$.
Then the homomorphism \ref{tipat3} maps ${\mathfrak C}_+(E)$ into 
${\mathfrak M}_+(E)$ hence the image of \ref{tipat3} is contained in ${\mathfrak M}(E)$.
Also all the elements of $J$ have left and right degree $0$ so we have  induced ring homomorphisms
\begin{equation}
\label{tipat2}
\text{deg}_l,\text{deg}_r:{\mathfrak C}(E)\ra \bZ.\end{equation}
Now since $J^t=J$ the ring ${\mathfrak C}(E)$ has an  involution induced by $c\mapsto c^t$.
Define a structure of ordered ring on ${\mathfrak C}(E)$ by taking the set of positive elements to be ${\mathfrak C}_+(E)\backslash \{0\}$; of course an element of ${\mathfrak C}(E)$ is $\geq 0$ if and only if, when written as a $\bZ$-linear combination of the canonical  basis, the coefficients are $\geq 0$. 
Define the {\it irreducible} elements of ${\mathfrak C}(E)$ to be the non-zero elements $\geq 0$ that cannot be written as a sum of two  non-zero elements $\geq 0$.
Then the canonical basis of ${\mathfrak C}(E)$ consists exactly of the {\it irreducible} elements of ${\mathfrak C}(E)$. If $c\in C_+$ is strictly symmetric then its image  $\tilde{c}\in{\mathfrak C}(E)$ lies in $\bZ_{>0}$; so if $c_1,c_2\in C(E)$, $c_1,c_2>0$ are compatible then their images $\tilde{c}_1,\tilde{c}_2\in {\mathfrak C}(E)$ satisfy
$\tilde{c}_1\tilde{c}_2^t\geq 1$.
\end{definition}

\begin{remark}\ 

1)  If $F$ is a finite Galois extension of ${\mathbb Q}$ then, of course, ${\mathfrak R}(F)={\mathfrak S}(F)$. 

We claim  that, in this case,  the injective ring homomorphism \ref{tipat1},
\begin{equation}
\label{tipat5}
\bZ{\mathfrak S}(F)\ra  {\mathfrak C}(F),
\end{equation}
is also surjective, so an isomorphism.
The claim follows from the fact that, due to the normality of $F/{\mathbb Q}$, for any two field homomorphisms $\pi,\varphi:F\ra L$ we must have $\pi(F)=\varphi(F)$.

We also claim that the injective homomorphism
\begin{equation}
\bZ{\mathfrak S}(F)\ra {\mathfrak M}(F)\end{equation}
is also surjective, so an isomorphism. Indeed let $\rho:F\ra {\mathfrak g}{\mathfrak l}_n(F)$ be a ring homomorphism and let $\chi=\text{tr}\circ\rho$. Since any $a\in F$ is a root of a separable polynomial with coefficients in ${\mathbb Q}$ the same is true for $\rho(a)$, hence the minimal polynomials of $\rho(a)$ are separable, hence $\rho(a)$ are  diagonalizable in
${\mathfrak g}{\mathfrak l}_n(\overline{F})$, where $\overline{F}$ is an algebraic closure of $F$. Since the family of matrices $\{\rho(a);\ a\in F\}$ is commuting it is simultaneously diagonalizable so there exists $U\in GL_n(\overline{F})$ and maps $\lambda_1,...,\lambda_n:F\ra \overline{F}$ such that 
$$\rho(a)=U^{-1}\cdot \text{diag}(\lambda_1(a),...,\lambda_n(a))\cdot U,\ \ \ a\in F.$$
One immediately gets that the $\lambda_i$'s are ring homomorphisms. Since $F$ is Galois $\lambda_i$ come from elements of ${\mathfrak S}(F)$. But 
$$\chi(a)=\lambda_1(a)+...+\lambda_n(a),\ \ \ a\in F.$$
 So $\chi$ is in the image of $\bZ{\mathfrak S}(F)$.

So we see that  in case $F$ is a Galois number field we have natural ring homomorphisms 
$$\bZ{\mathfrak S}(F)\simeq {\mathfrak C}(F)\simeq {\mathfrak M}(F).$$ 
 Hence, for an arbitrary field $E$ of characteristic zero
either of the rings ${\mathfrak C}(E),{\mathfrak M}(E)$ could be viewed as an analogue of the group ring of the Galois group. Of these two rings 
the ring ${\mathfrak C}(E)$ has the advantage of being equipped with a natural structure of an ordered ring so it is the ring ${\mathfrak C}(E)$ that we will view as the most natural generalization of the group ring of the Galois group.

\medskip

2) Consider  again an arbitrary field $E$ of characteristic zero and let $F\subset E$ be  a subfield which is a finite Galois extension of ${\mathbb Q}$. Consider the unique additive group homomorphism
\begin{equation}
\label{ploaieuh}
\text{res}:C(E) \ra \bZ{\mathfrak S}(F)={\mathfrak C}(F)\end{equation}
sending the class of any irreducible correspondence $\Gamma=(Spec\ L,\pi,\varphi)$, where   $$\pi,\varphi:E\ra L$$ are field homomorphisms, into
the element
$$\text{deg}(\pi)\cdot \sigma,$$
 where $\sigma\in {\mathfrak S}(F)$ is the composition
$$F\stackrel{\varphi}{\longrightarrow} \varphi(F)=\pi(F)\stackrel{\pi^{-1}}{\longrightarrow} F;$$
we refer to $\sigma$ as the {\it automorphism of $F$ induced by}
$\Gamma$ (or by its class in $C(E)$).
It is trivial to check that \ref{ploaieuh} is a ring homomorphism. Then clearly the homomorphism in \ref{ploaieuh} factors through a
ring (and hence Lie ring) homomorphism
\begin{equation}
\label{ploaieoi}
\text{res}:{\mathfrak C}(E) \ra {\mathfrak C}(F)\end{equation}
On the other hand note that there is a natural embedding
 $$
{\mathfrak S}(F) \ra  {\mathfrak S}(E)$$
 sending any $\sigma:F\ra F$ into the unique automorphism (still denoted by) $\sigma:E\ra E$
 that extends $\sigma$ on $F$ and satisfies $\sigma(x)=x$. 
So we get a natural injective ring homomorphism
$$\bZ{\mathfrak S}(F)\ra \bZ{\mathfrak S}(E).$$
Composing with the natural injective ring homomorphism $\bZ{\mathfrak S}(E)\ra C(E)$ we get an injective ring homomorphism
 $${\mathfrak C}(F)\ra C(E),$$
hence a ring homomorphism 
\begin{equation}
\label{ploaiebrr}
{\mathfrak C}(F)\ra {\mathfrak C}(E).\end{equation}
The latter composed with the homomorphisms  \ref{ploaieoi} equals the identity of ${\mathfrak C}(F)$ so in particular
\ref{ploaiebrr} is still injective.
 Note that ${\mathfrak S}(F)$ acts by ring automorphisms on $C(E)$ via conjugation. Also ${\mathfrak C}(F)$, viewed as a Lie ring
 acts by derivations on $C(E)$, viewed as a Lie ring, via the commutator bracket. These actions descend to actions 
of ${\mathfrak S}(F)$ and ${\mathfrak C}(F)$ on ${\mathfrak C}(E)$.

3)  For any field $E$  of characteristic zero consider the {\it ${\mathbb Q}$-algebra of correspondences on $E$},
$$C(E)_{\mathbb Q}:=C(E)\otimes {\mathbb Q}$$
and the {\it ${\mathbb Q}$-algebra of categorically reduced correspondences on $E$},
$${\mathfrak C}(E)_{\mathbb Q}:={\mathfrak C}(E)\otimes {\mathbb Q}$$
If $E\subset E'$ is a finite field extension then the ${\mathbb Q}$-linear map
\begin{equation}
\label{anitaanita}
C(E)_{\mathbb Q}\ra C(E')_{\mathbb Q}\end{equation}
defined by
\begin{equation}
\label{anita}
\Gamma \mapsto 
\frac{1}{[E':E]}(\Gamma\times_{Spec\ E} Spec\ E').\end{equation}
induces a ring homomorphism
\begin{equation}
\label{irinasweet}
{\mathfrak C}(E)_{\mathbb Q}\ra {\mathfrak C}(E')_{\mathbb Q}.\end{equation}
So $E\mapsto {\mathfrak C}(E)_{\mathbb Q}$ defines  a functor from the category of fields  of characteristic zero and their finite field extensions
to the category of ${\mathbb Q}$-algebras. Note, by the way, that the  natural additive map \ref{anitaanita} is not a ring homomorphism; this is one more a posteriori justification for considering ${\mathfrak C}(E)$ in place of $C(E)$ and taking tensorization with ${\mathbb Q}$.

4) If $F$ is a finite Galois extension of ${\mathbb Q}$ then recall that ${\mathfrak L}(F)$
is naturally  an $F$-linear space  and ${\mathfrak S}(F)$ is an $F$-basis for ${\mathfrak L}(F)$. So we have  natural  isomorphisms of $F$-vector spaces
\begin{equation}
\label{mammon}
{\mathfrak C}(F)_F:=F\otimes_{\mathbb Q} {\mathfrak C}(F)_{\mathbb Q}
\simeq F\otimes_{\bZ} \bZ{\mathfrak S}(F)=:F{\mathfrak S}(F)\simeq 
 {\mathfrak L}(F).\end{equation}

5) Let us record the following {\it relative curvature} construction. 
Assume $F\subset E$ is a subfield which is Galois over ${\mathbb Q}$
and let 
\begin{equation}
\label{ploaieoi7}
\text{res}:{\mathfrak C}(E)_{\mathbb Q} \ra {\mathfrak C}(F)_{\mathbb Q}
\end{equation}
be the ${\mathbb Q}$-algebra map induced by \ref{ploaieoi}.
Let ${\mathfrak C}(E/F)_{\mathbb Q}$ be the kernel of the ring homomorphism \ref{ploaieoi7}. Now assume we are given 
an ${\mathbb Q}$-linear map
\begin{equation}
\label{grabesc}
c:{\mathfrak C}(F)_{\mathbb Q}\ra {\mathfrak C}(E)_{\mathbb Q},\ \ \ u\mapsto c_u
\end{equation}
 which is a section of  \ref{ploaieoi7}.
Then one  can define the {\it relative curvature}  of \ref{grabesc} as the ${\mathbb Q}$-bilinear map
\begin{equation}
\label{capra}
{\mathfrak C}(F)_{\mathbb Q}\times {\mathfrak C}(F)_{\mathbb Q}\ra {\mathfrak C}(E/F)_{\mathbb Q},\end{equation}
$$
(u,v)\mapsto [c_u,c_v]-c_{[u,v]}=c_u c_v - c_v c_u -c_{uv}+c_{vu}.$$
Note that in our applications the maps \ref{grabesc} will not be induced, in general by maps ${\mathfrak C}(F)\ra {\mathfrak C}(E)$; this is another a posteriori 
motivation for considering tensorization with ${\mathbb Q}$. 

6) There is a way to interpret (and generalize) our construction of $ {\mathfrak M}(E)$ in the framework of Hochschild homology; we will address this elsewhere.

7) We summarize the various  rings and ring homomorphisms that we have attached to a field $E$ of characteristic zero in the following diagram:
\begin{equation}
\bZ{\mathfrak S}(E)  \ra  \bZ{\mathfrak R}(E)\ra  C(E)  \ra  {\mathfrak C}(E) \ra  {\mathfrak M}(E) \ra  {\mathfrak L}(E).\end{equation}
There are also  interesting ring homomorphisms with source ${\mathfrak C}(E)$ constructed using $0$-cycles \cite{curvature2}; for simplicity we will not discuss these in the present paper.
\end{remark}

The next definition relates  general  correspondences  to our Levi-Civita context.
Assume again we are in the global situation and a vertical gauge and a transversal gauge  are given.

\begin{definition}
\label{woken}
Let
$$q\in GL_n(\cO),\ \ \ \cO=\cO_{F,M},\ \ \ q^t=q,\ \ \ q_i=\sigma_i^{-1}(q).$$
For $p\in {\mathcal V}$ set  ${\mathfrak P}:={\mathfrak P}(p)$, 
 let 
 \begin{equation}
 \label{dush}
 (\d_1^{G_{\mathfrak P}},...,\d_n^{G_{\mathfrak P}})\end{equation}
 be the  Levi-Civita connection on $G_{\mathfrak P}$ over $\widehat{\cO_{\mathfrak P}}$ attached to $(q_1,...,q_n)$ and let  $$(\phi_1^{G_{{\mathfrak P}}},...,\phi_n^{G_{{\mathfrak P}}})$$ be the attached Frobenius lifts on $\widehat{G_{{\mathfrak P}}}$. 
 (So the family of all tuples in \ref{dush} is what we referred to as the {\it adelic Levi Civita connection} attached to $(q_1,...,q_n)$.)
 Consider an $n$-tuple of correspondences \ref{YpG}  satisfying Conditions 1, 2, 3  in Theorem \ref{algebraization}. 
Let $E=F(x)$ be the field of rational functions on $G$ and let 
\begin{equation}
\label{YpE}
\begin{array}{ccccc}
\ &\ & Y_{p/E} &\  & \  \\
\  &  \stackrel{\pi_{p/E}}{\swarrow} & \  & \stackrel{\varphi_{pi/E}}{\searrow} &\ \\
Spec\ E & \  & \  & \  & Spec\ E
\end{array}
\end{equation}
be the pull back of \ref{YpG} via $Spec\ E\ra G$. We have that $Y_{p/E}$ is the spectrum of a field and  $\pi_{p/E}$ is  finite; cf. Remark \ref{ff} below.  Also $\varphi_{pi/E}$ are finite because $E$ is finitely generated over ${\mathbb Q}$. So
we can consider  the classes
\begin{equation}
\label{fifa}
c_{pi}\in C(E)\end{equation}
 of the correspondences \ref{YpE};  these classes are therefore irreducible. These classes induce classes
 \begin{equation}
\label{fifafu}
\tilde{c}_{pi}\in {\mathfrak C}(E),\end{equation}
\begin{equation}
c_{pi}^*\in {\mathfrak M}(E).\end{equation}
 Finally we can define the (mixed) {\it curvature} of the Levi-Civita connection attached to $q$ as the family of commutators
\begin{equation}
\label{uite}
\varphi_{pp'ii'}:=[c_{pi},c_{p'i'}]\in C(E),
\end{equation}
where $i,i'$ run through $1,...,n$ and $p,p'$ run through ${\mathcal V}$. This family  induces a family 
\begin{equation}
\label{uitenui}
\tilde{\varphi}_{pp'ii'}=[\tilde{c}_{pi},\tilde{c}_{p'i'}] \in {\mathfrak C}(E),
\end{equation}
in the   ring of  categorically reduced correspondences,   and finally a family
\begin{equation}
\label{uitee}\varphi^*_{pp'ii'}=[c^*_{pi},c^*_{p'i'}]\in {\mathfrak M}(E),\end{equation}
which we refer to as the (mixed) {\it $*$-curvature}. \end{definition}

Of the $3$ rings $C(E), {\mathfrak C}(E), {\mathfrak M}(E)$  the most natural choice for a    recipient ring of our curvature is, probably, ${\mathfrak C}(E)$ which can be viewed as the generalization of the group ring of a Galois group; however in order to simplify our discussion we will concentrate in what follows on curvature with values in the other two rings which are easier to analyze.

\begin{remark}\label{ff}\ 

1)  Let us check that $Y_{p/E}$ 
is the spectrum of a field and $\pi_{p/E}$ are finite. 
 Note first that, since $\pi_p$ is \'{e}tale and $G$ is an integral regular scheme,
$Y_{p/G}$ is a disjoint union of integral regular schemes. Since $Y_{p/G}$ is irreducible it is an integral scheme. Hence 
\begin{equation}
\label{cumber}
Spec\ E\times_{G,\pi_{p}} Y_{p/G}\end{equation} is the spectrum of a field $L$ which is a finite extension of $E$ via $\pi_p$.
So $Y_{p/E}$ is the spectrum of a ring of fractions of $L$; hence either $Y_{p/E}$ is empty or equal to $Spec\ L$.
So we are left with proving  that $Y_{p/E}$ is non-empty; we check this in what follows.  Denote by an upper bar tensorization with $\overline{\cO}:=\cO/{\mathfrak P}$ over $\cO$. 
Let $Y$ be an affine open subset of $Y_{p/G}$ meeting
 $\overline{{\mathcal Y}_{p/G}}$ but not meeting any other connected component of 
 $\overline{Y_{p/G}}$. By Krull's intersection theorem $\cO(X)$ and $\cO(Y)$ (being domains in which $p$ is non-invertible)  embed into $\cO(X^{\widehat{\mathfrak P}})$
and $\cO(Y^{\widehat{\mathfrak P}})$ respectively. On the other hand the map 
$$\varphi_{pi}:\cO(X^{\widehat{\mathfrak P}})\ra 
\cO(Y^{\widehat{\mathfrak P}})$$
is injective (because its reduction mod $p$ is injective and $p$ is a non-zero divisor in both rings). So the map 
$$\varphi_{pi}:\cO(X)\ra \cO(Y)$$
 is injective. So the generic point of $Y_{p/G}$ is mapped by $\varphi_{pi}$ to the generic point of $G$; this implies that
$Y_{p/E}$ is non-empty.

\medskip

2) Setting $q_i=\sigma_i^{-1}(q)$ in Definition \ref{woken}
is (a posteriori)  justified by the way the global Theorem \ref{levi Civita} and the local Theorem \ref{LC} will turn out to be related.

\medskip

3) 
As already mentioned the proof of Theorem \ref{algebraization} will provide, for a given vertical gauge and transversal gauge, a {\it canonical} construction of
 correspondences \ref{YpG} on $G$  so, for a fixed vertical gauge and transversal gauge
 there is a  {\it canonical} way to construct  correspondences \ref{YpE}
 associated to our vertical Levi-Civita connection and, in particular, the mixed curvature can be  {\it canonically attached} to $q$, the vertical gauge, and the transversal gauge.
We will not make the canonicity of \ref{YpG} explicit here; it will become clear once the proof of Theorem \ref{algebraization} will be presented. 

On the other hand let us note that if one considers two choices of correspondences \ref{YpG} satisfying the conclusions 1, 2, 3 of Theorem \ref{algebraization} and if one denotes by 
\begin{equation}
\label{butter}
c_{pi}^{(1)}\ \ \ \text{and}\ \ \ c_{pi}^{(2)}\end{equation}
the classes in \ref{fifa} corresponding to these two choices then one can show that for each $p$ and $i$ the classes \ref{butter} are compatible; cf. the argument in the proof of \cite{foundations}, Lemma 3.97. We will not use this compatibility in what follows.

\medskip

4)  In  classical Riemannian geometry  ``standard" curvature can be viewed as a family of $n\times n$ matrices
indexed by $2$ indices each of which runs through $1,...,n$. In the arithmetic
 case the curvature \ref{uite} is a family of classes correspondences
on the generic point of $GL_n$
(rather than a family of $n\times n$ matrices) indexed by $4$ indices (rather than $2$ indices); $2$ of the $4$ indices still run through $1,...,n$ while the other two run
through ${\mathcal V}$. As in the case of connections one can reduce the $4$ indices to $2$ indices 
by considering
the family
$$(\varphi_{pp'11}),$$
 indexed by $2$ indices $p,p'\in {\mathcal V}$.
   This procedure is similar to a construction that can be introduced  in classical differential geometry (cf. our Appendix).  

\medskip

5) The above constructions lead to a natural context for {\it holonomy} in the correspondence setting. Indeed, assume we are in the setting of Definition \ref{woken}. Then define 
$${\mathfrak h}{\mathfrak o}{\mathfrak l}_E\subset {\mathfrak C}(E)$$
to be the $\bZ$-linear span
of the set of iterated commutators
$$[\tilde{c}_{pi},[\tilde{c}_{p'i'},[\tilde{c}_{p''i''},[....]]]]$$
of length $\geq 2$ (i.e. involving at least $2$ elements). Then ${\mathfrak h}{\mathfrak o}{\mathfrak l}_E$ is a Lie subring of ${\mathfrak C}(E)$. We can refer to ${\mathfrak h}{\mathfrak o}{\mathfrak l}_E$ as the {\it holonomy ring} of the Levi-Civita connection attached to $q$. It is a ``correspondence version" of the holonomy ring ${\mathfrak h}{\mathfrak o}{\mathfrak l}$ introduced   in \cite{foundations}. Now let 
${\mathfrak h}{\mathfrak o}{\mathfrak l}_E^0$ and ${\mathfrak h}{\mathfrak o}{\mathfrak l}_F$
be the kernel and the image of the composition
$${\mathfrak h}{\mathfrak o}{\mathfrak l}_E\ra {\mathfrak C}(E)\stackrel{\text{res}}{\longrightarrow} {\mathfrak C}(F).$$
We get an exact sequence of Lie rings
$$0\ra {\mathfrak h}{\mathfrak o}{\mathfrak l}_E^0 \ra {\mathfrak h}{\mathfrak o}{\mathfrak l}_E\ra {\mathfrak h}{\mathfrak o}{\mathfrak l}_F\ra 0$$
and hence a Lie ring homomorphism
$${\mathfrak h}{\mathfrak o}{\mathfrak l}_F\ra \text{Out}({\mathfrak h}{\mathfrak o}{\mathfrak l}_E^0)$$
where $\text{Out}$ stands for the Lie ring of outer derivations, i.e., derivations modulo inner derivations. Such a construction was conjectured in \cite{foundations} and could be viewed as an analogue of the classical presentation of the monodromy group  of a connection as the quotient of the holonomy group by its identity component \cite{KN}.

\medskip

6) Assume our vertical gauge and transversal gauge are  perfect. Then one can define a ``relative version" of curvature as follows. 
Recall that
$${\mathfrak C}(F)_{\mathbb Q}={\mathbb Q}{\mathfrak S}(F):=\bZ {\mathfrak S}(F)\otimes{\mathbb Q},\ \ \ {\mathfrak S}(F)=\{\sigma_1,...,\sigma_n\},$$
and recall the bijection ${\mathfrak S}(F)\ra {\mathcal V}$, $\sigma\mapsto p(\sigma)$.
One can  consider the unique ${\mathbb Q}$-linear map
$$c:{\mathfrak C}(F)_{\mathbb Q}\ra {\mathfrak C}(E)_{\mathbb Q}$$
sending each $\sigma_i\in {\mathfrak S}(F)$ into  
$$c_{\sigma_i}:=\frac{\tilde{c}_{p(\sigma_i)1}}{\text{deg}_l(\tilde{c}_{p(\sigma_i)1})}.$$
Then $c$ is a section of the map \ref{ploaieoi7}. So one can consider the {\it relative curvature} map;
cf. \ref{capra}. In particular we have the {\it relative curvature} of the Levi-Civita connection attached to $q$ defined as the collection $(\varphi_{ii'})$, $\varphi_{ii'}\in {\mathfrak C}(E/F)_{\mathbb Q}$,
$$\varphi_{ii'}=
c_{\sigma_i}c_{\sigma_{i'}}-c_{\sigma_{i'}}c_{\sigma_i}-c_{\sigma_i\sigma_{i'}}+c_{\sigma_{i'}\sigma_i}.$$
If $\sigma_i$ and $\sigma_{i'}$ commute then $\varphi_{ii'}$ is a ${\mathbb Q}$-multiple of the corresponding component 
$\varphi_{p(i)p(i')11}$ of the curvature as defined in Definition \ref{woken}.
 \end{remark}

\bigskip

We turn now to the case 
 $n=2$ of the above  correspondence story and show it is ``compatible with complex structure".  Recall our notation from \ref{popp}. In analogy with \ref{popp} let us introduce, again,  $2$ variables
$\alpha,\beta$ and 
set
\begin{equation}
\label{poppp}
\begin{array}{rclll}
G' & = & GL_1^c & = & Spec\ \cO_{F,M}[\alpha,\beta,(\alpha^2+\beta^2)^{-1}],\\
\  & \  & \  & \ & \ \\
G'_{\mathfrak P} & := & G'\otimes \widehat{\cO_{\mathfrak P}}
& = & Spec\ \widehat{\cO_{\mathfrak P}}[\alpha,\beta,(\alpha^2+\beta^2)^{-1}],\\
\  & \  & \ & \ & \  \\
\widehat{G'_{\mathfrak P}}
 & = & (G')^{\widehat{\mathfrak P}} & = &
Spf\ \widehat{\cO_{\mathfrak P}}[\alpha,\beta,(\alpha^2+\beta^2)^{-1}]^{\widehat{\ }}.
\end{array}\end{equation}
We have a closed embedding 
$G'\ra G$ defined by
$$x\mapsto \left(\begin{array}{cc}\alpha & \beta\\ -\beta & \alpha\end{array}\right).$$
Then we will prove:

\begin{theorem}
\label{algebraization1}
Let 
$$d_1,d_2\in \cO_{F,M}^{\times},\ \ \ q_1=d_1\cdot 1_2,\ \ \ q_2=d_2\cdot 1_2\in GL_2(\cO_{F,M}).$$
 Let $(\d'_1,\d'_2)$ be the  Levi-Civita connection on $G'_{\mathfrak P}$ over $\widehat{\cO_{\mathfrak P}}$ attached to $(q_1,q_2)$ and let  $(\phi'_1,\phi'_2)$ be the attached Frobenius lifts on $\widehat{G'_{\mathfrak P}}$.
Then there exists a pair of correspondences on $G'$, 
\begin{equation}
\label{YpGprime}
\begin{array}{ccccc}
\ &\ & Y'_{p/G'} &\  & \  \\
\  &  \stackrel{\pi'_p}{\swarrow} & \  & \stackrel{\varphi'_{pi}}{\searrow} &\ \\
G' & \  & \  & \  & G'
\end{array}
\end{equation}
where $i=1,2$,  such that the following hold:

1) The map $\pi'_p:Y'_{p/G'}\ra G'$ is \'{e}tale and $Y'_{p/G'}$ is affine and irreducible.

2) There is a connected component ${\mathcal Y}'_{p/G'}$ of  $(Y'_{p/G'})^{\widehat{\mathfrak P}}$ such that the induced map
 $(\pi'_p)^{\widehat{\mathfrak P}}:{\mathcal Y}'_{p/G'}\ra (G')^{\widehat{\mathfrak P}}$ is an isomorphism.

3) For each $i$ the 
following induced diagram is commutative:
\begin{equation}
\label{calYpGprime}
\begin{array}{ccccc}
\ &\ & {\mathcal Y}'_{p/G'} &\  & \  \\
\  &  \stackrel{(\pi'_p)^{\widehat{\mathfrak P}}}{\swarrow} & \  & \stackrel{(\varphi'_{pi})^{\widehat{\mathfrak P}}}{\searrow} &\ \\
(G')^{\widehat{\mathfrak P}} & \  & \stackrel{\phi'_i}{\longrightarrow}  & \  & (G')^{\widehat{\mathfrak P}}
\end{array}
\end{equation}
In addition the inclusion map $G'\ra G$ lifts to a (natural) morphism from the correspondence \ref{YpGprime} to the correspondence \ref{YpG} inducing a morphism from the correspondence of formal schemes  \ref{calYpGprime} to the correspondence of formal schemes \ref{calYpG}.
\end{theorem}

\begin{definition}
Assume we are in a global situation $F,M$ and that a vertical gauge and a transversal gauge with $n=2$ are given.
 Let 
$$d\in \cO^{\times},\ \ \cO:=\cO_{F,M},\ \ q:=d\cdot 1_2\in GL_2(\cO),\ \ d_i=\sigma_i^{-1}d.$$
(This situation  could be viewed as an analogue of the case of ``conformal coordinates for metrics on surfaces" in classical differential geometry.) Fix $p\in {\mathcal V}$, set ${\mathfrak P}={\mathfrak P}(p)$, and let
$$E'=F(\alpha,\beta)$$ be the field of rational functions on $G'$. Then
(by a formal argument as in Remark \ref{ff}, 1)) the pull-backs of the correspondences \ref{YpGprime} via $Spec\ E'\ra G'$ yield irreducible correspondences $\Gamma'_{pi}$,
\begin{equation}
\label{crying}
\begin{array}{ccccc}
\ &\ & Y'_{p/E'} &\  & \  \\
\  &  \stackrel{\pi'_p}{\swarrow} & \  & \stackrel{\varphi'_{pi}}{\searrow} &\ \\
Spec\ E' & \  & \  & \  & Spec\ E'
\end{array}
\end{equation}
whose isomorphism classes are irreducible elements of $C(E')$:
\begin{equation}
\label{coralmainebuneste}
c'_{pi}\in C(E').\end{equation}
  So,  in the setting of Definition \ref{woken}, with $q=d\cdot 1_2$, $d\in \cO_{F,M}^{\times}$,  we have a well defined  {\it curvature}  $(\varphi'_{pp'ii'})$ attached to $q$ with components  
$$\varphi'_{pp'ii'}=[c'_{pi},c'_{p'i'}]\in C(E').$$ 
The latter induces 
a family of elements
$$\tilde{\varphi}'_{pp'ii'}=[\tilde{c}'_{pi},\tilde{c}'_{p'i'}]\in {\mathfrak C}(E')$$
and a family of 
group endomorphisms
$$(\varphi'_{pp'ii'})^*=[(c'_{pi})^*,(c'_{p'i'})^*]\in {\mathfrak M}(E'),$$
which we refer to as the $*$-curvature.

If our vertical gauge and transversal gauge are perfect one can introduce, again, in the obvious way, the relative curvature
$$\varphi'_{ii'}\in {\mathfrak C}(E'/F)_{\mathbb Q}.$$\end{definition}

We have the following explicit description of the correspondence \ref{crying} in case $d$ is  vertical gauge invariant. In the statement below we fix $p$ and denote by $\phi_p=\phi_{{\mathfrak P}}\in {\mathfrak S}(F)$  the Frobenius element corresponding to ${\mathfrak P}:={\mathfrak P}(p)$.

 \begin{proposition}
\label{maro}
Assume $d$ is vertical gauge invariant
and set
$$\theta_p= \frac{d^p(\alpha^2+\beta^2)^p}{\phi_p(d)(\alpha^{2p}+\beta^{2p})}\in E'.$$
Then $Y'_{p/E'}$ in \ref{crying} is isomorphic to the spectrum of a field $L_p'$ which,
viewed as an extension of $E'$ via $\pi'_p$, is  generated by a root $v_p$  of the quadratic polynomial 
\begin{equation}
\label{comparable}
2z^2+2z+1-\theta_p.\end{equation}
On the other hand, for $u_p=1+v_p$, the  homomorphisms $\varphi'_{p1}, \varphi'_{p2}$ in \ref{crying} correspond to the homomorphisms (still denoted by)
$$\varphi'_{p1},\varphi'_{p2}:E'\ra L'_p$$
 that act on $F$ via $\phi_p$  and act on $\alpha,\beta$ 
via the formulae
\medskip
$$
\varphi'_{p1} \left( \begin{array}{cc} \alpha & \beta\\- \beta& \alpha\end{array}\right)=
\left(\begin{array}{cc} \alpha^p & \beta^p\\- \beta^p& \alpha^p\end{array}\right)
\left(\begin{array}{cc} u_p & v_p\\- v_p& u_p\end{array}\right),
$$
\medskip
$$
\varphi'_{p2} \left( \begin{array}{cc} \alpha & \beta\\- \beta& \alpha\end{array}\right)=
\left(\begin{array}{cc} \alpha^p & \beta^p\\- \beta^p& \alpha^p\end{array}\right)
\left(\begin{array}{cc} u_p & - v_p\\ v_p& u_p\end{array}\right).
$$
\end{proposition}

\medskip

\begin{remark}\label{eal}
 Assume the hypotheses of Proposition \ref{maro} and $d\in \bZ[1/M]^{\times}$. 
Note that the trace of $v_p$ in the extension $E'\subset L'_p$ is the trace of the matrix 
\begin{equation}
\label{ve}
V_p=\left(\begin{array}{cc} 0 & 1\\ \frac{\theta_p-1}{2} & - 1\end{array}\right).\end{equation}
One immediately gets that
for any (not necessarily distinct) primes $p,p'$, the following formulae hold:
\begin{equation}
\label{copper}
 (\varphi'_{pp'12})^*(\alpha)=
-\text{tr} \{(\beta^{p'}+(\beta^{p'}-\alpha^{p'})V_{p'})^p\}-
\text{tr}\{(\alpha^p+(\alpha^p-\beta^p)V_p)^{p'}\},
\end{equation}
\begin{equation}
\label{copperline}
 (\varphi'_{pp'12})^*(\beta)=
\text{tr} \{(\alpha^{p'}+(\alpha^{p'}+\beta^{p'})V_{p'})^p\}+
\text{tr}\{(\alpha^p+(\alpha^p-\beta^p)V_p)^{p'}\}.
\end{equation}
Note on the other hand that the matrix $V_p$ in \ref{ve}
has entries in the ring 
\begin{equation}
\label{oce}
\bZ[1/M][\alpha,\beta,(\alpha^{2p}+\beta^{2p})^{-1}]\end{equation}
hence
the right hand sides of \ref{copper} and \ref{copperline} are elements of \ref{oce}.
Since raising matrices to power $p$ and taking trace of matrices commute modulo  $p$ we immediately get
\begin{equation}
\label{garb}
\varphi'_{pp'12}(\alpha)\equiv -2(\alpha^{pp'}+\beta^{pp'})\ \ \ \text{mod}\ \ p\end{equation}
in the ring \ref{oce}. Similarly one gets, of course,
$$\varphi'_{pp'12}(\alpha)\equiv -2(\alpha^{pp'}+\beta^{pp'})\ \ \ \text{mod}\ \ p'.$$
So we get that for all  (not necessarily distinct) $p,p'$, the $*$-curvature satisfies
$$(\varphi'_{pp'12})^*\neq 0\ \ \ \text{in}\ \ \ {\mathfrak M}(E').$$
 In particular, for all  (not necessarily distinct) $p,p'$, we have
$$\tilde{\varphi}'_{pp'12}\neq 0\ \ \ \text{in}\ \ \ {\mathfrak C}(E'),\ \ \ \text{hence}\ \ \ \varphi'_{pp'12}\neq 0\ \ \ \text{in}\ \ \ C(E').$$
It is interesting to note a contrast between   formula \ref{garb} and the formula in Corollary \ref{musca}.
The former tells the ``correspondence story" (that involves the trace) whereas the latter is the $p$-adic story (that does not involve the trace);
the two stories turn out to be different even in this simple example. In particular the reduction mod $p$
of $(\varphi_{pp12}')^*$ in \ref{garb} does not depend on $d$ whereas the reduction mod $p$ of $\varphi_{12}$ in Corollary \ref{musca} {\it does depend} on $d$. Note, on the other hand that $(\varphi_{pp12}')^*$ itself (not reduced mod $p$) {\it still depends} on $d$.\end{remark}

Assume in what follows the hypotheses of Proposition \ref{maro} and consider the subfield 
\begin{equation}
\label{subfieldd}
E''':=F(t)\subset E'=F(\alpha,\beta),\ \ \ t=\alpha/\beta.\end{equation}
Note that, in case $\sqrt{-1}\in F$, the extension \ref{subfieldd}
is induced by the group homomorphism
$$\text{det}^{\perp}:G'=GL_1^c=Spec\ \cO[\alpha,\beta,(\alpha^2+\beta^2)^{-1}]\ra G'''
= GL_1=Spec\ \cO[z,z^{-1}]$$
given by 
$$z\mapsto s:=\frac{\alpha+\sqrt{-1}\beta}{\alpha-\sqrt{-1}\beta}=\frac{t+\sqrt{-1}}{t-\sqrt{-1}}.$$
Recalling the correspondences $\Gamma'_{pi}$ in \ref{crying} we have:

\begin{proposition}
\label{galben}
 There exist categorically reduced correspondences $\Gamma'''_{pi}$ in ${\mathcal C}_0(E''')$,
 \begin{equation}
 \label{magrab}
 \begin{array}{ccccc}
\ &\ & Y'''_{p/E'''} &\  & \  \\
\  &  \stackrel{\pi'''_p}{\swarrow} & \  & \stackrel{\varphi'''_{pi}}{\searrow} &\ \\
Spec\ E''' & \  & \  & \  & Spec\ E'''
\end{array}
 \end{equation}
 equipped with morphisms from 
 $\Gamma'_{pi}$ to $\Gamma'''_{pi}$,
  extending the inclusion \ref{subfieldd}.
 Moreover, up to isomorphisms in ${\mathcal C}_0(E''')$,
the categorically reduced correspondences $\Gamma'''_{pi}$ and  the morphisms  from $\Gamma'_{pi}$ to $\Gamma'''_{pi}$  
are unique.\end{proposition}

   The correspondences $\Gamma'''_{pi}$  define canonical classes 
\begin{equation}
\label{lindalinda}c'''_{pi}\in C(E'''),\ \ \ \tilde{c}_{pi}'''
\in {\mathfrak C}(E'''),\ \ \ c_{pi}'''^{*}
\in {\mathfrak M}(E''');\end{equation}
one can then consider the commutators of these elements in the corresponding rings as defining {\it curvatures} in these rings; in the case of a perfect vertical gauge and a perfect tansversal one can also introduce, as before, the {\it relative curvature} with values in
$ {\mathfrak C}(E'''/F)_{\mathbb Q}$.
 We will not pursue this here.

\medskip

\section{Proofs of the main results}

\subsection{Levi-Civita connections}
We begin by proving our results about the existence and uniqueness of Levi-Civita connections.

\bigskip

{\it Proof of Theorem \ref{LC}}.

We first prove  the existence of the tuple $(\d^G_1,...,\d^G_n)$.  
The argument is an extension of the argument in  the proof of Theorem 4.38 in \cite{foundations}.

Consider the matrices
$$A_i=x^{(p)t}\cdot \phi(q_i)\cdot x^{(p)},\ \ B_i=(x^tq_ix)^{(p)}.$$
Note that for any $n$-tuple $(\d^{G}_1,...,\d^{G}_n)$ of $p$-adic connections on $G$, if 
$$\Delta_i:=\d_i^G x,\ \ \ \Gamma_i:=\Delta_i^t \cdot \phi(q_i)\cdot x^{(p)}$$
 then
\begin{equation}
\label{xxo}
\Gamma_i
=\Delta^t_i\cdot (x^{(p)t})^{-1}\cdot A_i.\end{equation}
We will construct by induction a sequence of $n$-tuples
\begin{equation}
\label{lamb}
(\Lambda_1^{(\nu)},...,\Lambda_n^{(\nu)}),\ \ \nu\geq 1\end{equation}
of $n\times n$ matrices with entries in $\cO(\widehat{G})$
such that if 
$$\Delta_i^{(\nu)}:=\frac{1}{p}x^{(p)}(\Lambda_i^{(\nu)}-1),\ \ \Gamma_i^{(\nu)}:=\Delta^{(\nu)t}_i\cdot (x^{(p)t})^{-1}\cdot A_i,\ \ \Gamma_i^{(\nu)}=(\Gamma_{ijk}^{(\nu)}),$$
then the following properties hold:

\medskip

i) $\Lambda_i^{(1)}=1$,

\medskip

ii) $\Lambda_i^{(\nu+1)}\equiv \Lambda_i^{(\nu)}$ mod $p^{\nu}$,

\medskip

iii) $\Lambda_i^{(\nu)t}A_i\Lambda_i^{(\nu)}\equiv B_i$ mod $p^{\nu}$,

\medskip

iv) $\Gamma^{(\nu)}_{ijk}= \Gamma^{(\nu)}_{jik}$; equivalently $(A_i(\Lambda_i^{(\nu)}-1))_{kj}=(A_j(\Lambda_j^{(\nu)}-1))_{ki}$.

\medskip

\noindent We claim this  ends the proof of the existence
of the tuple $(\d^G_1,...,\d^G_n)$ in our theorem. Indeed one can then set
$$\Lambda_i:=\lim_{\nu\ra \infty} \Lambda_i^{(\nu)}$$
and one can define $\phi_i^G$ by setting
$$\phi^{G}_i(x):=x^{(p)}\Lambda_i.$$
By the way, with these definitions if
$$\Delta_i:=\lim_{\nu\ra \infty} \Delta_i^{(\nu)}$$
then $$\d^G_i x=\Delta_i,\ \ \ \Gamma_i=\lim_{\nu\ra \infty} \Gamma_i^{(\nu)}.$$
Now Condition iii) above implies 
$$\Lambda_i^t A_i \Lambda_i=B_i$$
 which is equivalent to
assertion 1 of the Theorem; also Condition iv) above
implies  assertion 2 of the Theorem, which ends our proof.

To construct of our sequence of $n$-tuples \ref{lamb} define the $n$-tuple for $\nu=1$ by Condition i), assume the $n$-tuple \ref{lamb}
was constructed for some $\nu$ and seek the $n$-tuple \ref{lamb}
corresponding to $\nu+1$ in the form
\begin{equation}
\label{alpha}
\Lambda_i^{(\nu+1)}=\Lambda_i^{(\nu)}+p^{\nu}Z_i.\end{equation}
Write 
\begin{equation}
\label{betabeta}
\Lambda_i^{(\nu)t}A_i\Lambda_i^{(\nu)}= B_i-p^{\nu}C^{(\nu)}_i.\end{equation}
Then
\begin{equation}
\label{cip cirip}
\begin{array}{rcl}
\Lambda_i^{(\nu+1)t}A_i\Lambda_i^{(\nu+1)} & \equiv & \Lambda_i^{(\nu)t}A_i\Lambda_i^{(\nu)}\\
\  & \  & \  \\
\  & \  & +p^{\nu}(\Lambda_i^{(\nu)t}A_iZ_i+Z_i^tA_i\Lambda_i^{(\nu)})\ \ \text{mod}\ \ p^{\nu+1}\\
\  & \  & \  \\
\  & \equiv & B_i+p^{\nu}(-C^{(\nu)}_i+A_iZ_i+Z_i^tA_i) \ \ \text{mod}\ \ p^{\nu+1}.
\end{array}\end{equation}
Now 
$A_i^t=A_i$ and $B_i^t=B_i$ so 
$C_i^{(\nu)t}=C^{(\nu)}_i.$
 So if 
$C^{(\nu)}_i=(C^{(\nu)}_{ijk})$
 then $C^{(\nu)}_{ijk}=C^{(\nu)}_{ikj}$.
Define
\begin{equation}
\label{gamma}
D^{(\nu)}_{ijk}:=\frac{1}{2}(C^{(\nu)}_{ijk}+C^{(\nu)}_{jik}-C^{(\nu)}_{kij}).\end{equation} 
Then
\begin{equation}
\label{66}
D^{(\nu)}_{ijk}=D^{(\nu)}_{jik}\end{equation}
 and
\begin{equation}
\label{67}
D^{(\nu)}_{ijk}+D^{(\nu)}_{ikj}=C^{(\nu)}_{ijk}.\end{equation}
So if we define  the matrices $$D^{(\nu)}_i=(D^{(\nu)}_{ijk})$$
 we have:
$$D^{(\nu)}_i+D_i^{(\nu)t}=C^{(\nu)}_i.$$
Setting 
\begin{equation}
\label{omega}
Z_i:=A_i^{-1}D^{(\nu)t}_i\end{equation}
 we get 
$$D^{(\nu)t}_i=A_iZ_i,\ \ \ D^{(\nu)}_i=Z_i^tA_i.$$
 So, by \ref{cip cirip},  
 $$\Lambda_i^{(\nu+1)t}A_i\Lambda_i^{(\nu+1)}\equiv B_i\ \ \ \text{mod}\ \ p^{\nu+1},$$
 hence Condition iii) holds for $\Lambda_i^{(\nu+1)}$ mod $p^{\nu+1}$.

 To check Condition iv) for $\nu$ replaced by $\nu+1$ note that
 $$
 \begin{array}{rcl}
(A_i(\Lambda_i^{(\nu+1)}-1))_{kj} & = & (A_i(\Lambda_i^{(\nu)}+p^{\nu}Z_i-1))_{kj}\\
\  & \  & \  \\
\  & = & (A_i(\Lambda^{(\nu)}_i-1))_{kj}+p^{\nu}D_{ijk}^{(\nu)}\\
\  & \  & \  \\
\  & = & (A_j(\Lambda^{(\nu)}_j-1))_{ki}+p^{\nu}D_{jik}^{(\nu)}\\
\  & \  & \  \\
\  & = & (A_j(\Lambda_j^{(\nu+1)}-1))_{ki}.
\end{array}
$$
This ends the proof of the existence part of our Theorem.
 
 We next prove the uniqueness of the tuple  $$(\d^G_1,...,\d^G_n)$$
  in our theorem. 
  
  Assume we have two such tuples which we denote by $$(\d_1,...,\d_n)\ \ \ \text{and}\ \ \  (\d'_1,...,\d_n').$$
 Let $\phi_i$ and $\phi'_i$ be the corresponding Frobenis lifts on $\widehat{G}$, write 
  $$\phi_i(x)=x^{(p)}\Lambda_i,\ \ \ \ \phi'_i(x)=x^{(p)}\Lambda'_i$$
   for matrices $\Lambda_i,\Lambda_i'$, and let 
   $$\d_ix=\Delta_i,\ \ \ \d'_ix=\Delta'_i,\ \ \ \Gamma_i=\Delta_i^t \cdot (x^{(p)t})^{-1}\cdot A_i,\ \ \ \Gamma'_i=(\Delta_i')^t \cdot (x^{(p)t})^{-1}\cdot A_i.$$ We have 
  \begin{equation}
  \label{fish}
  \Lambda_i^tA_i\Lambda_i=B_i,\ \ \ (\Lambda_i')^tA_i\Lambda'_i=B_i\end{equation}
  and 
  \begin{equation}
  \label{clouds}
  \Gamma_{ijk}=\Gamma_{jik},\ \ \ \Gamma'_{ijk}=\Gamma'_{jik}.\end{equation}
  We will prove that 
  \begin{equation}
  \label{ciri}
  \Lambda_i\equiv \Lambda'_i\ \ \ \text{mod}\ \ p^{\nu}
  \end{equation}
   by induction on $\nu$ and this will end the proof. The case $\nu=1$ is clear. Assume \ref{ciri} holds for some $\nu\geq 1$ and write
   \begin{equation}
   \label{toctoc}
   \Lambda'_i=\Lambda_i+p^{\nu}Z_i.\end{equation}
   From \ref{fish} we get
   $$B_i\equiv B_i+p^{\nu}A_iZ_i+p^{\nu}Z_i^tA_i\ \ \ \text{mod}\ \ \ p^{\nu+1},
   $$
   hence, setting
   $$E_i=Z_i^tA_i=(E_{ijk})$$
   we get
   $$E_i+E_i^t\equiv 0\ \ \ \text{mod}\ \ \ p,$$
   hence
   \begin{equation}
   \label{tuctuc}
   E_{ijk}\equiv - E_{ikj}\ \ \ \text{mod}\ \ \ p.
   \end{equation}
   On the other hand, from \ref{toctoc} we get
   $$\Gamma'_i=\Gamma_i+p^{\nu-1}E_i$$
   hence, by \ref{clouds}, 
   \begin{equation}
   \label{tctc}
   E_{ijk}=E_{jik}.
   \end{equation}
   Combining \ref{tuctuc} and \ref{tctc} we get
   \begin{equation}
   \label{tirgo}
   E_{ijk}\equiv E_{jik}\equiv - E_{jki}\ \ \ \text{mod}\ \ \ p.\end{equation}
   Applying \ref{tirgo} three times we get
   $$E_{ijk}\equiv - E_{jki}\equiv E_{kij}\equiv - E_{ijk}\ \ \ \text{mod}\ \ p,$$
 hence $$2E_{ijk}\equiv 0\ \ \ \text{mod}\ \ \ p.$$
 Since $p\neq 2$ we get $$E_{ijk}\equiv 0\ \ \ \text{mod}\ \ \ p,$$
hence
  $$Z_i\equiv 0\ \ \ \text{mod}\ \ \ p,$$ hence 
   $$\Lambda'_i\equiv \Lambda_i\ \ \ \text{mod}\ \ p^{\nu+1},$$
   and our induction step is proved.
  \qed

  \bigskip
  
  {\it Proof of Proposition \ref{congruences}}. Assume the notation in the proof of Theorem \ref{LC}.
  Then  $C_i$ in the statement of Proposition \ref{congruences} equals $C_i^{(1)}$ in the proof 
  of Theorem \ref{LC}. The proof of Theorem \ref{LC} shows that
\begin{equation}
\label{inter}
\Gamma_i^{(\nu+1)}=\Gamma_i^{(\nu)}+p^{\nu-1}D_i^{(\nu)}.
\end{equation}
  Since $\Lambda_i^{(1)}=1$ we have $\Gamma_i^{(1)}=0$ hence, by \ref{inter},
  $$\Gamma_{ijk}\equiv \Gamma_{ijk}^{(2)}\equiv D_{ijk}^{(1)}\ \ \ \text{mod}\ \ \ p$$
  and we are done by \ref{gamma}.
  
   \qed

\bigskip

{\it Proof of Proposition \ref{unichris}}. 
Define the following matrices
$$A^{\text{univ}}_i=x^{(p)t}\cdot (s_i^{(p)}+ps'_i)\cdot x^{(p)},\ \ \ B^{\text{univ}}_i=(x^ts_ix)^{(p)}\ \ i=1,...,n$$
with entries in the ring \ref{ncis}. 
To conclude the proof of the Proposition it is enough to find $n\times n$ matrices
$$\Lambda^{\text{univ}}_i,\ \ \ \ i=1,...,n$$
with entries in the ring \ref{ncis} and $\equiv 1$ mod $p$ such that if one sets
$$\Delta_i^{\text{univ}}:=\frac{1}{p}x^{(p)}(\Lambda_i^{\text{univ}}-1),\ \ \Gamma_i^{\text{univ}}:=(\Delta^{\text{univ}}_i)^t\cdot (x^{(p)t})^{-1}\cdot A_i^{\text{univ}}$$
then the following conditions are satisfied:

\medskip

1) $(\Lambda^{\text{univ}}_i)^t \cdot A^{\text{univ}}_i \cdot \Lambda^{\text{univ}}_i=B^{\text{univ}}_i$;

\medskip

2) $\Gamma^{\text{univ}}_{ijk}=\Gamma^{\text{univ}}_{jik}$.

\medskip

\noindent The existence of these matrices can be proved by redoing the existence part of the proof of Theorem \ref{LC} with the ring $\cO[x,\det(x)^{-1}]^{\widehat{\ }}$ replaced by the ring \ref{ncis} and with 
$q_i,\phi(q_i)$ replaced by $s_i,s_i^{(p)}+ps'_i$.
\qed

\bigskip

{\it Proof of Theorem \ref{levi Civita}}.
Consider any index $i$.
Consider the  bijection 
$$\d^i\mapsto \d_i:=\d^1_i:=\sigma_i^{-1}\circ \d^i \circ \sigma_i$$
between  $p$-adic connections 
on $G^i$ and $p$-adic connections on $G^1$.
The Frobenius lifts $\phi^i$ and $\phi_i^1$ attached to 
$\d^i$ and $\d_i^1$ are then related by
$$\phi^1_i:=\sigma_i^{-1}\circ \phi^i \circ \sigma_i.$$
Consider, on the other hand, the matrices
$$\alpha_i=x^{(p)t}\cdot \phi^i(q)\cdot x^{(p)},\ \ \beta=(x^tq x)^{(p)}$$
with entries in $\cO(\widehat{G^i})$. 
Set $$\phi^i(x)=x^{(p)}\lambda_i.$$ 
The condition that $\phi^i$ be ${\mathcal H}_q$-horizontal with respect to the trivial connection $\phi_0^i$ on $G^i$ is equivalent to
\begin{equation}
\label{mall1}
\lambda_i^t \alpha_i \lambda_i=\beta.
\end{equation}
Consider the matrices 
$$A_i=\sigma_i^{-1}\alpha_i,\ \ B_i=\sigma_i^{-1}\beta,\ \ \Lambda_i=\sigma_i^{-1}\lambda_i$$
with entries in $\cO(\widehat{G^1})$ and set
$$q_i=\sigma_i^{-1}(q).$$
 So \ref{mall1} is equivalent to
\begin{equation}
\label{mall2}
\Lambda_i^t A_i\Lambda_i=B_i.
\end{equation}
On the other hand we have
$$
\begin{array}{rcl}
A_i & = &  x^{(p)t} \cdot \phi_i^1(\sigma_i^{-1}(q)) \cdot x^{(p)}= x^{(p)t} \cdot \phi_i^1(q_i) \cdot x^{(p)};\\
\  & \  & \  \\
B_i & = & (\sigma_i^{-1}(x^tqx))^{(p)}= (x^t \sigma_i^{-1}(q) x)^{(p)}= (x^t q_i x)^{(p)};\\
\  & \  & \  \\
\phi^{G^1}_i(x) & = & (\sigma_i^{-1}\circ \phi^i \circ \sigma_i) (x)= \sigma_i^{-1}(\phi^i(x))= \sigma_i^{-1}(x^{(p)}\lambda_i)= x^{(p)}\Lambda_i.$$
\end{array}
$$
So \ref{mall2} is equivalent to $\phi^1_i$ being $\cH_{q_i}$-horizontal with respect to $\phi^1_0$; hence the latter condition 
is equivalent to the condition that $\phi^i$ 
be $\cH_{q}$-horizontal with respect to $\phi^i_0$.

To tackle torsion freeness consider the matrices
$$\gamma_i:=\d^i x^t \cdot \phi^i(q) \cdot x^{(p)},\ \ \Gamma_i:=\sigma_i^{-1}\gamma_i.$$
Note that
$$\Gamma_{ijk}=\sigma_i^{-1}\gamma_{ijk}$$
hence $(\d^1,...,\d^n)$ is torsion free relative to $q$ if and only if $\Gamma_{ijk}=\Gamma_{jik}$.
But on the other hand we have
$$\begin{array}{rcl}
\Gamma_i & = & \sigma_i^{-1}(\d^i(x^t)) \cdot
 \sigma_i^{-1}(\phi^i(q)) \cdot x^{(p)}\\
\  & \  & \  \\
\  & = & \d^1_i( \sigma_i^{-1} (x^t)) \cdot \phi^1(\sigma_i^{-1}(q)) \cdot x^{(p)}\\
\  & \  & \  \\
\ & = & \d^1_i x^t \cdot \phi^1(q_i) \cdot x^{(p)},
\end{array}
$$
so $\Gamma_i$ is the Christoffel symbol of $\d^1_i$ relative to $q_i$.
So $(\d^1,...,\d^n)$ is torsion free relative to $q$ if and only if $(\d^1_1,...,\d^1_n)$ is torsion free relative to $(q_1,...,q_n)$.
At this point it is clear that Theorem \ref{levi Civita} follows from Theorem
\ref{LC} applied to $\cO=\widehat{\cO_{{\mathfrak P}_1}}$ and $G=G^1$.\qed

\bigskip

{\it Proof of Proposition \ref{voices}}.
Assume the notation in the proof of Theorem \ref{levi Civita}.
Note that
$$
\begin{array}{rcl}
C_i & = & - x^{(p)t}\cdot \sigma_i^{-1}\d^i q \cdot x^{(p)} +\frac{1}{p}\{(x^t \cdot \sigma_i^{-1}q \cdot x)^{(p)}-x^{(p)t}\cdot \sigma_i^{-1}q^{(p)}\cdot x^{(p)}\}\\
\  & \  & \  \\
\  & = &  - x^{(p)t}\cdot \d_i^1 \sigma_i^{-1} q \cdot x^{(p)} +\frac{1}{p}\{(x^t \cdot \sigma_i^{-1}q \cdot x)^{(p)}-x^{(p)t}\cdot \sigma_i^{-1}q^{(p)}\cdot x^{(p)}\}\\
\  &  \  & \  \\
\  & = & - x^{(p)t}\cdot \d^1 q_i \cdot x^{(p)} +\frac{1}{p}\{(x^t \cdot q_i \cdot x)^{(p)}-x^{(p)t}\cdot q_i^{(p)}\cdot x^{(p)}\}.
\end{array}$$
So Proposition \ref{voices} follows directly from Proposition \ref{congruences} applied
to $\cO=\widehat{\cO_{{\mathfrak P}_1}}$ and $G=G^1$.
\qed

\bigskip

{\it Proof of assertion 1 in Theorem \ref{mock}}.
With the notation in the proof of Theorem \ref{LC}  view 
$$A_i,\ B_i,\  \Lambda_i$$
 as matrices with coefficients in $\widehat{\mathcal B}=\cO(\widehat{G})$, where
$$G=GL_2=Spec\ {\mathcal B},\ \ \ {\mathcal B}=\cO[x,\det(x)^{-1}].$$
Let $\alpha,\beta$ be two variables and recall that we view 
$$GL_1^c=Spec\ {\mathcal B}',\ \ \ {\mathcal B}'=\cO[\alpha,\beta,(\alpha^2+\beta^2)^{-1}]$$
embedded into 
$G$
via the map
\begin{equation}
\label{circus}
{\mathcal B}\ra {\mathcal B}',\ \ \ b\mapsto b',\end{equation}
defined by
$$x\mapsto x':= \left(\begin{array}{cc}\alpha & \beta\\- \beta & \alpha\end{array}\right).$$
Recall that we denoted by ${\mathfrak g}{\mathfrak l}_2$  the  functor that attaches to any ring the algebra of $2\times 2$ matrices with coefficients in that ring;   
let ${\mathfrak g}{\mathfrak l}_1^c$ the functor that attaches to any ring the commutator of $c$ in  ${\mathfrak g}{\mathfrak l}_2$ applied to that ring.  We still denote by
$${\mathfrak g}{\mathfrak l}_2({\mathcal B})\ra {\mathfrak g}{\mathfrak l}_2({\mathcal B}'),\ \ \ M\mapsto M'$$
the map induced by \ref{circus}.
It is enough to prove that for $i=1,2$, 
$$\Lambda_i'\in GL_1^c(\widehat{{\mathcal B}'}).$$
It is then enough to prove, by induction on $\nu$, that 
\begin{equation}
\label{sasolito}
(\Lambda^{(\nu)}_i)'\in GL_1^c(\widehat{{\mathcal B}'}).\end{equation}
This is clearly true for $\nu=1$. Assume \ref{sasolito} for some $\nu$. Now clearly 
$$A_i',\ B_i' \in GL_1^c(\widehat{{\mathcal B}'}).$$
By \ref{betabeta} we get
$$C_i':=(C^{(\nu)}_i)' \in {\mathfrak g}{\mathfrak l}_1^c(\widehat{{\mathcal B}'}).$$
Set 
$$D'_{ijk}:=\frac{1}{2}(C'_{ijk}+C'_{jik}-C'_{kij}).$$
We have
$$D'_{i11}=\frac{1}{2}C'_{i11}=\frac{1}{2}C'_{i22}=D'_{i22}$$
and
$$
D'_{i12}  =  \frac{1}{2}(C'_{i12}+C'_{1i2}-C'_{2i1})
=   \frac{1}{2}(- C'_{i21}+C'_{1i2}-C'_{2i1})
 =  - D'_{i21}.$$
So, with notation as in \ref{gamma},
$$D_i'=(D^{(\nu)}_i)'\in {\mathfrak g}{\mathfrak l}_1^c(\widehat{{\mathcal B}'}).$$
Hence, with the notation in \ref{omega} we have
$$Z_i'\in {\mathfrak g}{\mathfrak l}_1^c(\widehat{{\mathcal B}'})$$
and hence, by \ref{alpha},
$$(\Lambda^{(\nu+1)}_i)'\in {\mathfrak g}{\mathfrak l}_1^c(\widehat{{\mathcal B}'}),$$
which ends  our induction.
\qed

\bigskip

{\it Proof of Proposition \ref{fishatnoon}}.
With the notation in the proof of Theorem \ref{LC} set
 $\Phi_i=x^{(p)}\Lambda_i$ and view 
$$\Phi_i,\ A_i,\ B_i,\ \Lambda_i$$ as matrices with coefficients in $\cO(\widehat{G})$.
Let 
$$\Phi'_i, \ A'_i, \ B'_i, \ \Lambda'_i$$
 be the images of the corresponding matrices in the ring of matrices with coefficients in $\cO(\widehat{G'})$.
We then have
$$
A_i' =\phi(d_i)\cdot (\alpha^{2p}+\beta^{2p})\cdot 1_2,\ \ \ B'_i=d_i^p\cdot 
(\alpha^2+\beta^2)^p\cdot 
1_2$$
and recall that we defined
\begin{equation}
\label{russianrussian}
\theta_i:= \frac{d_i^p(\alpha^2+\beta^2)^p}{\phi(d_i)(\alpha^{2p}+\beta^{2p})}.\end{equation}
Now, by the proof of Theorem \ref{LC} and by Remark \ref{ozi},  the following hold  in $\cO(\widehat{G})$
for all $i,j,k\in \{1,2\}$:

\medskip 

1) $\Lambda_i\equiv 1_2$ mod $p$;

\medskip

2) $\Lambda_i^tA_i\Lambda_i= B_i$;

\medskip

3) $\phi(d_i)(\Lambda_i-1_2)_{kj}=\phi(d_j)(\Lambda_j-1_2)_{ki}$.

\medskip

\noindent We get that the following hold in $\cO(\widehat{G'})$:

\medskip 

1') $\Lambda'_i\equiv 1_2$ mod $p$;

\medskip

2') $(\Lambda'_i)^tA'_i\Lambda'_i= B'_i$;

\medskip

3') $\phi(d_i)(\Lambda'_i-1_2)_{kj}=\phi(d_j)(\Lambda'_j-1_2)_{ki}$.

\medskip

\noindent Of course 3')  only needs to be checked for $(i,j)=(1,2)$.
Now an argument similar to the one proving uniqueness in Theorem \ref{LC} shows that the Conditions 1', 2', 3' uniquely determine $$\Lambda_1',\Lambda'_2
\in GL_2(\cO(\widehat{G'})).$$
 So in order to conclude our proof  it is enough to show that 1', 2', 3' hold if one replaces $\Lambda'_1$ and $\Lambda'_2$ by the matrices 
\begin{equation}
\label{xox}
\left(\begin{array}{cc}
u_1 & v_1\\ -v_1 & u_1\end{array}\right) \ \ \ 
\text{and}\ \ \  \left(\begin{array}{cc}
u_2 & v_2\\ -v_2 & u_2\end{array}\right),\end{equation}
respectively, where $u_i,v_i$ are as in \ref{as} and \ref{bs} respectively.
The Conditions 1', 2', 3' for the matrices \ref{xox} translate into the following conditions:

\medskip

1') $u_i  \equiv  1$ and  $v_i  \equiv  0$ mod $p$,

\medskip

2') $u_i^2+v_i^2 =\theta_i$,

\medskip

3') $\phi(d_1)v_1=\phi(d_2)(u_2-1)$ and $\phi(d_1)(u_1-1)= - \phi(d_2) v_2$.

\medskip

Checking 1', 2', 3' is a trivial exercise left to the reader.
\qed

\begin{remark}
The argument in the proof of Proposition \ref{fishatnoon} can be used to  give an alternative proof of (the already proved) assertion 1 in Theorem \ref{mock}.
\end{remark}

\bigskip

{\it Proof of Proposition \ref{tirg}}. We adopt the notation in that Proposition and in Remark \ref{fishallday} preceding it.

Assertion 1 follows immediately by taking determinants in \ref{fid}. 

To check assertion 2 we may identify $z$ and $s$. Note that 
$$\theta=\frac{d^p}{\phi(d)}\cdot \frac{((s-1)^2-(s+1)^2)^p}{(s-1)^{2p}-(s+1)^{2p}}
\in 1+p\cO[s,s^{-1}]^{\widehat{\ }},$$
because 
$$(s-1)^{2p}-(s+1)^{2p}\equiv -4s^p\ \ \ \text{mod}\ \ p,$$
hence 
$$u,v\in \cO[s,s^{-1}]^{\widehat{\ }},$$
so 
$$\frac{u+\sqrt{-1}v}{u-\sqrt{-1}v}=\frac{(u+\sqrt{-1}v)^2}{\theta}
\in \cO[s,s^{-1}]^{\widehat{\ }}.$$
Also we have
$$\frac{\alpha^p+\sqrt{-1}\beta^p}{\alpha^p-\sqrt{-1}\beta^p} =
\frac{(\sqrt{-1})^{p-1}(s+1)^p+(s-1)^p}{(\sqrt{-1})^{p-1}(s+1)^p-(s-1)^p}\in \cO[s,s^{-1}]^{\widehat{\ }},$$
because
$$
\pm (s+1)^p-(s-1)^p\equiv 2\ \ \text{or}\ \ -2s^p\ \ \text{mod}\ \ p.
$$
On the other hand one immediately checks that
$$\phi^{G'}_1\left(\frac{\alpha+\sqrt{-1}\beta}{\alpha-\sqrt{-1}\beta}\right)=\frac{\alpha^p+\sqrt{-1}\beta^p}{\alpha^p-\sqrt{-1}\beta^p} \cdot \frac{u+\sqrt{-1}v}{u-\sqrt{-1}v},$$
hence
$$\phi^{G'}_1(s)\in \cO[s,s^{-1}]^{\widehat{\ }}.$$
Similarly one shows
$$\phi^{G'}_2(s)\in \cO[s,s^{-1}]^{\widehat{\ }},$$
which ends the proof of assertion 2.

To check assertion 3 note that,
putting the diagrams \ref{separatist} and \ref{separatistule}
together we get  commutative diagrams
\begin{equation}
\label{separatistulele}
\begin{array}{rcl}
\widehat{G'} & \stackrel{\phi^{G'}_i}{\longrightarrow} & \widehat{G'}\\
\text{det}\times \text{det}^{\perp} \downarrow & \  & \downarrow \text{det}\times \text{det}^{\perp}\\
\widehat{G'''} \times \widehat{G}'''& \stackrel{\phi^{G'''}\times \phi_i^{G'''}}{\longrightarrow} & \widehat{G'''}
\times \widehat{G'''}
\end{array}
\end{equation}
Now $\text{det}\times \text{det}^{\perp}$ above is induced by the degree $2$ isogeny
$$Spec\ \cO[\alpha,\beta,(\alpha^2+\beta^2)^{-1}]
\ra Spec\ \cO[z_1,z_1^{-1},z_2,z_2^{-1}]$$
defined by
$$z_1\mapsto \alpha^2+\beta^2,\ \ \ z_2\mapsto \frac{\alpha+\sqrt{-1}\beta}{\alpha-\sqrt{-1}\beta}.$$
If $\phi_1^{G'}$ and $\phi_2^{G'}$ commute
then clearly $\phi_1^{G'''}$ and $\phi_2^{G'''}$ commute. Conversely, if $\phi_1^{G'''}$ and $\phi_2^{G'''}$ commute then $\phi_1^{G'}\phi_2^{G'}$
and $\phi_2^{G'}\phi_1^{G'}$ coincide on 
$$\alpha^2+\beta^2\ \ \ \text{and}\ \ \ \frac{\alpha+\sqrt{-1}\beta}{\alpha-\sqrt{-1}\beta}.$$
Hence $\phi_1^{G'}\phi_2^{G'}$ and 
$\phi_2^{G'}\phi_1^{G'}$
coincide on   $$(\alpha+\sqrt{-1}\beta)^2.$$ 
Hence
\begin{equation}
\label{openletter}
\phi_1^{G'}\phi_2^{G'}(\alpha+\sqrt{-1}\beta)=\pm\phi_2^{G'}\phi_1^{G'}(\alpha+\sqrt{-1}\beta).\end{equation}
One cannot have the minus sign in \ref{openletter}
as one can see by reducing mod $p$. So, in \ref{openletter},  we have the plus sign.
But then we also have
$$\phi_1^{G'}\phi_2^{G'}(\alpha-\sqrt{-1}\beta)=\phi_2^{G'}\phi_1^{G'}(\alpha-\sqrt{-1}\beta).$$
We get that $\phi_1^{G'}\phi_2^{G'}$ and 
$\phi_2^{G'}\phi_1^{G'}$
coincide on $\alpha$ and $\beta$ and our claim is proved.
\qed

\bigskip

{\it Proof of assertions 2 and 3 in Theorem \ref{mock}}.

We start with assertion 2.
Taking determinants in \ref{fid}, using Condition 2' in the proof of Proposition \ref{fishatnoon}, and finally using \ref{russianrussian},  we get
\medskip
$$
\begin{array}{rcl}
\phi_i^{G'}(\alpha^2+\beta^2) & = & (\alpha^{2p}+\beta^{2p})(u_i^2+v_i^2) \\
\  & \  & \  \\
\  & = & (\alpha^{2p}+\beta^{2p})\cdot \theta_i\\
\  & \  & \  \\
\ & = & \frac{d_i^p}{\phi(d_i)}\cdot (\alpha^2+\beta^2)^p.
\end{array}
$$

\medskip

\noindent It is then trivial to check that 
$$\phi_i^{G'}(\alpha^2+\beta^2-1)$$ is in the ideal generated by $$\alpha^2+\beta^2-1$$ if and only if $$\frac{d_i^p}{\phi(d_i)}=1$$
hence if and only if $$\d d_i=0,$$
 which ends the proof of assertion 2.

We next address assertion 3, so assume 
$$d_1=d_2=:d\in \bZ,\ \ d\neq \pm 1,\ \ d\not\equiv 0\ \ \text{mod}\ \ p.$$ 
Without loss of generality we may assume $\sqrt{-1}\in \cO$.
So $GL_1^c$ is isomorphic to $GL_1\times GL_1$ where  the isomorphism is defined on points by
$$\left(\begin{array}{cc}a & b\\- b & a\end{array}\right)\mapsto (a+\sqrt{-1}b,a-\sqrt{-1}b).$$
Since the closed connected subgroup schemes of $GL_1\times GL_1$ are all kernels of characters it follows that any connected closed subgroup scheme of $GL_1^c$ is of the form $T_{k_1k_2}$ where the latter is given schematically by the equations
$$\left(\begin{array}{cc}\alpha & \beta\\- \beta & \alpha\end{array}\right)^{k_1}\left(\begin{array}{cc}\alpha & - \beta\\  \beta & \alpha\end{array}\right)^{k_2}=1_2,$$
and $k_1,k_2\in \bZ$ are coprime. Equivalently $T_{k_1k_2}$ is given schematically by the equations
\begin{equation}
\label{jazz}
\left(\begin{array}{cc}\alpha & \beta\\- \beta & \alpha\end{array}\right)^{k}=
(\alpha^2+\beta^2)^{l}\cdot 1_2\end{equation}
where $k=k_2-k_1$, $l=k_2$.

Assume now $T_{k_1k_2}$ is $\d^G_i$-horizontal for some $k_1,k_2$ for $i=1,2$; we will derive a contradiction.

 If $k=0$ then $l=\pm 1$ so $T_{k_1k_2}=U^c_1$  and  we are done by assertion 2 of the theorem.

 Assume
now $k\neq 0$. Applying $\phi^{G'}_i$ to \ref{jazz} and using Proposition \ref{fishatnoon} and Remark \ref{fishallday} we get
\begin{equation}
\label{johnny}
\left(\begin{array}{cc}\alpha^p & \beta^p\\- \beta^p & \alpha^p\end{array}\right)^k\cdot
\left(\begin{array}{cc}u & v\\- v & u\end{array}\right)^k= S+M,\end{equation}
where $u=v+1$, $v$ is the  root $\equiv 0$ mod $p$ of the equation
\begin{equation}
\label{shake}
2v^2+2v=\frac{d^p}{\phi(d)}\frac{(\alpha^2+\beta^2)^p}{\alpha^{2p}+\beta^{2p}}-1\end{equation}
in the ring $$\cO(\widehat{G'})=\cO[\alpha,\beta,(\alpha^2+\beta^2)^{-1}]^{\widehat{\ }},$$
the matrix $S$ is scalar,
$$S\in  \cO(\widehat{G'})\cdot 1_2,$$
and the matrix
$$M\in {\mathfrak g}{\mathfrak l}_n(\cO(\widehat{G'}))$$
has entries in the ideal defining $T_{k_1k_2}$.
Let $u_0,v_0\in \cO$ be obtained from $u,v$ by setting $\alpha=1$ and $\beta=0$. Then from \ref{johnny} and \ref{shake} we get
\begin{equation}
\label{rock}
\left(\begin{array}{cc}u_0 & v_0\\- v_0 & u_0\end{array}\right)^k\in \cO\cdot 1_2\end{equation}
and
\begin{equation}
\label{molly}
2v_0^2+2v_0+1-\frac{d^p}{\phi(d)}=0.
\end{equation}
Set $$J=\left(\begin{array}{cc} 0 & 1\\ -1 & 0\end{array}\right)$$ and write
\begin{equation}
\label{studenta}
\left(\begin{array}{cc}u_0 & v_0\\- v_0 & u_0\end{array}\right)^k=
(u_0\cdot 1_2+v_0\cdot J)^k.\end{equation}
Using the binomial formula to expand \ref{studenta} and looking at the upper right corner entry of the matrix in \ref{rock} we get that
\begin{equation}
\label{miss}
\left(\begin{array}{c} k\\ 1\end{array}\right) u_0^{k-1}v_0-\left(\begin{array}{c} k\\ 3\end{array}\right)  u_0^{k-3}v_0^3+\left(\begin{array}{c} k\\ 5\end{array}\right)  u_0^{k-5}v_0^5-...=0.
\end{equation}
Now set $e=|d^{(p-1)/2}|$, where $|\ |$ is the Archimedian absolute value,  and note that the discriminant  $2e^2-1$ of the polynomial in \ref{molly} is a positive rational number. 
So we can and will choose an embedding of ${\mathbb Q}(v_0)$ into ${\mathbb C}$ such that, in this embedding,
$$v_0=\frac{-1+ \sqrt{2e^2-1}}{2},\ \ u_0=\frac{1+ \sqrt{2e^2-1}}{2}$$
where the square root is the real positive one.  In particular $u_0$ and $v_0$ are real in this embedding.
Then by \ref{miss} we have that the complex number
$$(u_0+\sqrt{-1}\cdot v_0)^k$$
is real. So the complex number
\begin{equation}
\label{lucille}
\zeta:=\frac{u_0+\sqrt{-1}\cdot v_0}{|u_0+\sqrt{-1}\cdot v_0|}=\frac{u_0+\sqrt{-1}\cdot v_0}{\sqrt{u_0^2+ v_0^2}}=
\frac{u_0+\sqrt{-1}\cdot v_0}{e}\end{equation}
is a root of unity. Since
\begin{equation}
\label{quar}
\zeta\in {\mathbb Q}(\sqrt{-1},\sqrt{2e^2-1})\end{equation}
and the latter field has degree a divisor of $4$ the order $N$ of $\zeta$ 
must satisfy 
$$\varphi(N)\in \{1,2,4\}$$
 where $\varphi$ is the Euler function.
So the only possibilities for $N$ are:
$$N\in \{1, 2, 3, 4, 5, 6, 8, 10, 12\}.$$
Since the Galois group over ${\mathbb Q}$ of the field in \ref{quar} cannot be cyclic of order $4$ it follows that
$N$ cannot be $5$ or $10$. Now
 equation \ref{lucille} gives 
\begin{equation}
\label{enoughgh}
\text{Re}\ \zeta=\frac{u_0}{e}=\frac{1+\sqrt{2e^2-1}}{2e},
\end{equation}
so in particular $\text{Re}\ \zeta>0$. Hence the only possibilities for $\text{Re}\ \zeta$ are
\begin{equation}
\label{valuesofzeta}
\text{Re}\ \zeta\in \{1, \ \frac{1}{2},\ \frac{\sqrt{2}}{2},\ \frac{\sqrt{3}}{2}\}.\end{equation}
The case $\text{Re}\ \zeta=1$ of equation \ref{enoughgh} yields $e=1$, hence $d=\pm 1$, a contradiction. For the other $3$ values of $\text{Re}\ \zeta$ in \ref{valuesofzeta} equation \ref{enoughgh} gives values of $e$ that are not in ${\mathbb Q}$, which is again a contradiction.
This ends the proof.
\qed

\begin{remark}
Note that the above proof works if one replaces the hypotheses
$$d\in \bZ,\ \ \ d\neq \pm 1,\ \ \ d\not\equiv 0\ \ \ \text{mod}\ \ \ p$$
by the hypotheses
$$d\in \bZ_{(p)}^{\times},\ \ \ |d|>1.$$ 
\end{remark}

\bigskip

{\it Proof of Proposition  \ref{elisabeth}}.
Let us place ourselves,  in what follows, in the global situation and consider a symmetric matrix $q\in GL_n(\cO_{F,M})$,  the vertical Levi-Civita connection $$(\d^1,...,\d^n)$$ attached to $q$ at ${\mathfrak P}$, the attached Frobenius lifts $$(\phi^1,...,\phi^n),$$
  and the curvature $(\varphi_{ij})$. Let, as before, 
$$\phi^1_i=\sigma_i^{-1}\phi^i\sigma_i.$$
Then we have
\begin{equation}
\label{coffee soon}
\varphi_{ij}=\frac{1}{p}(\phi_i^1\phi_j^1-\phi_j^1\phi_i^1).\end{equation}
Set $$\Phi_{ij}:=\varphi_{ij}(x),\ \ \ \Delta_i:=\d_i^1 x,$$
 and recall that the Christoffel symbols $\gamma_i$ of the first kind are given by
the equalities
$$\sigma_i^{-1}\gamma_{i}=\Delta^t_i\cdot \phi^1(q_i)\cdot x^{(p)}.$$
 We have
$$
\begin{array}{rcl}
\phi_i^1\phi_j^1(x) & = & \phi_i^1(x^{(p)}+p\Delta_j)\\
\ & \ & \ \\
\ & = & (x^{(p)}+p\Delta_i)^{(p)}+ p\phi^1_i(\Delta_j)\\
\  & \  & \ \\
\ & \equiv & x^{(p^2)}+p\Delta_j^{(p)}\ \ \ \text{mod}\ \ p,
\end{array}
$$
hence
\begin{equation}
\label{richard}
\begin{array}{rcll}
\Phi_{ij} & \equiv & \Delta_j^{(p)}-\Delta_i^{(p)} &\ \\
\ &\ & \ & \ \\
\ & \equiv & ((\sigma_j^{-1}(q))^{(p^2)})^{-1} (x^{(p^2)t})^{-1} 
(\sigma_j^{-1}\gamma_j)^{(p)t}\\
\ &\ & \ & \ \\
\ & \  & 
- ((\sigma_i^{-1}(q))^{(p^2)})^{-1} (x^{(p^2)t})^{-1} (\sigma_i^{-1}\gamma_i)^{(p)t} &  \text{mod}\ \ p\\
\ & \ & \ & \ \\
\ & \equiv & ((\sigma_j^{-1}(q))^{(p^2)})^{-1}  (\sigma_j^{-1}\gamma_j)^{(p)t}\\
\ & \ & \ & \ \\
\ & \ & - ((\sigma_i^{-1}(q))^{(p^2)})^{-1}  (\sigma_i^{-1}\gamma_i)^{(p)t} &  \text{mod}\ \ (p,x-1).
\end{array}\end{equation}
Combining the congruences \ref{richard}
with the congruences \ref{greene} one immediately gets the congruences \ref{burton}.\qed

\bigskip

\subsection{Construction of an \'{e}tale cover}\ 

We discuss, in what follows, a construction that will be later used to prove the  existence of our correspondences.

Let 
$$y=(y_1,...,y_n),\ \ \ z=(z_1,...,z_n)$$ be two $n$-tuples of matrices of size $n\times n$ with intedeterminates as entries, 
$$y_i=(y_{ijk}),\ \ \ z_i=(z_{ijk}).$$ 
Consider the system of linear equations
in $n^3$ unknowns $z_{ijk}$, with coefficients in the ring $\bZ[y]$,
\begin{equation}
\label{systemncube}
\begin{array}{rcll}
(y_i^tz_i)_{jk}+(z_i^ty_i)_{jk} & = & 0, & i,j,k=1,...,n,\ \ \  j\leq k,\\
z_{ikj}-z_{jki} & = & 0, & i,j,k=1,...,n,\ \ \  i<j.\end{array}\end{equation}
There are $$\frac{n^2(n+1)}{2}$$ equations in the first row of \ref{systemncube}  and $$\frac{n^2(n-1)}{2}$$
 equations in the second row of \ref{systemncube}
so there are are $n^3$ equations in all. So 
the matrix of the system \ref{systemncube} is square and
one can consider the determinant of this matrix which we denote by
\begin{equation}
\label{detdetdet}
D(y):=D(y_1,...,y_n)\in \bZ[y].\end{equation}
Of course $D(y)$ is well defined only up to sign because the order of the variables and the order of the equations has not been specified.

 \begin{example}
 For $n=2$, $y=(y_1,y_2)$, 
$$y_1=
 \left(\begin{array}{cc}
 y_{111} & y_{112}\\ y_{121} & y_{122}\end{array}\right),\ \ \ 
y_2= \left(\begin{array}{cc}
 y_{211} & y_{212}\\ y_{221} & y_{222}\end{array}\right),\ \ \ 
y_{1|2}:= \left(\begin{array}{cc}
 y_{112} & y_{211}\\ y_{122} & y_{221}\end{array}\right),
$$
one gets
 \begin{equation}
\label{n2}
 D(y_1,y_2)=\pm \ \det(y_1)\cdot \det(y_2)\cdot \det (y_{1|2}).
 \end{equation}
 \end{example}
 
 Going back to an arbitrary $n$ and writing, as usual,  $1=1_n$ we may consider the integer
$D(1,...,1)\in \bZ$.

\begin{lemma}
\label{snoringg}
For any odd prime $p$ one has:
$$D(1,...,1)\not\equiv 0\ \ \ \text{mod}\ \ \ p.$$
\end{lemma}

In other words $D(1,...,1)$ is $\pm 1$ times (possibly) a power of $2$.

\medskip

{\it Proof}.
Assume an odd prime $p$ divides $D(1,...,1)$.
Then the system\begin{equation}
\label{system1}
\begin{array}{rcll}
z_{ijk}+z_{ikj} & = & 0, & i,j,k=1,...,n,\\
z_{ikj}-z_{jki} & = & 0, & i,j,k=1,...,n,\end{array}\end{equation}
 has a zero determinant in ${\mathbb F}_p$ so it has a non-trivial solution
 $(\zeta_{ijk})$ in ${\mathbb F}_p$.
So 
$$
\zeta_{ijk}=- \zeta_{ikj}=-\zeta_{jki}.$$
Using the latter $3$ times one gets $$2\zeta_{ijk}=0$$ hence $$\zeta_{ijk}=0,$$ a contradiction.
\qed

\bigskip

Assume now ${\mathcal B}$ is a Noetherian ring, fix an integer $n\geq 2$, and consider the polynomial $D(y)\in \bZ[y]$ in \ref{detdetdet}.
Also let  $$A_1,...,A_n,B_1,...,B_n$$ be $n\times n$ symmetric matrices with entries in ${\mathcal B}$, let 
$$b:=\det(B_1)\cdot ... \cdot \det(B_n)\in {\mathcal B},\ \ \ B_b=B[1/b],$$
 and define the ring ${\mathcal C}$
associated to the data $({\mathcal B},A,B)$
 by the formula
\begin{equation}
\label{theringC}
{\mathcal C} :={\mathcal C}({\mathcal B},A,B):= \frac{{\mathcal B}_b[y, D(y)^{-1}]}{((y_i^tA_iy_i-B_i)_{jk},(A_i(y_i-1))_{kj}-(A_j(y_j-1))_{ki})},
\end{equation}
where $A$ is the $n$-tuple $(A_i)$ and $B$ is the $n$-tuple $(B_i)$.

Note that the triples $({\mathcal B},A,B)$ are the objects of an obvious category: a morphism
$$({\mathcal B},A,B)\ra ({\mathcal B}',A',B')$$
 is a morphism of rings ${\mathcal B}\ra {\mathcal B}'$ which sends the matrices $A,B$ into $A',B'$ respectively. 
Then we clearly obtain a functor
$$\{({\mathcal B},A,B)\}\ra \{\text{rings}\},\ \ \ \ ({\mathcal B},A,B)\mapsto {\mathcal C}({\mathcal B},A,B).$$
For a morphism as above we have
$${\mathcal C}({\mathcal B}',A',B')\simeq {\mathcal C}({\mathcal B},A,B)\otimes_{\mathcal B}{\mathcal B}'.$$
So for any triple $({\mathcal B},A,B)$ we have
$${\mathcal C}({\mathcal B},A,B)\simeq {\mathcal C}({\mathcal B}^{\text{univ}},A^{\text{univ}},B^{\text{univ}})\otimes_{{\mathcal B}^{\text{univ}}}{\mathcal B},$$
where $A^{\text{univ}}$ and $B^{\text{univ}}$ are two $n$-tuples of  symmetric matrices with indeterminate coefficients on and above  the diagonal,
$${\mathcal B}^{\text{univ}}:=\bZ[A^{\text{univ}},B^{\text{univ}}],$$
is the polynomial ring in these variables, and ${\mathcal B}^{\text{univ}}\ra {\mathcal B}$ is given by $A^{\text{univ}}\mapsto A$, $B^{\text{univ}}\mapsto B$.
  
For ${\mathcal C}={\mathcal C}({\mathcal B},A,B)$ we have a natural map of  schemes
\begin{equation}
\label{themappi}
\pi:Y:=Spec\ {\mathcal C}\ra X:=Spec\ {\mathcal B}.\end{equation}

\begin{lemma} \label{gore}
The map $\pi:Y\ra X$ is \'{e}tale.
\end{lemma}

{\it Proof}.
Consider a diagram of rings
$$
\begin{array}{rcl}
{\mathcal B} & \stackrel{\pi}{\longrightarrow} & {\mathcal C}\\
v \downarrow &\ & \downarrow u\\
{\mathcal D} & \stackrel{\rho}{\longrightarrow} & {\mathcal D}/I
\end{array}
$$
where $I\subset {\mathcal D}$ is an ideal with $I^2=0$. We need to show that there is a unique map 
$$w:{\mathcal C}\ra {\mathcal D}$$ such that 
$$\rho\circ w=u,\ \ \ w\circ \pi=v.$$
 Set 
$$v(A_i)=a_i,\ \ \ v(B_i)=b_i,\ \ \ u(y_i)=\rho(\lambda_i),$$
 with $a_i,b_i,\lambda_i$ matrices with entries in ${\mathcal D}$. So we have that 
\begin{equation}\label{australia}
\begin{array}{rcl}
\lambda_i^ta_i\lambda_i-b_i & = & c_i\\
(a_i(\lambda_i-1))_{kj}-(a_j(\lambda_j-1))_{ki}
 & = & f_{ijk}
\end{array}
\end{equation}
for some symmetric matrices $c_i$ with coefficients in $I$ and some elements $f_{ijk}\in I$ with $$f_{ijk}=-f_{jik}.$$
 To find $w$ is the same as to find elements $$w(y_{ijk})=\lambda_{ijk}+\epsilon_{ijk},$$
 with $\epsilon_{ijk}\in I$, such that if $\epsilon_i=(\epsilon_{ijk})$ then
\begin{equation}
\label{atlantis}
\begin{array}{rcl}
(\lambda_i^t+\epsilon^t_i)a_i(\lambda_i+\epsilon_i)-b_i & = & 0,\\
(a_i(\lambda_i+\epsilon_i-1))_{kj}-(a_j(\lambda_j+\epsilon_j-1))_{ki} & =  & 0.
\end{array}\end{equation}
In view of \ref{australia}, if we set
\begin{equation}
\label{she}
\zeta_i=a_i\epsilon_i,
\end{equation}
with $\zeta_i=(\zeta_{ijk})$,
then
the equations \ref{atlantis} can be rewritten as
\begin{equation}
\label{chirpchirp}
\begin{array}{rcl}
(\lambda_i^t\zeta_i)_{jk}+(\zeta_i^t \lambda_i)_{jk} & = & -c_{ijk},\\
\zeta_{ikj}-\zeta_{jki} & = & -f_{ijk},
\end{array}\end{equation}
where $i,j,k=1,...,n$. Now the system \ref{chirpchirp} is, of course, equivalent to the system consisting of the same equations but where the indices satisfy, in addition,  $j\leq k$ for the equations in the first line and $i<j$ for the equations in the second line of \ref{chirpchirp}.
Since 
$$D(\lambda_1,...,\lambda_n)$$
 is invertible in ${\mathcal D}$ (because it is invertible mod $I$) the system 
\ref{chirpchirp} has a unique solution $(\zeta_{ijk})$ with entries in $I$.
Since $v(b)$ is invertible in ${\mathcal D}$ (because it is invertible mod $I$) it 
follows that $\det(b_i)$ and hence $\det(a_i)$ are invertible in ${\mathcal D}$ so  the system \ref{she} has a unique solution $(\epsilon_{ijk})$ with entries in $I$  
 and the Lemma is proved.
\qed

\bigskip

Assume, in what follows that $F$ is a number field, $0\neq M\in \bZ$ is an even integer, $\cO=\cO_{F,M}=\cO[1/M]$, let $p$ be an odd prime in $\bZ$ unramified in $F$ and not dividing $M$,
let ${\mathfrak P}$ be a prime in $\cO$ above $p$ and let $\overline{\cO}:=\cO/{\mathfrak P}$. Assume furthermore that $X=Spec\ {\mathcal B}$ is a smooth connected scheme over $\cO$ with geometrically irreducible fibers.  Also we denote by an upper bar tensorization over $\cO$ with $\overline{\cO}$.

\begin{lemma}
\label{teatime} Assume $\overline{b}\in \overline{\mathcal B}^{\times}$ and 
$$A_i\equiv B_i\ \ \ \text{mod}\ \ \ {\mathfrak P}{\mathcal B}$$
 for all $i$.
Then the map $\overline{\pi}:\overline{Y}\ra \overline{X}$ has a section for which the induced map between the corresponding rings pulls back  $y_i$ into $1$.
\end{lemma}

\medskip

In particular the map $\overline{\pi}:\overline{Y}\ra \overline{X}$ is surjective, hence an \'{e}tale cover.

\medskip

{\it Proof}.
 The map 
 \begin{equation}
 \label{themap}
 {\mathcal B}[y]\ra \overline{\mathcal B},\ \ \  y_i\mapsto 1\end{equation}
 sends 
$$\begin{array}{rcl}
y_i^tA_iy_i-B_i & \mapsto & 0\\
\  & \  & \  \\
(A_i(y_i-1))_{kj}-(A_j(y_j-1))_{ki} & \mapsto & 0\\
\  & \  & \  \\
b& \mapsto & \overline{b} \in \overline{\mathcal B}^{\times}\\
\  & \  & \  \\
D(y)& \mapsto & (D(1,...,1)\ \ \text{mod}\ \ {\mathfrak P})\in {\mathbb F}_p^{\times}\ \ \text{by Lemma \ref{snoringg}}.\end{array}
$$
So we have an induced map 
$\overline{\mathcal C}\ra \overline{\mathcal B}$; the latter induces a section
of the projection
$\overline{\pi}:\overline{Y}\ra \overline{X}$
and we are done.\qed

\bigskip

Denote now by
$$\overline{\sigma}:\overline{X}\ra \overline{Y}$$
the section of $\overline{\pi}:\overline{Y}\ra \overline{X}$ constructed in Lemma \ref{teatime}.
The image of the section $\overline{\sigma}$, 
$$\overline{Y}^1:=\overline{\sigma}(\overline{X}),$$
 is a closed subscheme of $\overline{Y}$, so by Lemma \ref{gore} and dimension considerations it is an irreducible component of $\overline{Y}$; since, again by Lemma \ref{gore},
$\overline{Y}$ is smooth, it follows that $\overline{Y}^1$ is a connected component of $\overline{Y}$. Let
$$\overline{Y}^2:=\overline{Y}\backslash \overline{Y}^1$$
and let $\overline{e}\in \cO(\overline{Y})$ be the idempotent which is $1$ on $\overline{Y}^1$ and $0$ on $\overline{Y}^2$.
 Finally let $e\in \cO(Y)={\mathcal C}$ be any lift of $\overline{e}$,  consider the 
scheme
$$Y^0:=Spec\ {\mathcal C}[1/e],$$
and the open immersion
$Y^0\subset Y$.
Also define the formal scheme
\begin{equation}
\label{close to the edge}
{\mathcal Y}:=(Y^0)^{\widehat{\mathfrak P}}.\end{equation}
Since 
$$\overline{{\mathcal Y}}=\overline{Y}^1$$
it follows that ${\mathcal Y}$ is a connected component of $Y^{\widehat{\mathfrak P}}$
where 
 $\widehat{\mathfrak P}$ means, as usual,  ${\mathfrak P}$-adic completion.
 Clearly

\begin{corollary}
\label{woodpecker}
 The induced map
\begin{equation}
\label{sculpt}
\overline{\pi}:\overline{Y}^1\ra \overline{X}\end{equation}
is an isomorphism and its inverse pulls back $y_i$ into $1$. 
\end{corollary}

\begin{corollary}
\label{lully}
The map of formal schemes 
$$\pi^{\widehat{\mathfrak P}}:{\mathcal Y}=(Y^0)^{\widehat{\mathfrak P}}\ra X^{\widehat{\mathfrak P}}$$
is an isomorphism.\end{corollary}

{\it Proof}.
By Corollary \ref{woodpecker} 
the map ${\mathcal B}\ra {\mathcal C}[1/e]$ induces an isomorphism after tensorization with $\overline{\cO}$. Hence the map ${\mathcal B}^{\widehat{\mathfrak P}}\ra {\mathcal C}[1/e]^{\widehat{\mathfrak P}}$  is an isomorphism because
 $p$ is a non-zero divisor in ${\mathcal C}[1/e]$; the latter fact follows from the fact that ${\mathcal C}$ is \'{e}tale, hence flat, over ${\mathcal B}$; cf. Lemma \ref{gore}.
\qed

\bigskip

\subsection{The case $X=GL_n$}
\label{CaseGLn}

\ 

The aim of this subsection is to prove Theorem \ref{algebraization}.

We continue to assume
$F$ is a number field, $0\neq M\in \bZ$ is an even integer, and $\cO=\cO_{F,M}=\cO[1/M]$. Let $p$ be an odd prime in $\bZ$ unramified in $F$ and not dividing $M$,
let ${\mathfrak P}$ be a prime in $\cO$ above $p$ and let $\overline{\cO}:=\cO/{\mathfrak P}$.
Furthermore let $\phi=\phi_{\mathfrak P}:\cO\ra \cO$ be the Frobenius element attached to  ${\mathfrak P}$,
let $q_1,...,q_n\in GL_n(\cO)$ be symmetric matrices,  and set
\begin{equation}
\label{hypo}
 {\mathcal B}=\cO[x,\det(x)^{-1}],\ \ \ A_i=x^{(p)t}\phi(q_i)x^{(p)},\ \ \ B_i=(x^tq_ix)^{(p)}.\end{equation}
With these data the ring ${\mathcal C}={\mathcal C}({\mathcal B},A,B)$ in \ref{theringC} becomes
\begin{equation}
\label{theringCp}
{\mathcal C}=\frac{\cO[x,\det(x)^{-1},b^{-1}, y, D(y)^{-1}]}{((y_i^tA_iy_i-B_i)_{jk},(A_i(y_i-1))_{kj}-(A_j(y_j-1))_{ki})}
\end{equation}
and the map of  schemes $\pi$ in \ref{themappi} becomes the map
\begin{equation}
\label{themappip}
\pi:Y:=Spec\ {\mathcal C}\ra X:=G:=GL_n=Spec\ {\mathcal B}\end{equation}
induced by ${\mathcal B}\ra {\mathcal C}$, $x\mapsto x$. By Lemma
\ref{gore} the map $\pi$ is \'{e}tale.

Consider now the maps
\begin{equation}
\label{sarat}
\varphi_{i}:Y\ra G\end{equation}
induced by the ring homomorphisms $\varphi_i:\cO(G)\ra \cO(Y)$
satisfying
\begin{equation}
\label{parbuclat}
\varphi_i(a)=\phi(a),\ \ \ a\in \cO,
\end{equation}
 and sending
\begin{equation}
\label{iron}
x\mapsto \varphi_{pi}(x):=\text{class}(x^{(p)}y_i)\in \cO(Y).\end{equation}
One can then consider the induced map between ${\mathfrak P}$-adic completions
$$\varphi_{i}^{\widehat{\mathfrak P}}:Y^{\widehat{\mathfrak P}}\ra G^{\widehat{\mathfrak P}}$$
and the restriction of the latter,
$$\varphi_{i}^{\widehat{\mathfrak P}}:{\mathcal Y}\ra G^{\widehat{\mathfrak P}}$$
where  ${\mathcal Y}$ is as in \ref{close to the edge}.
Then  we have:

\begin{lemma}
\label{ultima}
For each $i=1,...,n$ the map
$\varphi_{i}^{\widehat{\mathfrak P}}:{\mathcal Y}\ra G^{\widehat{\mathfrak P}}$
equals the composition 
$$\phi_i^{G_{\mathfrak P}} \circ \pi^{\widehat{\mathfrak P}}: {\mathcal Y}\ra G^{\widehat{\mathfrak P}}
\ra G^{\widehat{\mathfrak P}}.$$
\end{lemma}

{\it Proof}.
Let $\Lambda_i\in \cO(G^{\widehat{\mathfrak P}})$ be the pull back of  $\text{class}(y_i)\in \cO({\mathcal Y})$ via $(\pi^{\widehat{\mathfrak P}})^{-1}$.
Then clearly we have

\medskip

1) $\Lambda_i\equiv 1$ mod $p$ (by Corollary \ref{woodpecker})

2) $\Lambda_i^tA_i\Lambda_i=B_i$;

3) $(A_i(\Lambda_i-1))_{kj}=(A_j(\Lambda_j-1))_{ki}$. 

\medskip

\noindent But by the proof of Theorem \ref{LC} there is a unique tuple $\Lambda_i$ with properties 1, 2, 3 and the Frobenius lifts $\phi_i^{G_{\mathfrak P}}$ corresponding to the  Levi-Civita connection attached to $(q_1,,...,q_n)$ send $x$ into $x^{(p)}\Lambda_i$. This ends our proof.
\qed

\bigskip

{\it Proof of Theorem \ref{algebraization}}.
With notation as above we let $Y_{p/G}$ be the connected component of $Y$ containing $\overline{\mathcal Y}$ (which is a regular scheme hence irreducible).
Also we let $\pi_{p/G}:Y_{p/G}\ra G$ and $\varphi_{pi}:Y_{p/G}\ra G$ be the restrictions of $\pi:Y\ra G$ and $\varphi_{i}:Y\ra G$ respectively and we let ${\mathcal Y}_{p/G}={\mathcal Y}$. Then all assertions of Theorem \ref{algebraization} are satisfied.
\qed

\bigskip

\begin{remark}\label{nuuita}
\ 

1) It is clear that the conjunction of Lemma \ref{gore} and Corollaries  \ref{lully} and \ref{ultima} implies our Theorem \ref{algebraization}. It is also trivial to see that our arguments in the proof of Theorem \ref{algebraization} can also be used to prove the existence part of  our Theorem 
\ref{levi Civita}; however the proof that we already gave for 
the existence part of
Theorem \ref{levi Civita} has the advantage of also immediately yielding our proof of Proposition \ref{congruences}. 

2) The construction of  $Y_{p/G}$  
and of the maps $\pi_p,\varphi_{pi}$ in the proof of  Theorem \ref{algebraization} was entirely canonical/functorial. So our construction of mixed curvature in Definition \ref{woken}
is canonical.  \end{remark}

\subsection{The case $X=GL_1^c$}\ 

The aim of this subsection is to prove Theorem \ref{algebraization1} and Propositions \ref{maro}, \ref{galben}.

We  consider the situation in the previous section with $n=2$. In addition, we let $d_1,d_2\in \cO^{\times}$, we let $\alpha, \beta$ be $2$ indeterminates,  and set
$$a_i=(\alpha^{2p}+\beta^{2p})\cdot \phi(d_i),\ \ \ b_i:=(\alpha^2+\beta^2)^p\cdot d_i^p,$$
$${\mathcal B}':=\cO[\alpha,\beta,(\alpha^2+\beta^2)^{-1}],\ \ \ A'_i=a_i\cdot 1_2,\ \ \ B'_i:=b_i\cdot 1_2,\ \ \ b'=b_1^2b_2^2.$$
As usual, we set $$G'=GL_1^c=Spec\ {\mathcal B}'$$ viewed as embedded into $$G=GL_2=Spec\ {\mathcal B}.$$
With 
$$A'=(A'_1,A'_2),\ \ \ B'=(B'_1,B'_2)$$
 we consider the ring 
$${\mathcal C}':={\mathcal C}({\mathcal B}',A',B')$$
associated to the data $({\mathcal B}',A',B')$ as in \ref{theringC}. Note that $b_i$ are invertible in ${\mathcal C}'$ hence so are 
$b'$ and $a_i$ hence, setting 
$$\theta_i=\frac{b_i}{a_i}\in {\mathcal B}'_{a_1a_2},$$ we have
$$
{\mathcal C}' = \frac{{\mathcal B}'_{a_1a_2}[y,D(y)^{-1}]}{((y_i^ty_i-\theta_i)_{jk}, (a_i(y_i-1_2))_{kj}-a_j(y_j-1_2)_{ki})}.
$$
By Lemmas \ref{gore} and \ref{teatime} we have that the map   
$$Z':=Spec\ {\mathcal C}'\ra G'=Spec\ {\mathcal B}'$$
 is \'{e}tale and its reduction mod $p$ has a section defined by the map 
\begin{equation}
\label{aspir}
\overline{{\mathcal C}'} \ra \overline{{\mathcal B}'},\ \ \ y_i\mapsto 1_2.\end{equation}
Now  let 
$$t_i:=\text{class}(\text{tr}(y_i)):=\text{class}(y_{i11}+y_{i22})\in {\mathcal C}'$$
and
$$\tau_i:=\text{class}(\det(y_i)+\theta_i) \in {\mathcal C}'$$
Then \ref{aspir} sends $t_i$ and $\tau_i$ into
$$\overline{2} \in (\overline{{\mathcal B}'})^{\times}.$$
Setting $s=t_1t_2\tau_1\tau_2$
we get an induced map
$$\overline{{\mathcal C}'_{s}}\ra \overline{{\mathcal B}'}.$$
Let
$$Y':=Spec\ {\mathcal C}'_{s}$$
and denote by 
$$\pi':Y'\ra G'$$
 the induced morphism which is, of course, still \'{e}tale.
We get a section 
$$\overline{\sigma}:\overline{G'}\ra \overline{Y'}$$
of the projection
$$\overline{\pi'}:\overline{Y'}\ra \overline{G'}.$$
Exactly as in the case of $GL_n$, denoting by ${\mathcal Y}'$ the connected component of $(Y')^{\widehat{\mathfrak P}}$ containing
$\overline{\sigma}(\overline{G'})$ we get an isomorphism
$${\mathcal Y}'\ra (G')^{\widehat{\mathfrak P}}.$$
We will next construct for $i=1,2$ morphisms
$$\varphi'_i=\varphi'_{pi}:Y'\ra G'$$
as follows.
 We already have at our disposal the morphisms
$\varphi_i:Y\ra G$
in \ref{sarat}.
We want to construct $\varphi'_i$ so as to be induced by $\varphi_i$.
Note 
that the canonical map ${\mathcal B}\ra {\mathcal B}'$ sends 
$$A_i\mapsto A'_i,\ \ \ B_i\mapsto B'_i$$
so it induces a canonical map 
$$\text{can}:{\mathcal C}\ra {\mathcal C}'\ra {\mathcal C}'_s.$$
On the other hand we have  the following:

\begin{lemma}
\label{o3}
Let $J$ be the ideal in $\cO(G)$ defining $G'$; so $J$ is  generated by
$$x_{11}-x_{22},\ \ \ x_{12}+x_{21}.$$
Then $J$ is sent into $0$ by the map
$${\mathcal B}\stackrel{\varphi_i}{\longrightarrow} {\mathcal C}\stackrel{\text{can}}{\longrightarrow} {\mathcal C}'_s.$$
In particular the maps $\varphi_i:Y\ra G$ induce maps
$\varphi'_i:Y'\ra G'$.
\end{lemma}

{\it Proof}.
Recall from the proof of Theorem \ref{algebraization} that $\varphi_i(x)$ was defined as the class of $x^{(p)}\cdot y_i$
in ${\mathcal C}$. So in order to conclude we need to show that the image  $y_i'$ of $y_i$ in $GL_2({\mathcal C}'_s)$
belongs to $GL_1^c({\mathcal C}'_s)$. Pick an $i=1,2$ and write 
$$y_i'=\left(\begin{array}{cc} u & v\\ w & z\end{array}\right),\ \ \ u,v,w,z\in {\mathcal C}'_s.$$
The equality 
\begin{equation}
\label{fal}
(y_i')^t y_i'=\theta_i\end{equation}
gives
\begin{equation}
\label{tings}
\begin{array}{rcl}
u^2+ w^2 & = & \theta_i,\\
\  & \  & \  \\
uv+wz & = & 0,\\
\  & \  & \  \\
v^2+z^2 & =&  \theta_i.\end{array}
\end{equation}
A formal manipulation of the first $2$ equations in \ref{tings}
gives
\begin{equation}
\label{zag}
w(wv-uz)=\theta_i v.\end{equation}
On the other hand taking the determinant in \ref{fal} we get
$$(uz-wv)^2=\theta_i^2$$
so
$$(\det(y_i')+\theta_i)(\det(y_i')-\theta_i)=0.$$
Since $\det(y_i')+\theta_i$ is invertible in ${\mathcal C}'_s$ we get
$$uz-wv=\det(y'_i)=\theta_i.$$
Combining with \ref{zag} we get 
$$v=- w.$$
Subtracting the first and third equations in \ref{tings} we get
$$(u+z)(u-z)=0.$$
But now $u+z=\text{tr}(y_i')$ is invertible in ${\mathcal C}'_s$. So we get
$$u=z$$
which ends the proof of the fact that $y_i'$ belongs to $GL_1^c({\mathcal C}'_s)$.
\qed

\bigskip

{\it Proof of Theorem \ref{algebraization1}}.
 With the notation above we let
 ${\mathcal Y}'_{p/G'}={\mathcal Y}'$ and we let 
 $Y'_{p/G'}$ be the irreducible component of $Y'$ that contains $\overline{{\mathcal Y}'}$.
Furthermore we let $\pi'_p, \varphi'_{p/G'}:Y'_{p/G'}\ra G'$ be the restrictions of $\pi',\varphi_i:Y'\ra G'$. Then, clearly, all the assertions of Theorem \ref{algebraization1} follow.
\qed

\bigskip

{\it Proof of   Proposition \ref{maro}}. 
By our construction and the formula \ref{n2}, the tensor product ${\mathcal C}'_s \otimes_{{\mathcal B}'} E'$ (with ${\mathcal C}'$ over ${\mathcal B}'$ viewed via $\pi'$)
is isomorphic to
\begin{equation}\label{Gprime}
M':=\frac{E'[y,g(y)^{-1}]}{((y_i^ty_i-\theta_p)_{jk},\ y_{112}-y_{211}+1,\ y_{122}-y_{221}-1)}
\end{equation}
where $y=(y_1,y_2)$ and
$$g(y):=\det(y_1)\cdot \det(y_2)\cdot \text{tr}(y_1)\cdot \text{tr}(y_2)\cdot (\det(y_1)+\theta_p)\cdot (\det(y_2)+\theta_p).$$
Set
$$L':=L'_p:=\frac{E'[z]}{(2z^2+2z+1-\theta_p)}$$
where $z$ is a variable. The discriminant of $2z^2+2z+1-\theta_p$  is $2\theta_p-1$ which is not a square in $E'$ because $\alpha^{2p}+\beta^{2p}$ is a product of distinct linear factors. So $L'$ is a quadratic field extension of $E'$.

We will construct in what follows a natural isomorphism $L'\simeq M'$.

 Let $v=v_p\in L'$ be the class of $z$ and let $u=1+v$.
Then the homomorphism
\begin{equation}
\label{o1}
E'[y]\ra L',\ \ y_1\mapsto \left(\begin{array}{cc} u & v\\ - v& u\end{array}\right),\ \ y_2\mapsto \left(\begin{array}{cc} u & - v\\  v& u\end{array}\right)\end{equation}
is trivially seen to factor through a homomorphism $M'\ra L'$. We also claim that the homomorphism
\begin{equation}
\label{o2}
E'[z]\ra M',\ \ z\mapsto y'_{112}:=\text{class}(y_{112})\end{equation}
factors through a homomorphism $L'\ra M'$. This can be seen as follows. By an argument similar to the one in the proof of Lemma
\ref{o3} the classes $y'_i$ of $y_i$ in $M'$ have the form 
$$y'_i=\left(\begin{array}{cc} u_i & v_i\\ - v_i& u_i\end{array}\right)$$
with
$$u_2=1+v_1,\ \ \ u_1=1-v_2.$$
From the equations 
$$u_1^2+v_1^2=\theta_p,\ \ \ u_2^2+v_2^2=\theta_p$$
we get
$$1-2v_2+v_2^2+v_1^2=\theta_p,\ \ \ 1+2v_1+v_1^2+v_2^2=\theta_p.$$
Subtracting the last 2 equations we get $v_1=- v_2$ hence $u_1=u_2$. So $v_1=y'_{112}$ is a root of $2z^2+2z+1-\theta_p$ and our claim is proved. Finally, using the above considerations it is trivial to check that the two morphisms \ref{o1} and \ref{o2} are inverse to each other. This ends the construction of the isomorphism $L'\simeq M'$. Since $M'$ is a field we get, in particular, that $Y':=Spec\ {\mathcal C}'_s$ itself is irreducible, so $Y'_{p/G'}=Y'$.
The Proposition now follows easily by using  formula \ref{iron}.
\qed

\bigskip

{\it Proof of Proposition \ref{galben}}.  
To construct the correspondences $\Gamma'''_{pi}$
in Proposition \ref{galben} note that one  has 
 $$\theta_p=\frac{d^p(t^2+1)^p}{\phi_p(d)(t^{2p}+1)}\in E''',$$
 so the element $v_p\in L'_p$ is quadratic over $E'''$. Then  
  one
 can take
 $$Y'''_{p/E'''}=Spec\ L'''_p, \ \ L'''_p=E'''(v_p),$$
 one can take $\pi'''_p$ to be the inclusion
 $E'''\subset L'''_p$, and one can take
 $$\varphi'''_{p1},\varphi'''_{p2}:E'''\ra L'''_p$$
 to act on $F$ via $\phi_p$  and act on $t$ 
via the formulae
\medskip
$$
\varphi'''_{p1} (t)=\frac{u_pt^p-v_p}{v_pt^p+u_p},
\ \ \ 
\varphi'''_{p2} (t)=\frac{u_pt^p+v_p}{-v_pt^p+u_p}.
$$
In order to conclude the proof of Proposition \ref{galben} we need to check that:

\medskip

{\it Claim}. 
The correspondences $\Gamma'''_{pi}$ are categorically reduced.

\medskip

 Indeed if this is checked then the uniqueness of $\Gamma'''_{pi}$ is clear.

We check the Claim  for $i=1$; the case $i=2$ is similar.
Denote by $\overline{t}_1,\overline{t}_2\in L_p'''$ the images of $t_1,t_2$
hence 
$$\overline{t}_1=t,\ \ \ \overline{t}_2= \frac{u_p t^p-v_p}{v_pt^p+u_p}=\frac{v_p(\overline{t}_1^p-1)+\overline{t}_1^p}{v_p(\overline{t}^p_1+1)+1}.$$ 
One gets
$$(\overline{t}_2\overline{t}_1^p+\overline{t}_2-\overline{t}_1^p+1)v_p=\overline{t}_1^p-\overline{t}_2.$$
We claim that 
$\overline{t}_2\overline{t}_1^p+\overline{t}_2-\overline{t}_1^p+1\neq 0$. Indeed if
$\overline{t}_2\overline{t}_1^p+\overline{t}_2-\overline{t}_1^p+1= 0$ we get 
$\overline{t}_1^p=\overline{t}_2$
hence
$$t^p=\frac{u_p t^p-v_p}{v_pt^p+u_p},$$
which implies $v_p t^{2p}=-v_p$, a contradiction. So we can express
\begin{equation}
\label{buccilatto}
v_p=\frac{\overline{t}_1^p-\overline{t}_2}{\overline{t}_2\overline{t}_1^p+\overline{t}_2-\overline{t}_1^p+1}.
\end{equation}
 In particular 
 \begin{equation}
 \label{madel}
 L'''_p=E'''(v_p)=F(\overline{t}_1,\overline{t}_2),
 \end{equation} which ends the proof of the Claim, and hence of Proposition \ref{galben}.
 \qed
 
  \section{Appendix: Classical Levi-Civita connection revisited}
  
 The aim of this Appendix is to quickly revisit the classical theory of the Levi-Civita connection \cite{KN} with an emphasis on the analogy with the arithmetic case. This analogy  is rather ``indirect" in that it requires, as a preliminary,
 a re-thinking of 
the classical paradigm; cf., especially, our concepts of vertical and mixed Levi-Civita connection below. 
 
 \subsection{Connections and curvature}
 We are only interested in the  algebraic aspects of the classical theory so we 
 place ourselves in the context of differential algebra \cite{kolchin} by considering 
a ring $A$ equipped with and $n$-tuple 
$$(\d^A_1,...,\d^A_n)$$ of commuting derivations. Recall that a derivation is an additive map that satisfies the usual Leibniz rule. For convenience we assume $A$ contains ${\mathbb Q}$. The example we have in mind is, of course, the ring $A$ of smooth functions on ${\mathbb R}^n$ equipped with the partial derivations with respect to the coordinates.  Due to the commutativity requirement for our derivations the setting above is an analogue of the arithmetic situation only in case our number field $F$ is an abelian extension of ${\mathbb Q}$; this was, by the way, the situation considered in \cite{foundations}. 
Following the Introduction to \cite{foundations}
we consider an $n\times n$ matrix of indeterminates $x=(x_{ij})$ and the ring 
$$B=A[x,\det(x)^{-1}].$$
 By a $(\d^A_1,...,\d^A_n)$-{\it connection}  (or simply a {\it connection}) on $GL_n:=Spec\ B$ (or on $B$) we mean an $n$-tuple 
\begin{equation}
\label{one}
(\d_1^B,...,\d_n^B)\end{equation}
 of derivations  on $B$ extending the corresponding derivations $(\d^A_1,...,\d^A_n)$. 
The {\it curvature} of the connection is the family $(\varphi_{ij})$ of commutators
\begin{equation}
\label{cur}
\varphi_{ij}:=[\d^B_i,\d^B_j]=\d_i^B\d_j^B-\d_j^B\d_i^B:B\ra B.\end{equation}
We say that the connection is {\it linear} if 
$$\d^B_i x=A_i x$$ for some $n\times n$ matrices $$A_i=(A_{ijk})$$ with coefficients in $A$. 
For a linear connection the curvature satisfies
 $$\varphi_{ij}(x)=F_{ij}x$$ where $F_{ij}$ is the matrix given by the classical formula
$$
F_{ij}:=\d_i^AA_j-\d_j^AA_i-[A_i,A_j];$$
we still refer to $(F_{ij})$ as the curvature of the connection.
There is one distinguished connection $(\d_{01}^B,...,\d^B_{0n})$ called {\it trivial}, defined by
\begin{equation}
\label{t}
\d_{0i}^Bx=0.\end{equation}

\subsection{Transversal Levi-Civita}
By a {\it metric} we understand a symmetric matrix 
$$q=(q_{ij})\in GL_n(A),\ \ \ q^t=q.$$ 
We define  the {\it Christoffel symbols of the first kind} of the connection $(\d^B_1,...,\d^B_n)$ with respect to the metric $q$ by
\begin{equation}
\label{cris}
\Gamma_{ijk} :=(-A^t_iq)_{jk},\end{equation}
the $(j,k)$-entry of the matrix $-A_i^tq$.
 Passing from the $A_{ijk}$'s to the $\Gamma_{ijk}$'s (and later passing from the entries of the curvature matrices $F_{ij}$ to the components of the covariant Riemann tensor $R_{ijkl}$) is accounted for by our starting with a connection that is {\it dual} to the classical Levi-Civita connection; we adopted this approach simply in order to match the conventions in \cite{foundations}.

Consider the unique $A$-algebra homomorphism 
$$\cH_q:B\ra B$$ such that $$\cH_q(x)=x^tqx.$$
 Say that a $(\d^A_1,...,\d^A_n)$-connection $(\d^B_1,...,\d^B_n)$ is {\it metric} with respect to $q$ if the following diagrams are commutative:
\begin{equation}
 \label{gott}
 \begin{array}{rcl}
 B & \stackrel{\d_i^B}{\longrightarrow} & B\\
 \cH_q \downarrow &\ &\downarrow \cH_q\\
 B & \stackrel{\d_{0i}^B}{\longrightarrow} & B\end{array}
\end{equation}
It is trivial to check that a linear connection is metric with respect to $q$, in the sense of the above (somewhat non-conventional) definition, if and only
if the following classical equalities hold:
\begin{equation}
\label{somewhatnon}\d_i q_{jk}=\Gamma_{ijk} + \Gamma_{ikj}.\end{equation}
Note, by the way, that \ref{somewhatnon} implies the following formula
\begin{equation}
\label{someyes}
\text{tr}(A_i)=- \frac{1}{2}\text{tr}(q^{-1}\d_i q).
\end{equation}
Say that a connection is {\it torsion free} if the following diagrams of $A$-algebras are commutative:
\begin{equation}
\label{risju}
\begin{array}{rcl}
\cO(G) & \stackrel{s_i}{\longleftarrow} & \cO(\mathfrak g)\\
s_j \uparrow & \ & \uparrow r_j\\
\cO({\mathfrak g}) & \stackrel{r_i}{\longleftarrow} & \cO({\mathbb A}^n)
\end{array}
\end{equation}
$${\mathfrak g}:=Spec\ A[x],\ \ {\mathbb A}^n:=Spec\ A[z_1,...,z_n],$$
$$r_i(z_k):=x_{ki},\ \ s_i(x):=-(\d_i x \cdot x^{-1})^t.$$
The commutativity of \ref{risju} is analogous to the commutativity of \ref{risj} and has an invariant meaning involving the 
Lie algebra  of $G$ and Kolchin's logarithmic derivative \cite{kolchin}
(alternatively, the Maurer-Cartan connection); we will not review this interpretation here. 
There is a  minus sign and a transpose in \ref{risju} that do not appear
in \ref{risj}; the discrepancy comes again from the fact that the two situations are ``dual" to each other.
It is trivial to see that a linear connection is torsion free if and only if
the following classical symmetry holds:
\begin{equation}
\label{tf}
\Gamma_{ijk}=\Gamma_{jik}.
\end{equation}

The ``Fundamental Theorem of Riemannian Geometry" is the following statement that can be checked by easy  algebraic manipulations:

  \begin{theorem}
\label{fundfund}
Let $q$ be a metric and $\d_1^A,...,\d^A_n$ commuting derivations on $A$. Then there is a unique linear 
$(\d^A_1,...,\d^A_n)$-connection $(\d^B_1,...,\d^B_n)$
 which is metric with respect to $q$ and torsion free.
It is given by the following formulae:
\begin{equation}\label{windy}
\Gamma_{ijk}=\frac{1}{2}\left(\d_i^A q_{jk}+\d_j^Aq_{ki}-\d_k^Aq_{ij}\right).
\end{equation}\end{theorem}

We refer to $(\d^B_1,...,\d^B_n)$ as the {\it transversal Levi-Civita connection} attached to $(\d^A_1,...,\d^A_n)$ and $q$. This is the ``standard" notion of Levi-Civita connection in classical differential geometry.
At this point it is not clear why we are using the term {\it transversal} for it;
the implication is, of course, that this connection is an analogue of the transversal Levi-Civita connection introduced in our arithmetic theory. 
Will will see that this is the case presently.

For $(F_{ij})$ the curvature 
of the transversal Levi-Civita connection  
 we set:
\begin{equation}
\label{notation monkey}
F_{ij}=(F_{ijkl}),\ \ R_{lij}^k  :=  - F_{ijkl},\ \ \ 
R_{ijkl} := q_{im} R^m_{jkl},\end{equation}
where the repeated index $m$ is summed over.
One refers to $R_{ijkl}$ as the {\it covariant Riemann tensor}; then one shows by easy  algebraic manipulations that:

\begin{proposition}
 The covariant Riemann tensor 
has the following symmetries:
\begin{equation}
\label{macouttt}
\begin{array}{rcl}
R_{ijkl} & =  &- R_{ijlk},\\
\ & \ &  \ \\
R_{ijkl} & = & - R_{jikl},\\
\ & \ &  \ \\
R_{lijk}+R_{ljki}+R_{lkij} & = & 0,\\
\ & \ &  \ \\
R_{ijkl} &= & R_{klij}.
\end{array}
\end{equation}\end{proposition}

In particular if one defines the {\it Ricci tensor}  by the formula
\begin{equation}
R_{ik}:=R^j_{ijk}=q^{jl}R_{jilk},
\end{equation}
where the repeated indeces $j,l$ are summed over, then one gets
the following formal consequence of \ref{notation monkey} and \ref{macouttt}:
\begin{equation}
\label{teago}
R_{jk}=R_{kj}.
\end{equation}

The next Proposition is  a version of a classical formula that appears when one considers {\it normal coordinates}; its proof is, again, a trivial algebraic manipulation.

\begin{proposition} Assume that $J$ is an ideal in $A$ and we are given a metric $q=q^t$ such that $q\equiv 1$ mod $J^2$
where $1$ is, as usual,  the identity matrix.
Then 
 the covariant Riemann tensor satisfies the following congruences:
\begin{equation}
\label{traintrain}
R_{ijkl}\equiv \frac{1}{2}(\d_j \d_k q_{il}+
\d_i \d_l q_{jk} - \d_i \d_k q_{jl} - \d_j \d_l q_{ik})\ \ \ \text{mod}\ \ J.
\end{equation}\end{proposition}

We end our discussion of the classical Levi-Civita connection by recording some classical formulae for the case $n=2$
which are the classical analogues of our formulae \ref{fidx} and \ref{fidxx}
and of Proposition \ref{tirg}.

Indeed assume $n=2$ and $q=d\cdot 1_2$ is a scalar matrix, with $d\in A^{\times}$; this is the case of ``conformal coordinates".  Then it is trivial to check that the transversal Levi-Civita connection attached to $(\d_1^A,\d^A_2)$ and $q$ is defined by
$\d_1^Bx=A_1 x$, $\d_2^Bx=A_2 x$, with
\begin{equation}\label{gogu}
A_1 =-\frac{1}{2}\left(\begin{array}{rr}
\frac{\d_1^A d}{d} & \frac{\d^A_2 d}{d}\\
\ & \ \\
- \frac{\d^A_2 d}{d} & \frac{\d^A_1 d}{d}
\end{array}\right),\ \ \ \ \ A_2 =-\frac{1}{2}\left(\begin{array}{rr}
\frac{\d_2^A d}{d} & - \frac{\d^A_1 d}{d}\\
\ & \ \\
 \frac{\d^A_1 d}{d} & \frac{\d^A_2 d}{d}
\end{array}\right).
\end{equation}
If one considers the algebraic group
$$G':=Spec\ B',\ \ \ B':=A[\alpha,\beta,(\alpha^2+\beta^2)^{-1}]$$
embedded in $G=GL_2=Spec\ B$ via $x\mapsto \left(\begin{array}{cc}
\alpha & \beta\\ -\beta & \alpha\end{array}\right)$, 
then $\d_1^B$ and $\d^B_2$ induce derivations $\d^{B'}_1$ and $\d^{B'}_2$ on $B'$. Furthermore if one considers the algebraic group 
$G'''=Spec\ A[z,z^{-1}]$
and the homomorphism
$\det:G'\ra G'''$, $z\mapsto \alpha^2+\beta^2$,
one trivially checks that 
\begin{equation}
\label{gogu1}
\d_i^{B'}(\alpha^2+\beta^2)=
-\frac{\d_i^A d}{d}\cdot (\alpha^2+\beta^2),\end{equation}
hence $\d_i^{B'}$ induce derivations on $A[z,z^{-1}]$ compatible with $\det$, which are trivially seen to commute on $A[z,z^{-1}]$.

Similarly if one considers the homomorphism
$$\text{det}^{\perp}:G'\ra G''',\ \ \ z\mapsto s:=\frac{\alpha+\sqrt{-1}\beta}{\alpha-\sqrt{-1}\beta},$$
 (defined for $\sqrt{-1}\in A$) then one trivially checks that
 \begin{equation}
 \label{gogu2}\d_1^{B'} s=-\sqrt{-1}\cdot \frac{\d_2^A d}{d} \cdot s,
 \ \ \ \d_2^{B'} s=\sqrt{-1}\cdot \frac{\d_1^A d}{d} \cdot s;\end{equation}
 hence $\d_i^{B'}$ induce derivations $\d_i^{B'''}$ on $B''':=A[z,z^{-1}]$ compatible with $\text{det}^{\perp}$. The derivations $\d_i^{B'''}$ do not commute on $B'''$ in general; indeed we have the following classical formula involving the ``Laplacian of the logarithm":
 \begin{equation}
 \label{gogugogu}
(\d^{B'''}_1\d_2^{B'''}-\d^{B'''}_2\d_1^{B'''})(z)=\sqrt{-1}\cdot (\Delta \log d) \cdot z, \end{equation}
where
$$\Delta \log d:=((\d^A_1)^2+(\d^A_2)^2)\log d:=
\left(\d^A_1\left(\frac{\d_1^A d}{d}\right)+ \d^A_2\left(\frac{\d_2^A d}{d}\right)\right).$$

As explained in previous sections a number of formulae in  the classical setting, especially \ref{someyes}, \ref{windy}, \ref{macouttt}, \ref{traintrain},  \ref{gogu}, \ref{gogu1}, \ref{gogu2} have corresponding arithmetic analogues.

\subsection{Vertical and mixed Levi-Civita}
In what follows we will introduce, in the classical differential geometric setting discussed here, a couple of non-conventional concepts that  we shall call {\it vertical and mixed Levi-Civita connections};  they can be viewed as blueprints of our 
vertical and mixed Levi-Civita connections in the arithmetic case. 
With these concepts at hand it is easier to see why our transversal
Levi-Civita connection in the arithmetic case can be viewed as an analogue
of the transversal Levi-Civita connection in the classical differential geometric case.

Let us start, again, with a ring $A$ equipped, this time, with a single derivation $\check{\d}^A$ and consider
the $n$-tuple of derivations
$$(\check{\d}^A,...,\check{\d}^A).$$
Also consider symmetric matrices 
$$\check{q}_1,...,\check{q}_n\in GL_n(A),\ \ \ \check{q}_i^t=\check{q}_i.$$
For any $(\check{\d}^A,...,\check{\d}^A)$-connection $(\check{\d}^B_1,...,\check{\d}^B_n)$ we may consider 
the {\it Christoffel symbols of the first kind} with respect to $(\check{q}_1,...,\check{q}_n)$ defined by the formulae
\begin{equation}
\label{cris}
\check{\Gamma}_{ijk} :=(-\check{A}^t_iq_i)_{jk},\end{equation}
where
$$\check{\d}^B_i x=\check{A}_i x.$$
We say that a $(\check{\d}^A,...,\check{\d}^A)$-connection $(\check{\d}^B_1,...,\check{\d}^B_n)$ is {\it metric} with respect to $(\check{q}_1,...,\check{q}_n)$  if the following diagrams are commutative:
\begin{equation}
 \begin{array}{rcl}
 B & \stackrel{\check{\d}_i^B}{\longrightarrow} & B\\
 \cH_{q_i} \downarrow &\ &\downarrow \cH_{q_i}\\
 B & \stackrel{\check{\d}_{0}^B}{\longrightarrow} & B\end{array}
\end{equation}
where $\check{\d}_0^B$ is $\check{\d}^A$ on $A$ and $\check{\d}_0^Bx=0$. 

We say that 
$(\check{\d}^B_1,...,\check{\d}^B_n)$ is {\it torsion free} with respect to $(\check{q}_1,...,\check{q}_n)$ if
\begin{equation}
\check{\Gamma}_{ijk}=\check{\Gamma}_{jik}.
\end{equation}
\bigskip

The following is trivial to check:

  \begin{theorem}\label{fundulet}
Assume we are given  an $n$-tuple of metrics $(\check{q}_1,...,\check{q}_n)$, $\check{q}_i=(\check{q}_{ijk})$,  and a derivation $\check{\d}^A$ on $A$. Then
there is a unique linear 
$(\check{\d}^A,...,\check{\d}^A)$-connection $(\check{\d}^B_1,...,\check{\d}^B_n)$
 which is metric and torsion free with respect to $(\check{q}_1,...,\check{q}_n)$.
It is given by the following formulae:
\begin{equation}\label{windy}
\check{\Gamma}_{ijk}=\frac{1}{2}\left(\check{\d}^A \check{q}_{ijk}+\check{\d}^A\check{q}_{jki}-\check{\d}^A\check{q}_{kij}\right).
\end{equation}\end{theorem}

Let us refer to $(\check{\d}^B_1,...,\check{\d}^B_n)$ as the {\it vertical Levi-Civita connection} attached to $\check{\d}^A$ and $(\check{q}_1,...,\check{q}_n)$.

A link  between Theorems \ref{fundfund} and \ref{fundulet} can be established as follows.
Assume one is given commuting derivations 
$\d_1^A,...,\d_n^A$ on $A$, a derivation $\check{\d}^A$ on $A$ commuting with all $\d_i^A$'s,  and  a metric $q\in GL_n(A)$, $q^t=q$. Assume moreover
that one can find an $n\times n$ symmetric matrix 
with coefficients in $A$, which we abusively denote by $(\check{\d}^A)^{-1} q$, such that
\begin{equation}
\check{\d}^A(\check{\d}^A)^{-1}q=q
\end{equation} and set
$$\check{q}_i:=\d_i^A (\check{\d}^A)^{-1} q.$$
In particular
\begin{equation}
\check{\d}^A \check{q}_i=\d_i^Aq,\ \ \ \check{q}_i^t=\check{q}_i.
\end{equation}
Assume in addition one can choose $(\check{\d}^A)^{-1}q$ such that
$$\check{q}_i\in GL_n(A).$$
(Such a matrix $(\check{\d}^A)^{-1}q$ can, of course, be found under very general conditions in the context of smooth functions.)
Consider the Christoffel symbols $\Gamma_{ijk}$ of the transversal Levi-Civita connection $(\d_1^B,...,\d_n^B)$ attached to $(\d_1^A,...,\d_n^A)$ and $q$; cf.  Theorem \ref{fundfund}. Also consider the Christoffel 
symbols $\check{\Gamma}_{ijk}$ of 
the vertical Levi-Civita connection
$(\check{\d}_1^B,...,\check{\d}_n^B)$ attached to $\check{\d}$ and $(\check{q}_1,...,\check{q}_n)$; cf. Theorem \ref{fundulet}. Then, clearly,
$$\check{\Gamma}_{ijk}={\Gamma}_{ijk};$$
in other words, for
$$\d_i^B x=A_i x,\ \ \ \check{\d}_i^Bx=\check{A}_i x$$
 we have equalities of matrices:
\begin{equation}
\check{q}_i \check{A}_i=q A_i.
\end{equation}
The above construction can be considered, of course, in the special case when $\check{\d}^A$ is the derivation
$$\check{\d}^A:=\d_k^A,$$
 where $k$ is any of the indices $1,...,n$. In this case write
 $$ \check{\d}_i^B=
\check{\d}_{ki}^B,
 \ \ \ \ (\check{\d}^A)^{-1}q=(\d_k^A)^{-1}q,\ \ \ \ \ \check{A}_i=\check{A}_{ki},\ \ \ 
\check{q}_i=\check{q}_{ki},$$
 so we have that 
 \begin{equation}
 \label{pies}
 \d_k^A (\d^A_k)^{-1} q=q,\ \ \ \ 
 \check{q}_{ki}=\d_i^A(\d^A_k)^{-1}q,
 \ \ \ \check{q}_{ki}  \check{A}_{ki}=qA_i,\end{equation}
 and
 $$(\check{\d}_{ki}^B)_{|A}=\d_k^A,\ \ \ \check{\d}_{ki}^Bx=\check{A}_{ki}x.$$
 The family 
 \begin{equation}
 \label{heehee}
 (\check{\d}_{ji}^B)\end{equation}
 indexed by $i,j=1,...,n$ can be referred to as the {\it mixed Levi-Civita connection}. 
Note that 
\begin{equation}
\label{mordecai}
\check{q}_{ii}=q,\ \ \ \check{A}_{ii}=A_i,\ \ \  \check{\d}^B_{ii}=\d^B_i.\end{equation}
In other words the transversal Levi-Civita connection $(\d_i^B)$ attached to $q$  can be extracted from the 
mixed Levi-Civita connection \ref{heehee}
 by ``taking the diagonal",
 \begin{equation}
 \label{diagonall}
 (\check{\d}^B_{ii}).\end{equation}
 This ``taking the diagonal" procedure is the analogue of ``using a transversal gauge"
 in our Definition \ref{papu}; 
  our (arithmetic) vertical,  mixed, and transversal Levi-Civita connections are therefore the analogues of the  vertical, mixed, and transversal Levi-Civita connections that we discussed in this Appendix;
  cf. also  Definition \ref{papu}. Note that there is a discrepancy between taking the diagonal, i.e., indices $ii$ 
  in \ref{diagonall} and setting the upper index in \ref{asa} equal to $1$ instead of equal to $i$; this discrepancy should be viewed as  a mere artifact of our normalizations.  Indeed the choice of the indices $ii$ in \ref{diagonall} corresponds to the fact that $\check{q}_{ii}=q$ in \ref{mordecai} while the choice 
  of the upper index $1$ in \ref{asa} corresponds to the fact that $\sigma_1^{-1}q=q$ in \ref{secret}. So the two choices in these two contexts are analogous to each other.
     
  Note also that the  curvature 
 $$\varphi_{ij}:=[\d_i^B,\d_j^B],\ \ \ \varphi_{ij}(x)=F_{ij}x,$$
 of the transversal Levi-Civita connection can be read off the collection of commutators
 $$\check{\varphi}_{klij}:=[\check{\d}_{ki}^B,\check{\d}_{lj}^B]$$
 between the derivations appearing in   the mixed Levi-Civita connection:
 $$\varphi_{ij}=\check{\varphi}_{ijij}.$$
 The  collection $\check{\varphi}_{klij}$ has, as analogue in our arithmetic theory, the mixed curvature; cf Definition \ref{woken}.

\end{document}